\documentclass[a4paper,12pt]{amsart}

\newcommand{\supp}{\text {\rm supp}}

\newcommand{\sgn}{\text{ \rm sgn}}

\def\<{\langle} 
\def\>{\rangle}

\def\Om{\Omega}

\def\ge{\geqslant}
\def\le{\leqslant}
\def\a{\alpha}
\def\b{\beta}
\def\g{\gamma}

\def\d{\delta}

\def\e{\epsilon}

\def\o{\omega}
\def\p{\pi}

\def\s{\sigma}
\def\t{\tau}
\def\th{\theta}

\def\l{\lambda}

\def\i{^{-1}}
\def\ZZ{\mathbb Z}
\def\NN{\mathbb N}
\def\QQ{\mathbb Q}

\def\cd{\mathcal D}

\def\co{\mathcal O}
\def\cp{\mathcal P}

\def\cw{\mathcal W}

\def\aa{\mathbf a}
\def\ua{\underline {\mathbf a}}

\def\tS{\tilde S}

\def\tx{\tilde x}
\def\ty{\tilde y}

\def\tW{\tilde W}
\def\tw{\tilde w}

\def\lto{\hookrightarrow}

\usepackage{amssymb} 
\usepackage{latexsym} 
\usepackage{amsfonts} 
\usepackage{amsmath} 
\usepackage{eucal} 
\usepackage{bm} 
\usepackage{bbm} 
\usepackage{graphicx} 
\usepackage[english]{varioref} 
\usepackage[nice]{nicefrac} 
\usepackage[all]{xy}

\newcommand{\kk}{\Bbbk}
\newcommand{\PP}{\mathbbm{P}}

\newcommand{\Q}{\mathbbm{Q}}

\usepackage{amsthm}

\theoremstyle{plain}
\newtheorem{thm}{Theorem}[section] 
\newtheorem*{thm*}{Theorem} 
 \newtheorem{prop}[thm]{Proposition}
 \newtheorem{lem}[thm]{Lemma}
 \newtheorem{cor}[thm]{Corollary}

\theoremstyle{definition}

\theoremstyle{remark}

\newtheorem*{rmk}{Remark}
\newtheorem*{claim*}{Claim}

\begin{document}

\author{Xuhua He}
\address{Department of Mathematics, The Hong Kong University of Science and Technology, Clear Water Bay, Kowloon, Hong Kong}
\email{maxhhe@ust.hk}
\thanks{The author is partially supported by HKRGC grants 601409.}
\title[]{Minimal length elements in conjugacy classes of extended affine Weyl group}

\begin{abstract}
We study the minimal length elements in an integral conjugacy class of a classical extended affine Weyl group and we show that these elements are quite ``special'' in the sense of Geck and Pfeiffer \cite{GP93}. We also discuss some application on extended affine Hecke algebras and loop groups. 
\end{abstract}

\maketitle

\section*{Introduction}

\subsection{}\label{start} The main purpose of this paper is to investigate some ``special'' elements in an extended affine Weyl groups and to discuss some applications on extended affine Hecke algebras and loop groups. 

More precisely, let $\tW$ be an extended affine Weyl group and $W_a \subset \tW$ be the affine Weyl group. Let $\tW_{m}$ be the subset of elements that are of minimal length in their conjugacy classes of $\tW$. We call an element $\tw \in \tW$ a good element if $l(\tw^n)=n l(\tw)$ for all $n \in \NN$. Let $\tW_{good}$ be the subset of good elements in $\tW$. Then $$\tW_{good} \subset \tW_m \subset \tW.$$ 

We will prove the following results:

(1) Let $\tW$ be of classical type. Then any element in $W_a$ can be reduced to an element in $\tW_m$ and any two minimal length element in the same conjugacy class are ``strongly conjugate'' in the sense of $\S$\ref{to}. 

(2) If $\tW$ is of type A or C, then any element $\tw \in W_a$ can be ``reduced'' further to a good element. 

\subsection{} In this subsection, we make a short digression to finite Weyl groups, which provide some motivation for our study on the affine case. 

Minimal length element in a conjugacy class of a finite Weyl group was first studied by Geck and Pfeiffer in \cite{GP93}. They showed that these elements are ``special'' in the sense of $\S$\ref{start}(1). Later in \cite{GKP}, Geck, Kim and Pfeiffer generalized the result to ``twisted'' conjugacy classes in finite Coxeter groups. In \cite{He072}, we gave a new approach to study the minimal length element using partial conjugation action. 

These minimal length elements have many interesting applications in representation theory. 

Geck and Pfeiffer in \cite{GP93} used the minimal length elements to define ``character tables'' for Iwahori-Hecke algebras. The ``character tables'' were later defined by Geck, Kim and Pfeiffer in \cite{GKP} for the twisted Iwahori-Hecke algebras. The knowledge of character tables can be used to study modular representations of Iwahori-Hecke algebras (see Geck-Pfeiffer \cite{GP93}) and to determine the L-functions of Deligne-Lusztig varieties (see Lusztig \cite{Lu84a} and Digne-Michel \cite{DM85}).

The minimal length elements were also used in the study of unipotent representations and character sheaves (see Lusztig \cite{Lu02}), in the study of $G$-stable pieces and parabolic character sheaves (see \cite{He072}), in the proof of Orlik-Rapoport conjecture \cite{OR08} on the affineness of certain Deligne-Lusztig varieties (see \cite{He082} and Bonnaf\'e-Rouquier \cite{BR08}) and in the study of the finer Deligne-Lusztig varieties (see \cite{He09}). 

\subsection{} We come back to the affine case and discuss some basic ideas for the proof of our main result. 

Our aim is to find minimal length elements for a conjugacy class of $\tW$ and to show that any other elements in that conjugacy class can be ``reduced'' to a minimal one. 

A big difference between finite Weyl groups and extended affine Weyl groups is that a conjugacy class in an extended affine Weyl group has, in general, infinitely many elements. For example, finding a representative of minimal length elements for a given conjugacy class in the extended affine Weyl group of type $\tilde A_2$, seems not an easy task. 

To overcome the difficulty, we combine the method of ``partial conjugation action'' developed in \cite{He072} and the $P$-operators introduced in section 4. Roughly speaking, each conjugacy class is a disjoint union of $W$-conjugacy classes, where $W$ is the finite Weyl group. The ``partial conjugation action'' allows us to find a minimal length element in each $W$-conjugacy class and the $P$-operator changes a $W$-conjugacy class to another one. By combining the $P$-operator and the ``partial conjugation action'', we are able to construct a list of representatives of elements of minimal length in a given conjugacy class and to show that any other elements in this conjugacy class can reach one of the representatives by a sequence of conjugation by simple reflections which weakly decrease the length. This proves $\S$\ref{start}(1). See Theorem \ref{61} and  \ref{thA}.

Based on these explicit representatives, any element can be reduced further to a minimal length element in a ``distinguished'' conjugacy class in the sense of Theorem \ref{62}. Moreover, as we will see in Corollary \ref{xx}, in type A and C, any minimal length element in a distinguished conjugacy class is a good element. This proves $\S$\ref{start}(2). 

\subsection{}\label{sa} Now we discuss some application on extended affine Hecke algebra and loop groups. 

Based on $\S$\ref{start}(1), we introduce the class polynomials for classical extended affine Hecke algebras, generalizing the construction of Geck and Pfeiffer for finite Hecke algebras in \cite{GP93}. 

Before discussing application on loop groups, we need to introduce some more notations. 

Let $G$ be a connected reductive algebraic group over an algebraically closed field $\kk$. Let $L=\kk((\e))$ is the field of formal Laurent power series and $o=\kk[[\e]]$ be the ring of formal power series. 

We consider a ``twisted'' conjugation action of $G(L)$ on itself as $g \cdot_\d h=g h \d(g) \i$ for $g, h \in G(L)$. Here $\d$ is a bijective group homomorphism on $G(L)$ of one of the following type:

(1) For any nonzero element $a \in \kk$, define $\d_a(p(\e))=p(a \cdot \e)$ for any formal Laurent power series $p(\e)$. We extend $\d_a$ to a group homomorphism on $G(L)$, which we still denote by $\d_a$.

(2) If $\kk$ is of positive characteristic and $F: \kk \to \kk$ be a Frobenius morphism. Then set $F(\sum a_n \e^n)=\sum F(a_n) \e^n$. We extend $F$ to a group homomorphism on $G(L)$, which we denote by $\d_F$. 

The $\d_a$-conjugacy classes are studies by Baranovsky and Ginzburg in \cite{BG96}. The $\d_F$-conjugacy classes are studied by Kottwitz in \cite{Ko97}. 

\subsection{} Let $I$ be a Iwahori subgroup of the loop group $G(L)$. The quotient $G(L)/I$ is called an affine flag variety. Let $\tW$ be the extended affine Weyl group of $G(L)$. For $b \in G(L)$ and $\tw \in \tW$, we set $$X_{\tw, \d}(b)=\{g I \in G(L)/I; g \i b \d(g) \in I \dot \tw I\}.$$ Then we have a partition of affine flag variety $G(L)/I=\sqcup_{\tw \in \tW} X_{\tw, \d}(b)$. In the case that $\d=\d_F$, $X_{\tw, \d}(b)$ is called an {\it affine Deligne-Lusztig variety}. In the case that $\d$ is the identity map, $X_{1, \d}(b)$ is called an {\it affine Springer fiber}. $X_{\tw, \d}(b)$ is also considered in \cite[Section 7]{L3}. 


Let $J_{b, \d}=\{g \in G(L); g \i b \d(g)=b\}$ be the centralizer of $b$ for the twisted conjugation action. Then $J_{b, \d}$ acts on $X_{\tw, \d}(b)$ and on the homology of $X_{\tw, \d}(b)$ in the natural way. 

\subsection{} For simplicity, we only mention here some applications for loop group of type A and C. Some weaker results are also obtained for other classical types. 

We consider $G(L) \cdot_\d I \dot \tw I$ and $X_{\tw, \d}(b)$ for $\tw \in W_a$. By some reduction method discussed in Section 9, these objects are ``constructed'' from the objects corresponding to $\tw \in \tW_m$. On the other hand, the objects corresponding to a good element can be fairly understood. Now based on $\S$\ref{start}(1) \& (2), we have that 

\

(1) A stratification of $G(L)$ into locally closed subschemes that are equivariant for the $\d$-conjugation action $$G(L)=\sqcup_{[\tw] \in \tW_{good}/\asymp} G(L) \cdot_\s I \dot \tw I.$$

If $\d=\d_F$, this follows from Kottwitz's classification of $\s_F$-conjugacy classes. In the case that $\d=\d_a$, this is proved in Prop \ref{str} for $GL_n(L)$ and the identity component of $PSP_{2n}(L)$. We expect that such stratification holds for other types. We also believe that this stratification is a necessary ingredient for establishing the (conjectural) theory of character sheaves of loop groups for the $\d_a$-conjugation action, where $a$ is not a root of unity. 

\

(2) A formula on the dimension of $X_{\tw, \d}(b)$. 

(i) If $\d=\d_F$, then the dimension of affine Deligne-Lusztig variety $X_{\tw, \d}(b)$ can be expressed in terms of the degree of certain class polynomials of the extended affine Hecke algebra. See Proposition \ref{class3}(1). 

(ii) If $\d=\d_a$ with $a$ not a root of unity, then $X_{\tw, \d}(b)$ is always finite dimensional. See Proposition \ref{class3}(2). 

(iii) If $\d$ is the identity map and $b$ is a regular semisimple integral element, then $X_{\tw, \d}(b)$ is finite dimensional. See Proposition \ref{class5}. This answers a question of Lusztig for $GL_n$ and $PSP_{2n}$. 

\

(3) Some reduction on the homology of affine Deligne-Lusztig varieties. 

Given $b \in G(L)$ and an open compact subgroup $K$ of $J_{b, \d_F}$, there are infinitely many $\tw \in \tW$ with $X_{\tw, \d_F}(b) \neq \emptyset$. Then a priori, to get all the simple $K$-modules that are obtained from homological construction, one needs to calculate the homology of $X_{\tw, \d_F}(b)$ for any $\tw \in \tW$ with $X_{\tw, \d_F}(b) \neq \emptyset$. 

However, we show in Theorem \ref{BM} that for type A and C, one only needs to calculate the homology of finitely many affine Deligne-Lusztig varieties. Moreover, if $b=\e^\chi$ for some regular coweight in the coroot lattice or $b$ is a superbasic element in $GL_n(L)$, then one needs to calculate the homology of only one affine Deligne-Lusztig varieties and $J_{b, \d_F} \cap I$ acts trivially on the homology of $X_{\tw, \d_F}(b)$ for any $\tw \in \tW$. See Corollary \ref{ddd} and \ref{superbasic}. 

\subsection{} We now review the content of this paper in more detail. 

The notations will be introduced in section 1. In section 2, we introduce the groups $\tW^!$ of classical type and give a parametrization of their integral conjugacy classes in terms of pairs of double partitions. These groups are closely related to the classical extended affine Weyl groups and we consider $\tW^!$ instead of $\tW$ mainly for some technical reason. In section 3, we recall the result of ``partial conjugation action'' in \cite{He072}, which leads to the definition of $\cw$, a key ingredient in the proof of our main theorems. We then give an explicit description of $\cw$ for $\tW^!$ of classical types. In section 4, we introduce the $P$-operators, another key ingredient in the proof of our main theorems. Section 5 is the most technical part of this paper, in which we establish some fundamental properties of the $P$-operators. In section 6, we provide explicit representatives of minimal length elements in each integral conjugacy class of $\tW^!$ and prove $\S$\ref{start}(1). In section 7, we introduce distinguished conjugacy classes and prove that any element in $W_a$ can be ``reduced'' further to a minimal length element in a distinguished conjugacy class. We also introduce some partial order on the set of distinguished conjugacy classes. In section 8, we introduce the good elements for extended affine Weyl groups and prove their existence. We also show that good elements are exactly minimal length elements in distinguished conjugacy classes for type A and C. 

The rest of this paper is the applications of the main theorems. In section 9, we discuss some reduction method in the study of affine Hecke algebra and affine flag varieties. We introduce the class polynomials for classical extended affine Hecke algebras, generalizing \cite{GP93} for finite Hecke algebras. We also give a formula which relates the dimension of affine Deligne-Lusztig varieties with the degree of class polynomials. In section 10, we discuss the dimensions of some subvarieties of affine flag varieties corresponding to a good element. In section 11, we focus on the loop groups of type A and C. We first give a stratification of the loop group into locally closed subvarieties, stable under the $\d$-conjugation action. We then give a dimension formula for the affine Deligne-Lusztig varieties, sharpening the formula in section 9. In the end, we discuss the homology of affine Deligne-Lusztig varieties and representation of $J_{b, \s}$ obtained from homological construction. 

\section{Notations}

\subsection{} Let $B$ be a Borel subgroup of $G$ and $B^-$ be an opposite Borel subgroup. Let $T=B \cap B^-$ be a maximal torus of $G$. Let $W=N_G(T)/T$ be the finite Weyl group of $G$. For $w \in W$, we choose a representative $\dot w \in N_G(T)$. Let $P^\vee$ be the coweight lattice and $Q^\vee$ be the coroot lattice of $G$. 

Let $\tW=P^\vee \rtimes W=\{t^\chi w; \chi \in P^\vee, w \in W\}$ be the extended affine Weyl group of the loop $G(L)$. The multiplication is given by the formula $(t^\chi w) (t^{\chi'} w')=t^{\chi+w \chi'} w w'$. For $\tw=t^\chi w \in \tW$, we choose a representative $\dot \tw=\e^\chi \dot w$ in $G(L)$. 

Let $I$ be the inverse image of $B^-$ under the projection map $G(o) \mapsto G$ sending $t$ to $0$. The we have the Bruhat-Tits decomposition $G(L)=\sqcup_{\tw \in \tW} I \dot{\tw} I$. If $\t \in \tW$ with $l(\t)=0$, then $\dot \t I \dot \t \i=I$. 

\subsection{} Let $\Phi$ be the set of roots of $G$ and $\Phi^+$ (resp. $\Phi^-$) be the set of positive (resp. negative) roots of $G$. Let $(\a_i)_{i \in S}$ be the set of simple roots of $G$. For any $i \in S$, let $s_i$ be the corresponding simple reflection in $W$ and $\o_i$ be the corresponding fundamental coweight. Set $\tS=S \cup \{0\}$ and $s_0=t^{\th^\vee} s_\th$, where $\th$ is the largest positive root of $G$. 

For any $\a \in \Phi$, let $u_\a: \kk \to G$ with $t u_\a(k) t\i=u_\a(\a(t) k)$ for $k \in \kk$. For any $J \subset S$, let $\Phi_J$ be the set of roots spanned by $(\a_i)_{i \in J}$. 

Let $W_a=Q^\vee \rtimes W$ be the affine Weyl group. It is a Coxeter group  with generators $s_i$ (for $i \in \tS$) and is a normal subgroup of $\tW$. Following \cite{IM65}, we define the length function on $\tW$ by $$l(t^\chi w)=\sum_{\a \in R^+, w \i(\a) \in R^+} |<\chi, \a>|+\sum_{\a \in R^+, w \i(\a) \in R^-} |<\chi, \a>-1|.$$ For any coset of $W_a$ in $\tW$, there is a unique element of length $0$. Moreover, there is a natural group isomorphism between $\{\t \in \tW; l(\t)=0\}$ and $\tW/W_a \cong P^\vee/Q^\vee$. 

\subsection{} In general, let $\Om$ be a group of automorphisms on $W_a$ that sends $\tS$ to $\tS$. Set $\tW^!=W_a \rtimes \Om$. We regard the elements in $\Om$ as length $0$ elements and extend the length function to $\tW^!$ by $l(w \t)=l(w)$ for $w \in W_a$ and $\t \in \Om$. 

For $\t, \t' \in \Om$ and $w, w' \in W_a$, we say that $w \t \le \t w' \t'$ if $\t=\t'$ and $w \le w'$ for the Bruhat order on $W_a$. 

For any $J \subsetneqq \tS$, let $W_J$ be the subgroup of $W_a$ generated by $s_j$ (for $j \in J$), $w_J$ be the longest element in $W_J$ and ${}^J \tW^!$ be the set of minimal elements for the cosets $W_J \backslash \tW^!$. 

For any subset $C$ of $\tW^!$, set $$C_{\min}=\{w \in C; l(w) \le l(w') \text{ for any } w' \in C\}$$ and $l(C)=l(w)$ for any $w \in C_{\min}$. 

\subsection{}\label{to} 

For $\tw, \tw' \in \tW^!$ and $i \in \tS$, we write $\tw \xrightarrow{s_i} \tw'$ if $\tw'=s_i \tw s_i$ and $l(\tw') \le l(\tw)$. We write $\tw \to \tw'$ if there is a sequence $\tw=\tw_0, \tw_1, \cdots, \tw_n=\tw'$ of elements in $\tW^!$ such that for any $k$, $\tw_{k-1} \xrightarrow{s_i} \tw_k$ for some $i \in \tS$.  

We write $\tw \tilde \to \tw'$ if there is a sequence $\tw=\tw_0, \tw_1, \cdots, \tw_n=\tw'$ of elements in $\tW^!$ such that for any $k$, $\tw_k=a \tw_{k-1} a \i$ for some $a \in \Om$ or $\tw_{k-1} \xrightarrow{s_i} \tw_k$ for some $i \in \tS$. 

We write $\tw \approx \tw'$ if $\tw \to \tw'$ and $\tw' \to \tw$ and $\tw \tilde \approx \tw'$ if $\tw \tilde \to \tw'$ and $\tw' \tilde \to \tw$. 

We call $\tw, \tw' \in \tW^!$ {\it elementarily strongly conjugate} in $\tW^!$ (resp. in $W_a$) if $l(\tw)=l(\tw')$ and there exists $\tx \in \tW^!$ (resp. $\tx \in W_a$) such that $\tw'=\tx \tw \tx \i$ and either $l(\tx \tw)=l(\tx)+l(\tw)$ or $l(\tw \tx \i)=l(\tx)+l(\tw)$. We call $\tw, \tw'$ {\it strongly conjugate} in $\tW^!$ (resp. in $W_a$) if there is a sequence $\tw=\tw_0, \tw_1, \cdots, \tw_n=\tw'$ such that for each $i$, $\tw_{i-1}$ is elementarily strongly conjugate to $\tw_i$ in $\tW^!$ (resp. in $W_a$). We write $\tw \tilde \sim \tw'$ if $\tw$ and $\tw'$ are strongly conjugate in $\tW^!$ and $\tw \sim \tw'$ if $\tw$ and $\tw'$ are strongly conjugate in $W_a$. It is easy to see that if $\tw \approx \tw'$, then $\tw \sim \tw'$ and if $\tw \tilde \approx \tw'$, then $\tw \tilde \sim \tw'$. 

Since conjugation by an element in $\Om$ preserves the length function and in particular permutes $\tS$, we have that 

(1) $\tw \tilde \to \tw'$ if and only if $\tw \to a \tw' a \i$ for some $a \in \Om$;

(2) $\tw \tilde \sim \tw'$ if and only if $\tw \sim a \tw' a \i$ for some $a \in \Om$.


If $\tw' \approx \tw$ for all $\tw' \in \tW^!$ with $\tw \to \tw'$, then we call $\tw$ {\it terminal}. It is easy to see by induction on length that for any $\tx \in \tW^!$, there exists a terminal element $\tx'$ such that $\tx \to \tx'$. 

\subsection{} Let $S'_n=(\ZZ/2 \ZZ)^n \rtimes S_n$ be the set of permutations $\s$ on $\{\pm 1, \cdots, \pm n\}$ with $\s(-i)=-\s(i)$ for all $i$. If $\s \in S'_n$ and there is only one or two orbits on $\{\pm 1, \cdots, \pm n\}$ consisting more than one element and the  orbit(s) are of the form $$i_1 \to i_2 \to \cdots \to i_k \to i_1 \text{ and/or } -i_1 \to -i_2 \to \cdots \to -i_k \to -i_1,$$ then we simply write $(i_1 i_2 \cdots i_k)$ for $\s$. 

A {\it partition} $\l$ is a sequence of positive integers $[a_1, a_2, \cdots, a_k]$ with $a_1 \ge a_2 \ge \cdots \ge a_k$. We write $|\l|$ for $a_1+a_2+\cdots+a_k$. We write $\O$ for the empty partition and set $|\O|=0$. 

For a pair of partitions $(\l, \mu)$ with $\l=[a_1, \cdots, a_k]$ and $\mu=[a_{k+1}, \cdots, a_l]$ and $\sum_{i=1}^l a_i=n$, we set  \begin{align*} w_{(\l, \mu)}=& (|\l|+a_{k+1}, -|\l|-a_{k+1}) (|\l|+a_{k+1}+a_{k+2}, -|\l|-a_{k+1}-a_{k+2}) \cdots (n, -n) \\ & (1 2 \cdots a_1) \cdots (n-a_l+1, n-a_l+2, \cdots, n). \end{align*}

\subsection{}\label{cla} 

A {\it double partition} $\tilde \l$ is a sequence $[(b_1, c_1), \cdots, (b_k, c_k)]$ with $(b_i, c_i) \in \NN \times \ZZ$ for all $i$ and $(b_1, c_1) \ge \cdots \ge (b_k, c_k)$ for the lexicographic order on $\NN \times \ZZ$. We write $|\tilde \l|$ for $(\sum b_i, \sum c_i)$ and $\l$ for $[b_1, \cdots, b_k]$. We write $\tilde \O$ for the empty double partition and set $|\tilde \O|=(0, 0)$. 

We call a double partition $\tilde \l=[(b_1, c_1), \cdots, (b_k, c_k)]$ {\it positive} if $c_i \ge 0$ for all $i$. 

We call a double partition $\tilde \mu=[(b_1, c_1), \cdots, (b_k, c_k)]$ {\it special} if $c_i \in \{0, 1\}$ for all $i$. In this case, let $\underline{\tilde \mu}$ be the double partition whose entries are $(b_1, 1-c_1), \cdots, (b_k, 1-c_k)$. Then $\underline {\tilde \mu}$ is also special. 

Let \begin{gather*} \cd \cp=\{(\tilde \l, \tilde \mu); \tilde \mu \text{ is special and } |\tilde l|+|\tilde \mu|=(n, *)\}, \\ \cd \cp_{\ge 0}=\{(\tilde \l, \tilde \mu) \in \cd \cp; \tilde \l \text{ is  positive}\}. \end{gather*} Let $\sim$ be the equivalent relation on $\cd \cp$ defined by $(\tilde \l, \tilde \mu) \sim (\tilde \l, \underline{\tilde \mu})$ for all $(\tilde \l, \tilde \mu) \in \cd \cp$. 

For $\tilde \l=[(b_1, c_1), \cdots, (b_k, c_k)]$ and $\tilde \mu=[(b_{k+1}, c_{k+1}), \cdots, (b_l, c_l)]$ with $(\tilde \l, \tilde \mu) \in \cd \cp$, we set $$\tw_{(\tilde \l, \tilde \mu)}=t^{[a_1, \cdots, a_n]} w_{(\l, \mu)} \in \ZZ^n \rtimes S'_n,$$ where \[a_i=\begin{cases} c_j, & \text{ if } i=b_1+\cdots+b_j \text{ for some } j; \\ 0, & \text{ otherwise}. \end{cases}\]

For any $(\tilde \l, \tilde \O) \in \cd \cp$, let $[(\tilde \l, \tilde \O)]=(\ZZ^n \rtimes S_n) \cdot \tw_{(\tilde \l, \tilde \O)}$ be the conjugacy class of $\ZZ^n \rtimes S_n$ that contains $\tw_{(\tilde \l, \tilde \O)}$. For any $(\tilde \l, \tilde \mu) \in \cd \cp_{\ge 0}$, let $[\tilde \l, \tilde \mu)]'=(\ZZ^n \rtimes S'_n) \cdot \tw_{(\tilde \l, \tilde \mu)}$ be the conjugacy class of $\ZZ^n \rtimes S'_n$ that contains $\tw_{(\tilde \l, \tilde \mu)}$.

\section{Parametrization of conjugacy classes} 

\subsection{} We first recall the parametrization of conjugacy classes of classical finite Weyl groups.

{\bf Type $A_{n-1}$.} We may regard $W(A_{n-1})=S_n$ as $\{\s \in S'_n; \s(i)>0, \forall i>0\}$. There is a bijection between the conjugacy classes of $W(A_{n-1})$ and the pairs of partitions $(\l, \O)$ with $|\l|=n$ and the element $w_{(\l, \O)}$ is a representative for the conjugacy class corresponding to $(\l, \O)$. 

{\bf Type $B_n$ and $C_n$.} We may regard $W(B_n)=W(C_n)$ as $S'_n$. There is a bijection between the conjugacy classes of $W(B_n)=W(C_n)$ and the pair of partitions $(\l, \mu)$ with $|\l|+|\mu|=n$ and the element $w_{(\l, \mu)}$ is a representative for the conjugacy class corresponding to $(\l, \mu)$. 

{\bf Type $D_n$.} We may regards $W(D_n)$ as $\{\s \in S'_n; \prod_{i=1}^n \s(i)=n!\}$. A conjugacy class in $W(C_n)$ intersects $W(D_n)$ if and only if it corresponds to the set of pair of partitions $(\l, \mu)$ with $2 \mid |\mu|$. In this case, $w_{(\l, \mu)} \in W(D_n)$ and the intersection is a single conjugacy class of $W(D_n)$ and the representative is given by $w_{(\l, \mu)}$ except the case that $\mu$ is the empty partition and all the entries of $\l$ are even numbers. In the last case, there are two conjugacy classes of $W(D_n)$ and the conjugation of $(n, -n) \in W(C_n)-W(D_n)$ permutes these two conjugacy classes.

\subsection{}\label{tww} Let $V=\mathbb R^n$ and $e_1, \cdots, e_n$ be a standard basis of $V$. We identify $V$ with $V^*$ in such a way that $<e_i, e_j>=\d_{i j}$. Set $$\PP=\oplus_{i=1}^n \ZZ e_i=\{[a_1, \cdots, a_n]; a_i \in \ZZ\}.$$ The action of $S'_n$ on $\PP$ is defined by $$\s \cdot [a_1, \cdots, a_n]=[a_{\s \i(1)}, \cdots, a_{\s \i(n)}],$$ here $a_{-i}=-a_i$ for $i \in \{1, \cdots, n\}$. 

{\bf Type A.} The roots are $\{e_i-e_j; i \neq j\}$ and the simple roots are $\a_1=e_1-e_2, \cdots, \a_{n-1}=e_{n-1}-e_n$. The coweight lattice is $P^\vee(A_{n-1})=\PP/\mathbb Z (e_1+e_2+\cdots+e_n)$. Set $\tW^!_A=\PP \rtimes S_n \cong \ZZ^n \rtimes S_n$. The natural projection $\PP \to P^\vee(A_{n-1})$ extends to a surjective group homomorphism $p: \tW^!_A \to \tW(A_{n-1})$. For $\tw \in \tW^!_A$, we define the length $l(\tw)$ of $\tw$ as $l(p(\tw))$. So the element $t^{[a, \cdots, a]} \in \tW^!_A$ if of length $0$ for any $a \in \ZZ$. 

{\bf Type B.} The roots are $\{\pm e_i; 1 \le i \le n\} \sqcup \{\pm e_i \pm e_j; 1 \le i<j \le n\}$ and the simple roots are $\a_1=e_1-e_2, \cdots, \a_{n-1}=e_{n-1}-e_n, \a_n=e_n$. The coweight lattice is $P^\vee(B_n)=\PP$. We set $\tW^!_B=\tW(B_n) \cong \ZZ^n \rtimes S'_n$. 

{\bf Type C.} The roots are $\{\pm 2 e_i; 1 \le i \le n\} \sqcup \{\pm e_i \pm e_j; 1 \le i<j \le n\}$ and the simple roots are $\a_1=e_1-e_2, \cdots, \a_{n-1}=e_{n-1}-e_n, \a_n=2 e_n$. The coweight lattice is $P^\vee(C_n)=\PP \sqcup (\PP+\o)$, where $\o=[\frac{1}{2}, \cdots \frac{1}{2}]$. We set $\tW^!_C=\tW(C_n)$. 

{\bf Type D.} The roots are $\{\pm e_i \pm e_j; 1 \le i<j \le n\}$ and the simple roots are $\a_1=e_1-e_2, \cdots, \a_{n-1}=e_{n-1}-e_n, \a_n=e_{n-1}+e_n$. The coweight lattice is $P^\vee(C_n)=\PP \sqcup (\PP+\o)$, where $\o=[\frac{1}{2}, \cdots \frac{1}{2}]$. Let $\iota=(n, -n)$ be the outer diagram automorphism interchanging $\a_{n-1}$ and $\a_n$. Set $\tW^!_D=\tW(D_n) \rtimes \<\iota\>$ and regard $\iota$ as a length-$0$ element in $\tW^!_D$. 

We say that $\tW^!$ is {\it of classical type} if $\tW^!=\tW^!_*$ for $*$ is A, B, C or D. 

For type BCD, $\ZZ^n \rtimes S'_n$ is a subgroup of $\tW^!$ and equals to $W_a$ for type C and is a union of two cosets of $W_a$ for type B or D. 

Notice that the group structure on $\tW^!_C$ and on $\tW^!_D$ are the same. However, the length functions are different.

\subsection{}\label{psc} An element $\tw \in \tW^!$ is called {\it integral} if $\tw \in \ZZ^n \rtimes S'_n$, i.e., the translation part of $\tw$ lies in $\mathbb P$. We denote by $\tW^!_{int}=\tW^! \cap (\ZZ^n \rtimes S'_n)$ the set of all integral elements in $\tW^!$. It is easy to see that $\tW^!_{int}$ is a union of conjugacy classes. We have that

(1) $\ZZ^n \rtimes S_n=\sqcup_{(\tilde \l, \tilde \O) \in \cd \cp} [(\tilde \l, \tilde \O)]$ is the union of conjugacy classes. 

(2) $\ZZ^n \rtimes S'_n=\sqcup_{(\tilde \l, \tilde \mu) \in \cd \cp_{\ge 0}} [(\tilde \l, \tilde \mu)]'$ is the union of conjugacy classes. 

(3) Let $\o=[\frac{1}{2}, \frac{1}{2}, \cdots, \frac{1}{2}]$. Then for any $(\tilde \l, \tilde \mu) \in \cd \cp_{\ge 0}$, $t^{-\o} [(\tilde \l, \tilde \mu)]' t^\o=[(\tilde \l, \underline{\tilde \mu})]'$. 

The parametrization of integral conjugacy classes of $\tW^!$ follows easily from (1)-(3) above. 

{\bf Type A}. All the elements in $\tW^!_A$ are integral. There is a bijection between the set of (integral) conjugacy classes of $\tW^!_A=\ZZ^n \rtimes S_n$ and the set of pairs of double partitions $(\tilde \l, \tilde \O) \in \cd \cp$. 

{\bf Type B}.  All the elements in $\tW^!_B$ are integral. There is a bijection between the set of (integral) conjugacy classes of $\tW^!_B=\ZZ^n \rtimes S'_n$ and $\cd \cp_{\ge 0}$.

{\bf Type C}.  Here $\tW^!_{int}=W_a$. There is a bijection between the set of integral conjugacy classes of $\tW^!_C$ and $\cd \cp_{\ge 0}/\sim$. 

{\bf Type D}.  There is a bijection between the set of integral conjugacy classes of $\tW^!_D$ and $\cd \cp_{\ge 0}/\sim$. An integral conjugacy class of $\tW^!_D$ intersects $\tW(D_n)$ if and only if it is represented by $(\tilde \l, \tilde \mu)$ with $|\mu|$ even. In this case, the intersection is a single conjugacy class of $\tW(D_n)$ except the case where $\tilde \mu$ is the empty double partition and no entry of $\tilde \l$ is of the form $(b, 0)$ with $b$ an odd number. In the last case, the intersection is a union of two conjugacy classes of $\tW(D_n)$ and the automorphism $\iota=(n, -n)$ on $\tW(D_n)$ interchanges these two conjugacy classes.

\section{The partial conjugation action and the subset $\cw$}

\subsection{}\label{iw} We first discuss the case where only the simple reflections $s_i$ for $i \in S$ are used, i.e., the conjugation action of the finite Weyl group $W$ on $\tW^!$. This was done in \cite{He072}.

For $\tw \in {}^S \tW^!$, let $$I(\tw)=\max\{J \subset S; \forall j \in J, \exists j' \in J, \text{ such that } s_j \tw=\tw s_{j'}\}.$$ Set $$\cw=\{x \tw; \tw \in {}^S \tW', x \in W_{I(\tw)}\}.$$ We have that $W \subset \cw$. 

The element $t^\chi \in {}^S \tW$ if and only if $\chi$ is dominant, i.e., $$\chi \in P^\vee_+=\{\g \in P^\vee; \<\g, \a\> \ge 0, \forall \a \in \Phi^+\}.$$ In this case, $I(t^\chi)=\{j \in S; \<\chi, \a_j\>=0\}$. For $\chi \in P^\vee, w \in W$, $t^\chi w \in {}^S \tW$ if and only if $\chi \in P^\vee_+$ and $w \in {}^{I(t^\chi)} W$. In this case, \[\tag{1} I(t^\chi w)=\max\{J \subset I(t^\chi); \forall j \in J, \exists j' \in J, \text{ such that } s_j w=w s_{j'}\}.\]

\begin{prop}\label{par}
Consider the (partial) conjugation action of $W$ on $\tW^!$. Let $\co$ be a $W$-conjugacy class on $\tW^!$ and $\co_{min}$ be the set of minimal length elements in $\co$. Then there exists $\tx \in {}^S \tW^!$ such that $\co \cap \cw=W_{I(\tx)} \tx=\tx W_{I(\tx)}$. Moreover, 

(1) For each $\tw \in \co$, there exists $\tw' \in \co_{\min} \cap \cw$ such that $\tw \to \tw'$. 

(2) Let $\tw, \tw' \in \co_{\min}$, then $\tw \sim \tw'$. 

(3) if $\co \cap {}^S \tW \neq \emptyset$, then $\co \cap {}^S \tW^!=\{\tx\}$ is a single element. In this case, $\tw \to \tx$ for any $\tw \in \co$. 
\end{prop}

There exists $\t \in \Om$ with $\co \subset W_a \t$. Define the action of $W$ on $W_a$ by $(x, \tw) \mapsto x \tw \t x \i \t \i$. Then  the map $W_a \to \tW^!$ defined by $\tw \mapsto \tw \t$ is a $W$-equivariant length-preserving map. Notice that $(W_a, \tS)$ is a Coxeter group and $x \mapsto \t x \t \i$ is an automorphism on $W_a$ which maps $S$ to another subset of $\tS$. By \cite[Corollary 2.6]{He072}, $\co \cap \cw=W_{I(\tx)} \tx=\tx W_{I(\tx)}$ for some $\tx \in {}^S \tW^!$. Part (1) and (2) follows from \cite[Corollary 3.8]{He072}. Part (3) follows from \cite[Corollary 2.5]{He072} and \cite[Corollary 3.7]{He072}.  

\

In the rest of this section we will give an explicit description of $\cw$ for $\tW^!$ of classical type. 

\begin{lem}\label{qqq}
Let $\chi \in P^\vee_+$ and $w \in {}^{I(t^\chi)} W$. Let $\a \in \Phi$. Then $\a \in \Phi_{I(t^\chi w)}$ if and only if $\<\chi, w^{-n} \a\>=0$ for all $n \ge 0$. 
\end{lem}

If $\a \in \Phi_{I(t^\chi w)}$, then $w^{-n} \a \in \Phi_{I(t^\chi w)}$ for all $n$. So we have that $\<\chi, w^{-n} \a\>=0$. Now we prove the other direction. 

For $\b=\sum_{i \in S} a_i \a_i \in \Phi$, set $\supp(\b)=\{i \in S; a_i \neq 0\}$ and $ht(\b)=\sum_{i \in S} |a_i|$. Notice that if $i \in I(t^\chi)$, then $w \i \a_i \in \Phi^+$. In particular, $ht(w \i \a_i) \ge ht(\a_i)=1$. Thus if $\b \in \Phi_{I(t^\chi)}$, then $wt(w \i \b) \ge wt(\b)$ and the equality holds if and only if for any $i \in \supp(\b)$, $w \i \a_i=\a_j$ for some $j \in S$. 

Assume that $\<\chi, w^{-n} \a\>=0$ for all $n \ge 0$. Then $w^{-n} \a \in \Phi_{I(t^\chi)}$ for all $n \ge 0$. Hence there exists $N \ge 0$ such that $ht(w^{-N} \a)=ht(w^{-N-1} \a)=\cdots$. So for $i \in \supp(w^{-N} \a)$ and $n \in \NN$, $w^{-n} \a_i=\a_j$ for some $j \in S$. Also we have that $\supp(w^{-N} \a) \subset I(t^\chi)$. By $\S$\ref{iw}(1), $\supp(w^{-N} \a) \subset I(t^\chi w)$. Hence $w^{-N} \a \in \Phi_{I(t^\chi w)}$ and $\a \in \Phi_{I(t^\chi w)}$. 

\begin{prop}\label{tw}
Let $\chi \in P^\vee$ and $w \in W$. Then $t^\chi w \in \cw$ if and only if for any $\a \in \Phi^+$, either $\<\chi, w^{-n} \a\>=0$ for all $n \ge 0$ or there exists $n \ge 0$ such that $\<\chi, \a\>=\<\chi, w \i \a\>=\cdots=\<\chi, w^{-n+1} \a\>=0$ and $\<\chi, w^{-n} \a\> >0$.
\end{prop}

If $t^\chi w \in \cw$, then $\chi \in P^\vee_+$ and $w=x w_1$ for some $w_1 \in {}^{I(t^\chi)} W$ and $x \in W_{I(t^\chi w_1)}$. Then $w^2=x (w_1 x w_1 \i) (w_1)^2$. Since $x \in W_{I(t^\chi w_1)}$, $w_1 x w_1 \i \in W_{I(t^\chi w_1)}$. Thus $w^2=x_2 w_1^2=w_1^2 x'_2$ for some $x_2, x'_2 \in W_{I(t^\chi w_1)}$. One can show in the same way that in general, $w^n=x_n w_1^n=w_1^n x'_n$ for some $x_n, x'_n \in W_{I(t^\chi w_1)}$. Hence for any $\a \in \Phi$, $$\<\chi, w^{-n} \a\>=\<\chi, (x'_n) \i w_1^{-n} \a\>=\<x'_n \chi, w_1^{-n} \a\>=\<\chi, w_1^{-n} \a\>.$$

By Lemma \ref{qqq}, if $\a \in \Phi^+_{I(t^\chi w_1)}$, then $w_1^{-n} \a \in \Phi^+_{I(t^\chi w_1)}$ and $\<\chi, w^{-n} \a\>=\<\chi, w_1^{-n} \a\>=0$. If $\a \in \Phi^+-\Phi^+_{I(t^\chi w_1)}$, then $\<\chi, w_1^{-n} \a\> \neq 0$ for some $n \ge 0$. In other words, there exists $N \ge 0$, such that $\<\chi, \a\>=\<\chi, w_1 \i \a\>=\cdots=\<\chi, w_1^{-N+1} \a\>=0$ and  $\<\chi, w_1^{-N} \a\> \neq 0$. We have that $\a, w_1 \i \a, \cdots, w_1^{-N+1} \a \in \Phi_{I(t^\chi)}$. Since $w_1 \in {}^{I(t^\chi)} W$, $w_1 \i \Phi^+_{I(t^\chi)} \subset \Phi^+$. Therefore $\a, w_1 \i \a, \cdots, w_1^{-N+1} \a \in \Phi^+_{I(t^\chi)}$ and $w_1^{-N} \a \in \Phi^+$. Since $\<\chi, w_1^{-N} \a\> \neq 0$, $\<\chi, w_1^{-N} \a\> >0$. 

On the other hand, assume that for $\a \in \Phi^+$, either $\<\chi, w^{-n} \a\>=0$ for all $n \ge 0$ or there exists $n \ge 0$ such that $\<\chi, \a\>=\<\chi, w \i \a\>=\cdots=\<\chi, w^{-n+1} \a\>=0$ and $\<\chi, w^{-n} \a\> >0$. Then $\chi \in P^\vee_+$. 

Let $\a \in \Phi^+$ and $n \ge 0$ with $\<\chi, \a\>=\<\chi, w \i \a\>=\cdots=\<\chi, w^{-n+1} \a\>=0$ and $\<\chi, w^{-n} \a\> >0$. If $w^{-n+1} \a \in \Phi^-$, then by our assumption for $-w^{-n+1} \a$, we have that $\<\chi, w^{-n} \a\> \le 0$, which is a contradiction. So $w^{-n+1} \a \in \Phi^+$. One can show by induction that $w^{-i} \a \in \Phi^+$ for $i \le n$. In particular, $w \i \a \in \Phi^+$. 

We have that $w=x w_1$ for some $w_1 \in {}^{I(t^\chi)} W$ and $x \in W_{I(t^\chi)}$. Let $x=s_{i_1} \cdots s_{i_l}$ be a reduced expression of $x$. Then for $\a=\a_{i_1}, s_{i_1} \a_{i_2}, \cdots, \allowbreak s_{i_1} \cdots s_{i_{l-1}} \a_{i_l}$, $w \i \a=w_1 \i x \i \a \in w_1 \i \Phi^-_{I(t^\chi)} \subset \Phi^-$. Thus for $n \ge 0$, $\<\chi, w^{-n} \a\>=0$. Notice that $s_{i_1} \a_{i_2}=\a_{i_2}+a \a_{i_1}$ for some $a \in \mathbb Z$. Since $\<\chi, w^{-n} \a_{i_1}\>=\<\chi, w^{-n} (\a_{i_2}+a \a_{i_1})\>=0$, we have that $\<\chi, w^{-n} \a_{i_2}\>=0$. One can show in the same way that $\<\chi, w^{-n} \a_i\>=0$ for $i=i_1, \cdots, i_l$ and $n \ge 0$. 

We have that $w^n=(x w_1)^n=(w_1)^n y$, where $$y=(w_1^{-n} x w_1^n) (w_1^{-n+1} x w_1^{n-1})  \cdots (w_1 \i x w_1)$$ is generated by $s_\a$ for $<\chi, \a>=0$. In particular, $y \chi=\chi$. Now $$\<\chi, w^{-n} \a\>=\<\chi, y \i w_1^{-n}\a\>=\<y \chi, w_1^{-n} \a\>=\<\chi, w_1^{-n} \a\>.$$ By the previous lemma, $i \in I(t^\chi w_1)$ for $i=i_1, \cdots, i_l$. Thus $x \in W_{I(t^\chi w_1)}$ and $t^\chi w \in \cw$. 

\subsection{}\label{xxxx} Let $\tw=t^{[a_1, \cdots, a_n]} \s$ for $a_i \in \mathbb Z$ and $\s \in S'_n$. Set $$\ua_i(\tw)=(a_i, a_{\s \i(i)}, \cdots) \in \mathbb Z^\infty,$$ here $a_{-i}=-a_i$ for $i \in \{1, \cdots, n\}$. Then for $\t \in S'_n$, $\t \tw \t \i=t^{[a_{\t \i(1)}, \cdots, a_{\t \i(n)}]} \t \s \t \i$ and \[\tag{1} \ua_{\t(i)}(\t \tw \t \i)=(a_i, a_{\s \i(i)}, \cdots, )=\ua_i(\tw).\]

The explicit description of $\cw$ for classical types follows from Proposition \ref{tw}. 

{\bf Type $A$.} If $\s \in S_n$, then for $\t \in S_n$, $\t \tw \t \i \in \cw$ if and only if for $1 \le i, j \le n$ with $\ua_i(\tw)>\ua_j(\tw)$, we have that $\t(i)<\t(j)$, in other words, $\ua_{\t \i(1)}(\tw) \ge \cdots \ge \ua_{\t \i(n)}(\tw)$. 

{\bf Type $B$ and $C$.} For $\t \in S'_n$, $\t \tw \t \i \in \cw$ if and only if $\ua_{\t \i(1)}(\tw) \ge \cdots \ge \ua_{\t \i(n)}(\tw) \ge (0, 0, \cdots)$.

{\bf Type $D$.} For $\t \in W(D_n)$, $\t \tw \t \i \in \cw$ if and only if $\ua_{\t \i(1)}(\tw) \ge \cdots \ge \ua_{\t \i(n)}(\tw)$ and $\ua_{\t \i(n-1)}+\ua_{\t \i(n)} \ge (0, 0, \cdots)$. This condition is equivalent to $\ua_{\t \i(1)}(\tw)  \ge \cdots \ge \ua_{\t \i(n-1)}(\tw) \ge (0, 0, \cdots)$ and $\ua_{\t \i(n-1)}(\tw) \ge \pm \ua_{\t \i(n)}(\tw)$.

\subsection{}\label{qs} Let $\tw=t^\chi \s$ with $\chi \in \ZZ^n$ and $\s \in S'_n$. We say that $\tw$ is {\it quasi-positive} if it satisfies

(1) $\ua_i(\tw) \ge (0, 0, \cdots)$ for all $i$. 

(2) If $\ua_i(\tw)=(0, 0, \cdots)$, then $\s \i(j)>0$ for $j<\max\{ |(\s)^l(i)|; l \in \ZZ\}$. 

So if $\tw$ is quasi-positive, then in particular all the entries of $\chi$ are non-negative. 

Let $\tw \in \ZZ^n \rtimes S'_n$ and $\t=(t_1, \cdots, t_n) \in (\ZZ/2 \ZZ)^n \subset S'_n$. If $\t \tw \t \i$ is quasi-positive, then by condition (1), $t_i>0$ if $\ua_i(\tw)>(0, 0, \cdots)$ and $t_i<0$ if $\ua_i(\tw)<(0, 0, \cdots)$. 
The reason for condition (2) here is to guarantee that each orbit of $(\ZZ/2 \ZZ)^n$ on $\ZZ^n \rtimes S'_n$ (for the conjugation action) contains a unique quasi-positive element. 

For example, the element $(1, 2, -1, -2), (1, -2, -1, 2) \in S'_2 \subset \ZZ^2 \rtimes S'_2$ are in the same orbit of $(\ZZ/2 \ZZ)^2$ and both satisfy condition (1), but only $(1, -2, -1, 2)$ satisfies condition (2). Thus $(1, -2, -1, 2)$ is quasi-positive, but $(1, 2, -1, -2)$ is not quasi-positive. 

It is also easy to see that 

(3) Any two quasi-positive elements in a given $S'_n$-conjugacy class are conjugated by an element in $S_n$. 

However, not all the element conjugated to a quasi-positive element by $S_n$ are still quasi-positive. In the above example, $(1, 2, -1, -2)=(1 2) (1, -2, -1, 2) (1 2)$ is not quasi-positive. 

\

The following variation of Proposition \ref{par} for classical types is very useful in the study of minimal length elements. 

\begin{cor}\label{xxxxx}
Let $\tW^!$ be of classical type and $\tw \in \tW^!_{int}$.

(1) If $\tW$ is of type A, then there exists $\t \in S_n$ such that $\tw \to \t \tw \t \i$ and $\ua_{\t \i(1)}(\tw) \ge \cdots \ge \ua_{\t \i(n)}(\tw)$. If $\tw$ is quasi-positive, then $\t \tw \t \i$ is automatically quasi-positive. 

(2) If $\tW^!$ is of type B, C or D, then there exists $\t \in S'_n$ such that $\tw \tilde \to \t \tw \t \i$, $\t \tw \t \i$ is quasi-positive and $\ua_{\t \i(1)}(\tw) \ge \cdots \ge \ua_{\t \i(n)}(\tw) \ge (0, 0, \cdots)$. If moreover, $\tw$ is quasi-positive, then we may choose such $\t$ in $S_n$. 
\end{cor}

By Proposition \ref{par}, there exists $\t \in W$ such that $\tw \to \t \tw \t \i$ and $\t \tw \t \i \in \cw$. If $\tW^!$ is of type $A$, by $\S$\ref{xxxx} we have that $\ua_{\t \i(1)}(\tw) \ge \cdots \ge \ua_{\t \i(n)}(\tw)$. If moreover, $\tw \in \tW^!$ is quasi-positive, then $\ua_i(\tw), \cdots, \ua_n(\tw) \ge (0, 0, \cdots)$ and $\ua_1(\t \tw \t \i), \cdots, \ua_n(\t \tw \t \i) \ge (0, 0, \cdots)$ for all $\t \in S_n$. This finishes the proof for type A.

Now assume that $\tW^!$ is of type B, C or D. We first show that 

(a) There exists $\s \in S'_n$ such that $\tw \tilde \to \s \tw \s \i$ and $\ua_{\s \i(1)}(\tw) \ge \cdots \ge \ua_{\s \i(n)}(\tw) \ge (0, 0, \cdots)$.

Let $\t \in W$ with $\tw \to \t \tw \t \i$ and $\t \tw \t \i \in \cw$. By $\S$\ref{xxxx}, $\ua_{\t \i(1)}(\tw) \ge \cdots \ge \ua_{\t \i(n)}(\tw) \ge (0, 0, \cdots)$ for type B and C. 

Now we consider type D. By $\S$\ref{xxxx}, $\ua_{\t \i(1)}(\tw)  \ge \cdots \ge \ua_{\t \i(n-1)}(\tw) \ge (0, 0, \cdots)$ and $\ua_{\t \i(n-1)}(\tw) \ge \pm \ua_{\t \i(n)}(\tw)$. The statement holds if $\ua_{\t \i(n)}(\tw) \ge (0, 0, \cdots)$. If $\ua_{\t \i(n)}(\tw)<(0, 0, \cdots)$, then $\ua_{\t \i(n-1)}(\tw) \ge -\ua_{\t \i(n)}(\tw) \ge (0, 0, \cdots)$. 

Notice that $S'_n=W(D_n) \rtimes <\iota>$ and by definition, $\t \tw \t \i \tilde \approx \iota \t \tw \t \i \iota \i$. Since $(\iota \t) \i(i)=\t \i \iota \i(i)=\begin{cases} i, & \text{ if } i \neq n \\ -n, & \text{ if } i=n \end{cases}$, we have that $\ua_{\t \i(i)}(\tw)=\ua_{(\iota \t) \i(i)}(\tw)$ for $i \neq n$ and $\ua_{(\iota \t) \i(n)}(tw)=-\ua_{\t \i(n)}(\tw)$. So $\tw \tilde \to (\iota \t) \tw (\iota \t) \i$ and $\ua_{(\iota \t) \i(1)}(\tw) \ge \cdots \ge \ua_{(\iota \t) \i(n)}(\tw) \ge (0, 0, \cdots)$. 

(a) is proved. 

Let $\tw_1=\s \tw \s \i$ with $\s$ in (a). We may assume that $\tw_1=t^{[a_1, \cdots, a_n]} w$ with $a_i \in \ZZ$ and $w \in S'_n$. Then there exists $0 \le i \le n$ such that $$\ua_1(\tw_1) \ge \cdots \ge \ua_i(\tw_1)>(0, 0, \cdots)=\ua_{i+1}(\tw_1)=\cdots=\ua_n(\tw_1).$$ We show that 

(b) $a_{i+1}=\cdots=a_n=0$ and $w \in S'_i \times S'_{n-i} \subset S'_n$. Here $S'_i$ is the the group of permutations on $\{\pm 1, \cdots, \pm i\}$ and $S'_{n-i}$ is the group of permutations on $\{\pm(i+1), \cdots, \pm n\}$. 

By definition, $\ua_j(\tw_1)=(a_j, \ua_{w \i(j)}(\tw_1))$. Thus for $j>i$, $a_j=0$ and $\ua_{w \i(i)}(\tw_1)=(0, \cdots, 0)$. Since $\ua_k(\tw_1)>(0, 0, \cdots)$ for $k \le i$, we must have that $w \i(i) \in \{\pm(i+1), \cdots, \pm n\}$. Since $w$ is a permutation of finite order on $\{\pm 1, \cdots, \pm n\}$, each orbit of $w$ is a finite set and is contained in $\{\pm(i+1), \cdots, \pm n\}$ or contains no elements in $\{\pm(i+1), \cdots, \pm n\}$. Thus $w \in S'_i \times S'_{n-i}$. 

(b) is proved. 

Let $w=w_1 w_2$ with $w_1 \in S'_i$ and $w_2 \in S'_{n-i}$. Set $\tw'=t^{[a_1, \cdots, a_n]} w_1$. Then for any $w' \in S'_{n-i}$, $\tw' w'=w' \tw'$ and $l(\tw' w')=l(\tw')+l(w')$. Here $l$ is the length function on $\tW^!$. So for any $w', w'' \in S'_{n-i}$ and a simple reflection of $W$ that is in $S'_{n-i}$, $w' \xrightarrow i w''$ if and only if $\tw' w' \xrightarrow i \tw' w''$. 

By \cite[Proposition 2.3]{GP93}, $w_2 \to w'_2$ for some $w'_2$ of the form $(i+b_1, -(i+b_1))^{\e_1} (i+b_1+b_2, -(i+b_1+b_2))^{\e_2} \cdots (n, -n)^{\e_k} (i+1, i+2, \cdots, i+b_1) (i+b_1+1, i+b_1+2, \cdots, i+b_1+b_2) \cdots (n-b_k+1, n-b_k+2, \cdots, n)$, with $b_1+\cdots+b_k=n-i$ and $\e_1, \cdots, \e_k \in \{0, 1\}$. Here $(i+1, i+2, \cdots, i+b_1)$ is a positive block and $(i+b_1, -(i+b_1)) (i+1, i+2, \cdots, i+b_1)$ is a negative block in the sense of \cite[2.2]{GP93}. 

Let $\s' \in S'_{n-i}$ with $w'_2=\s' w_2 (\s') \i$. Then $\tw' w'_2=\s' \tw_1 (\s') \i=(\s' \s) \tw (\s' \s) \i$ and $\tw \tilde \to \tw_1 \tilde \to \tw' w'_2$. Moreover, $\ua_i(\tw' w'_2)=\ua_i(\tw_1)$ for any $i$. So $\ua_1(\tw' w'_2) \ge \cdots \ua_n(\tw' w'_2) \ge (0, 0, \cdots)$. The condition (2) in the definition of quasi-positivity is also satisfied for $\tw' w'_2$. Thus $\tw' w'_2$ is quasi-positive. The ``moreover'' part follows from $\S$\ref{qs} (3). 

\section{$P$-operators}

In this section, we introduce $P$-operators on $\ZZ^n \rtimes S_n$ and $\ZZ^n \rtimes S'_n$ and discuss some relations between $P$-operators and terminal elements in classical affine Weyl groups. We first discuss some general results on terminal elements. 

\begin{lem}\label{gen}
Let $\tw, \tx \in \tW$ with $l(\tx \i \tw)=l(\tw)-l(\tx)$. If $\tw$ is a terminal element, then $\tw \tilde \approx \tx \i \tw \tx$. 
\end{lem}

We argue by induction on $l(\tx)$. If $l(\tx)=0$, then the statement is obvioius. Now assume that $\tx=s_i \tx'$ for some $i \in \tilde I$ with $\tx'<\tx$. Since $l((\tx') \i s_i \tw)=l(\tw)-l(\tx')-1$, we have that $l(s_i \tw)=l(\tw)-1$ and $\tw \tilde \to s_i \tw s_i$. Since $\tw$ is terminal, $l(s_i \tw s_i)=l(\tw)$. Notice that \begin{align*} l(\tx \i \tw)+1 & \ge l(\tx \i \tw s_i)=l((\tx') \i s_i \tw s_i) \ge l(s_i \tw s_i)-l(\tx') \\ &=l(\tw)-l(\tx)+1=l(\tx \i \tw)+1.\end{align*} Therefore $l((\tx') \i s_i \tw s_i)=l(\tw)-l(\tx')$. By induction hypothesis, $$\tw \approx s_i \tw s_i \tilde \approx (\tx') \i s_i \tw s_i \tx'=\tx \i \tw \tx.$$ 

\subsection{} Let $\tw=t^\chi w$ with $\chi \in P^\vee_+$ and $w \in W$. For $\g \in P^\vee_+$, we say that $\g$ is {\it minuscule} for $\tw$ if for any $\a \in \Phi^+$ with $\<\g, \a\> \ge 2$, $$\<\chi-\g, \a\> \ge \begin{cases} -1, & \text{ if } w \i \a \in \Phi^+ \\ 0, & \text{ if } w \i \a \in \Phi^- \end{cases}.$$

If $\g$ is a minuscule coweight, then $\<\g, \a\> \le 1$ for all $\a \in \Phi^+$ and the condition above is automatically satisfied. On the other hand, for $\tw=1$, we have that $\<\chi-\g, \a\>=-\<\chi, \a\>$. In this case, $\g$ is minuscule for $\tw=1$ if and only if $\g=0$ or is a minuscule coweight. 

The following result provides some element $\t$ with $l(\t \i \tw)=l(\tw)-l(\t)$. 

\begin{lem}\label{-}
Let $\tw=t^\chi w$ with $\chi \in P^\vee_+$ and $w \in W$. Let $\g \in P^\vee_+$ such that $\g$ is minuscule for $\tw$. Set $\t=t^\g w_{I(t^\g)} w_S$. Then $$l(\t \i \tw)=l(\tw)-l(\t).$$
\end{lem}

We have that $\t \i \tw=w_S w_{I(t^\g)} t^{\chi-\g} w$. Notice that for any $\a \in \Phi^+$, $w_S w_{I(t^\g)} \a \in \Phi^+$ if and only if $\a \in \Phi^+_{I(t^\g)}$. By definition, \begin{align*} l(\t \i \tw)=& \sum_{\a \in \Phi^+_{I(t^\g)}, w \i \a \in \Phi^+} |\<\chi-\g, \a\>|+\sum_{\a \in \Phi^+_{I(t^\g)}, w \i \a \in \Phi^-} |\<\chi-\g, \a\>-1| \\ &+\sum_{\a \in \Phi^+-\Phi^+_{I(t^\g)}, w \i \a \in \Phi^+} |\<\chi-\g, \a\>+1| \\ &+\sum_{\a \in \Phi^+-\Phi^+_{I(t^\g)}, w \i \a \in \Phi^-} |\<\chi-\g, \a\>| \end{align*}

For $\a \in \Phi^+_{I(t^\g)}$, $|\<\chi-\g, \a\>|=|\<\chi, \a\>|$ and $|\<\chi-\g, \a\>-1|=|\<\chi, \a\>-1|$. 

For $\a \in \Phi^+-\Phi^+_{I(t^\g)}$ with $w \i \a \in \Phi^+$, $|\<\chi-\g, \a\>+1|=\<\chi, \a\>-(\<\g, \a\>-1)=|\<\chi, \a\>|-|\<\g, \a\>-1|$. 

For $\a \in \Phi^+-\Phi^+_{I(t^\g)}$ with $w \i \a \in \Phi^-$, if $\<\chi-\g, \a\> \ge 0$, then $\<\chi, \a\> \ge \<\g, \a\> \ge 1$. In this case, $|\<\chi-\g, \a\>|=\<\chi-\g, \a\>=|\<\chi, \a\>-1|-|\<\g, \a\>-1|$. If $\<\chi-\g, \a\><0$, then $\<\chi, \a\>=0$ and $\<\g, \a\>=1$ since $\g$ is minuscule for $\tw$. In this case, $|\<\chi-\g, \a\>|=1=|\<\chi, \a\>-1|-|\<\g, \a\>-1|$.

Therefore $l(\t \i \tw)=l(\tw)-l(\t)$. 

\subsection{}\label{se1} We recall some basic properties on the lexicographic order. 

Let $(P, \le)$ be a total order set. 

For two sequences $\mu=(a_1, a_2, \cdots, a_n)$ and $\mu'=(a'_1, a'_2, \cdots, a'_n)$ in $P^n$, we say that $\mu>\mu'$ if there exists $1 \le i \le n$ such that $a_j=a'_j$ for $j<i$ and $a_i>a'_i$. For any sequence $\mu=(a_1, \cdots, a_n) \in P^n$, define $|\mu|=\sum_{i=1}^n a_i$ and $\mu[i]=(a_{i+1}, a_{i+2}, \cdots, a_n, a_1, a_2, \cdots, a_i)$. 

For two sequences $\th=(a_1, a_2, \cdots,)$ and $\th'=(a'_1, a'_2, \cdots,)$ in $P^\infty$, we say that $\th>\th'$ if there exists $i \in \NN$ such that $a_j=a'_j$ for $j<i$ and $a_i>a'_i$. We define $\th[i]=(a_{i+1}, a_{i+2}, \cdots) \in P^\infty$. 

It is easy to see that 

(1) If $\th \in P^n$, then $\th[i]=\th$ if and only if $\th$ is of the from $(\g, \g, \cdots, \g)$ for some $\g \in P^{\gcd(n, i)}$. 

(2) If $\th \in P^\infty$, then $\th[i]=\th$ if and only if $\th$ is of the form $(\g, \g, \cdots)$ for some $\g \in \ZZ^i$. 

(3) If $\a, \a' \in P^n$ and $\b, \b' \in P^m$ with $\a \ge \a'$ and $\b \ge \b'$, we have that $(\a, \b) \ge (\a', \b')$. 

(4) If $\th=(a_1, a_2, \cdots) \in P^\infty$ and $\g \in P^n$ is an upper bound of length-$n$ subsequences in $\th$ (i.e., $\g \ge (a_l, a_{l+1}, \cdots, a_{l+n-1})$ for all $l$). Then $(\g, \th) \ge \th$. 

\begin{lem}
Let $(P, \le)$ be a total order set. For any $n \in \NN$, define the action of $S_n$ on $P^n$ by $\t \cdot (p_1, \cdots, p_n)=(p_{\t \i(1)}, \cdots, p_{\t \i(n)})$. If $a, a' \in P^n$, $b, b' \in P^m$ with $\max\{S_n \cdot a\} \ge \max\{S_n \cdot a'\}$ and $\max\{S_m \cdot b\} \ge \max\{S_m \cdot b'\}$. Then $\max\{S_{n+m} \cdot (a, b)\} \ge \max\{S_{n+m} \cdot (a', b')\}$. 
\end{lem}

We argue by induction on $n+m$. The lemma holds for $n+m=1$. Now assume that $n+m=k>0$ and the lemma holds for $n+m=k-1$. Let $a=(a_1, \cdots, a_n), a'=(a'_1, \cdots, a'_n)$ and $b=(b_1, \cdots, b_m), b'=(b'_1, \cdots, b'_m)$ with $a_i, a'_i, b_j, b'_j \in P$. Without loss of generality, we may assume that $a_1 \ge a_i$ for any $i$ and $a_1 \ge b_j$ for any $j$. Then $a_1 \ge \max_i a'_i$ and $a_1 \ge \max_j b'_j$. 

If $a_1>\max_i a'_i$ and $a_1>\max_j b'_j$, then $\max\{S_k \cdot (a, b)\}>\max\{S_k \cdot (a', b')\}$. 

If $a_1=\max_i a'_i$, then we may assume that $a'_i=a_1$. Set $c=(a_2, \cdots, a_n)$ and $c'=(a'_1, \cdots, \hat{a'_i}, \cdots, a'_n)$. Then $$\max(S_n \cdot a)=(a_1, \max\{S_{n-1} \cdot c\}) \text{ and } \max(S_n \cdot a')=(a_1, \max\{S_{n-1} \cdot c'\}).$$ So $\max\{S_{n-1} \cdot c\} \ge \max\{S_{n-1} \cdot c'\}$. By induction hypothesis, $\max\{S_{k-1} \cdot (c, b)\} \ge \max\{S_{k-1} \cdot (c', b')\}$. So \begin{align*} \max\{S_k \cdot (a, b)\} &=(a_1, \max\{S_{k-1} \cdot (c, b)\}) \ge (a_1, \max\{S_{k-1} \cdot (c', b')\}) \\ &=\max\{S_k \cdot (a', b')\}.\end{align*}

If $a_1>\max_i a'_i$ but $a_1=\max_j b'_j$. Then $\max_j b_j=\max_j b'_j=a_1$. By the same argument as we did in the previous case (for $b, b'$ instead of $a, a'$), we also have that $\max\{S_k \cdot (a, b)\}>\max\{S_k \cdot (a', b')\}$.  

\

In the rest of this section, we introduce $P$-operators for classical types. This is related to Lemma \ref{-}. 

\subsection{}\label{ps} Let $\tw=t^\chi \s$ be a quasi-positive element with $\chi=[a_1, \cdots, a_n]$ for $a_i \in \ZZ$ and $\s \in S'_n$. Let $\th \in \ZZ^\infty$ with positive leading entry. In this case, we define $e_\th(\tw)=[c_1, \cdots, c_n]$, where $$c_i=\begin{cases} 1, & \text{ if } \ua_i(\tw) \ge \th; \\ 0, & \text{ otherwise}.\end{cases}$$ Then

(1) $\text{ for } \t \in S_n, \t e_\th(\tw)=[c_{\t \i(1)}, \cdots, c_{\t \i(n)}]=e_\th(\t \tw \t \i)$.

We let $P_{\th}(\tw)$ be the unique quasi-positive element in $(\ZZ/2 \ZZ)^n$-conjugacy class of $t^{-e_\th(\tw)} \tw t^{e_\th(\tw)}=t^{\chi-e_\th(\tw)+\s e_\th(\tw)} \s$. We call $P_\th$ a {\it $P$-operator}. 

In the special case where $\th>\ua_i(\tw)$ for all $i$, we have that $e_\th(\tw)=[0, \cdots, 0]$ and $P_\th(\tw)=\tw$. 

For any $a \in \NN$, It is easy to see that $e_{\th+[a, a, \cdots]}(t^{[a, \cdots, a]} \tw)=e_\th (\tw)$. Thus

(2) $t^{[a, \cdots, a]} P_\th(\tw)=P_{\th+[a, a, \cdots]}(t^{[a, \cdots, a]} \tw)$. 

By (1), for any $\t \in S_n$, $\t(\chi-e_\th(\tw)+\s e_\th(\tw))=\t \chi-e_\th(\t \tw \t \i)+\t \s \t \i e_\th(\t \tw \t \i)$. Hence by $\S$\ref{qs} (3),

(3) For any $\t \in S_n$, $P_\th(\t \tw \t \i) \in S_n \cdot P_\th(\tw)$. 

\begin{lem}\label{sn}
Let $\tw=t^\chi \s$ be a quasi-positive element with $\chi=[a_1, \cdots, a_n]$ for $a_i \in \ZZ$ and $\s \in S_n$. Let $\th \in \ZZ^\infty$ with positive leading entry. Then $P_\th(\tw)=t^{-e_\th(\tw)} \tw t^{e_\th(\tw)}=t^{\chi-e_\th(\tw)+\s e_\th(\tw)} \s$.
\end{lem}

We assume that $e_\th(\tw)=[c_1, \cdots, c_n]$. Then $\chi-e_\th(\tw)+\s e_\th(\tw)=[a'_1, \cdots, a'_n]$, where $a'_i=a_i-c_i+c_{\s \i(i)}$. Since $\s \in S_n$, $\s \i(i)>0$ and $c_{\s \i(i)} \ge 0$ for $1 \le i \le n$. 

If $a_i>0$, then $a_i-c_i \ge 0$ and $a'_i \ge 0$. 

If $a_i=0$, then $\ua_i(\tw)<\th$ since the leading entry of $\th$ is positive. Hence $c_i=0$ and $a'_i \ge 0$. 

Therefore $\ua_i(t^{\chi-e_\th(\tw)+\s e_\th(\tw)} \s) \ge (0, 0, \cdots)$. Since $\s \in S_n$, condition (2) in the definition of quasi-positive element is automatically satisfied. Thus $t^{\chi-e_\th(\tw)+\s e_\th(\tw)} \s$ is quasi-positive and $P_\th(\tw)=t^{\chi-e_\th(\tw)+\s e_\th(\tw)} \s$.



\begin{prop}\label{red}
Let $\tW^!$ be of classical type. Let $\tw \in \ZZ^n \rtimes S'_n$ be quasi-positive and $\th \in \ZZ^\infty$ with positive leading entry. If $\tw' \in W \cdot \tw$ is a terminal element in $\tW^!$, then then there exists a terminal element $\tw_1 \in S_n \cdot P_\th(\tw)$ such that $\tw' \tilde \approx \tw_1$. 
\end{prop} 

We first prove for type A. By Corollary \ref{xxxxx}, then there exists $\t \in S_n$ such that $\tw' \to \t \tw \t \i$ and $\ua_{\t \i(1)}(\tw) \ge \cdots \ge \ua_{\t \i(n)}(\tw)$. Since $\tw'$ is terminal, $\t \tw \t \i$ is also terminal and $\tw' \approx \t \tw \t \i$. 

Let $0 \le i \le n+1$ with $\ua_{\t \i(i)}(\tw) \ge \th>\ua_{\t \i(i+1)}$. Set $\g=\sum_{j \le i} e_j$. Then $\g=e_\th(\t \tw \t \i)$ and $\t \i \g=e_\th(\tw)$. By definition, $\g$ is either $0$ or $[1, 1, \cdots, 1]$ or a fundamental coweight. Hence $\g$ is minuscule for $\t \tw \t \i$. 

By Lemma \ref{-}, there exists $x \in S_n$ such that $l(x \i t^{-\g} \t \tw \t \i)=l(\t \tw \t \i)-l(t^\g x)$. By Lemma \ref{sn}, $t^{-\t \i \g} \tw t^{\t \i \g}=P_\th(\tw)$. So by Lemma \ref{gen}, $$\t \tw \t \i \tilde \approx (t^\g x) \i \t \tw \t \i (t^\g x)=x \i \t P_\th(\tw) \t \i x \in S_n \cdot P_\th(\tw).$$

Since $\t \tw \t \i$ is terminal, $x \i \t P_\th(\tw) \t \i x$ is also terminal. The Proposition is proved for Type A. 

Now we assume that $\tW^!$ is of type B, C or D. By Corollary \ref{xxxxx}, then there exists $\tw'' \in S'_n \cdot \tw'=S'_n \cdot \tw$ such that $\tw' \tilde \to \tw''$, $\tw''$ is quasi-positive and $\ua_1(\tw'') \ge \cdots \ge \ua_n(\tw'') \ge (0, 0, \cdots)$. Since $\tw'$ is terminal, $\tw''$ is also terminal and $\tw' \tilde \approx \tw''$. Since $\tw$ is quasi-positive, by $\S$\ref{qs} (3), there exists $\t \in S_n$ such that $\tw''=\t \tw \t \i$. 

Let $0 \le i \le n+1$ with $\ua_{\t \i(i)}(\tw) \ge \th>\ua_{\t \i(i+1)}$. Set $\g=\sum_{j \le i} e_j$. Then $\g=e_\th(\t \tw \t \i)$ and $\t \i \g=e_\th(\tw)$. 

Assume that $\t \tw \t \i=t^{[a_1, \cdots, a_n]} \s$. Then $a_1 \ge \cdots \ge a_i \ge 1$. Notice that for any $\a \in \Phi^+$ with $\<\g, \a\> \ge 2$, we have that $\a=e_j+e_{j'}$ for some $j, j' \le i$ and $\<\g, \a\>=2$. It is easy to see that $\g$ is minuscule for $\t \tw \t \i$.

By Lemma \ref{-}, there exists $x \in W$ such that $l(x \i t^{-\g} \t \tw \t \i)=l(\t \tw \t \i)-l(t^\g x)$. So by Lemma \ref{gen}, $$\t \tw \t \i \tilde \approx (t^\g x) \i \t \tw \t \i (t^\g x)=x \i \t (t^{-e_\th(\tw)} \tw t^{e_\th(\tw)}) \t \i x \in S'_n \cdot P_\th(\tw).$$

By Corollary \ref{xxxxx}, there exists a quasi-positive $\tw_1 \in S_n \cdot P_\th(\tw)$ such that $x \i \t (t^{-e_\th(\tw)} \tw t^{e_\th(\tw)}) \t \i x \tilde \to \tw_1$. Since $\t \tw \t \i$ is terminal, $\tw_1$ is also terminal. 

\section{Properties of $P$-operators}

In this section, we study elements $t^\chi \s$ with $\s \in (1 2 \cdots n) (\ZZ/2 \ZZ)^n \subset S'_n$. As we will see in the next section, the general case can be essentially reduce to these elements. 

\subsection{}\label{chi} For any $m \in \ZZ$, set $\chi_{n, m}=(a_1, \cdots, a_n)$ and $\chi'_{n, m}=(a_n, a_{n-1}, \cdots, a_1)$, where $a_i=\lceil \frac{i m}{n} \rceil-\lceil \frac{(i-1) m}{n} \rceil$. Here $\lceil x \rceil=\min\{k \in \ZZ; k \ge x\}$. Set $$\tw_{n, m}=t^{\chi'_{n, m}} (1 2 \cdots n).$$ 
Then it is easy to see that $\ua_n(\tw_{n, m})=(\chi_{n, m}, \chi_{n, m}, \cdots)$. 

Set $\tw'_{n, 0}=(n, -n) (1 2 \cdots n)$. Then by definition $\tw'_{n, 0}$ is quasi-positive and for any $\th \in \ZZ^\infty$ with positive leading entry, $P_\th(\tw'_{n, 0})=\tw'_{n, 0}$. 

For any $m \in \NN$ with $m \mid n$, set $$\tw'_{n, m}=t^{\chi'_{n, m}}  (n/m, -n/m) (2 n/m, -2 n/m) \cdots (n, -n) (1 2 \cdots n) .$$ 
It is easy to see that $\tw'_{n, m}$ is quasi-positive and $\ua_n(\tw'_{n, m})=(\g, -\g, \g, -\g, \cdots)$, here $\g=(1, 0, 0, \cdots, 0) \in \ZZ^{n/m}$. 

\subsection{} We give two examples of $\tw'_{n, m}$ here. 

(1) $\chi_{6, 2}=(1, 0, 0, 1, 0, 0)$ and $\tw'_{6, 2}=t^{[0, 0, 1, 0, 0, 1]} (1, 2, -3, -4, -5, 6)$. For $\th \in \ZZ^\infty$ with positive leading entry and $\th \le \ua_6(\tw'_{6, 2})$, we have that $e_\th(\tw'_{6, 2})=[0, 0, 1, 0, 0, 1]$ and $P_\th(\tw'_{6, 2})$ is the unique quasi-positive element in $(\ZZ/2 \ZZ)^6$-conjugacy class of $t^{[1, 0, 0, 1, 0, 0]} (1, 2, -3, -4, -5, 6)$. Hence $$P_\th(\tw'_{6, 2})=t^{[1, 0, 0, 1, 0, 0]} (1, 2, 3, -4, -5, -6)=(1 2 3 4 5 6) \tw'_{6, 2} (1 2 3 4 5 6) \i.$$ 

(2) $\chi_{3, 3}=(1, 1, 1)$ and $\tw'_{3, 3}=t^{[1, 1, 1]} (1, -2, 3, -1, 2, -3)$. For $\th \in \ZZ^\infty$ with positive leading entry and $\th \le \ua_3(\tw'_{3, 3})$, we have that $e_\th(\tw'_{3, 3})=[1, 1, 1]$ and $P_\th(\tw'_{3, 3})$ is the unique quasi-positive element in $(\ZZ/2 \ZZ)^3$-conjugacy class of $t^{[-1, -1, -1]} (1, -2, 3, -1, 2, -3)$. Hence $$P_\th(\tw'_{3, 3})=t^{[1, 1, 1]} (1, -2, 3, -1, 2, -3)=(1 2 3) \tw'_{3, 3} (1 2 3) \i.$$ 

In both examples, if $\th>\ua_n(\tw'_{n, m})$, then by $\S$\ref{ps}, $P_\th(\tw'_{n, m})=\tw'_{n, m}$. 

In general, $\tw'_{n, m}$ is ``stable under the $P$-operator in the following sense by direct calculation. 

\begin{prop}\label{stable1}
Let $m \in \NN$ with $m \mid n$ and $\th \in \ZZ^\infty$ with positive leading entry positive. Then $\tw'_{n, m}$ is quasi-positive and $P_\th(\tw'_{n, m})=\tw'_{n, m}$ or $(1 2 \cdots n) \tw'_{n, m} (1 2 \cdots n) \i$. 
\end{prop}


\

The element $\tw_{n, m}$ is also ``stable'' under the $P$-operators. The proof will be given in $\S$\ref{pstable}.

\begin{prop}\label{stable}
Let $m \in \NN$ and $\th \in \ZZ^\infty$. Then there exists $i \in \ZZ$ such that $$P_\th(\tw_{n, m})=(1 2 \cdots n)^i \tw_{n, m} (1 2 \cdots n)^{-i}.$$
\end{prop}

\begin{lem} For any $m \in \ZZ$, we have that 

(1) $\chi_{n, m} \ge \chi_{n, m}[i]$ for all $i$.

(2) For any $\chi' \in \ZZ^n$ with $|\chi'|=m$, $\max_i \chi'[i] \ge \chi_{n, m}$. 
\end{lem}

Let $a_i=\lceil \frac{i m}{n} \rceil-\lceil \frac{(i-1) m}{n} \rceil$. By definition, \[\tag{a} \lceil x+y \rceil \le \lceil x \rceil +\lceil y \rceil \le \lceil x+y \rceil+1.\] If $\chi_{n, m}[i]>\chi_{n, m}$, then there exists $1 \le j \le n$ such that $a_{i+k}=a_k$ for $k<j$ and that $a_{i+j}>a_j$. Therefore $$\lceil \frac{(i+j) m}{n} \rceil-\lceil \frac{i m}{n} \rceil=a_{i+1}+\cdots+a_{i+j}>a_1+\cdots+a_j=\lceil \frac{j m}{n} \rceil.$$ This contradicts (a). So $\chi_{n, m}[i] \le \chi_{n, m}$ for all $i$ and part (1) is proved. 

We prove part (2) by induction on $n$. The case where $n=1$ is trivial. Assume that part (2) holds for all $n'<n$. Let $\chi'=(a'_1, \cdots, a'_n)$ with $\sum a'_i=m$. After rearranging $\{1, \cdots, n\}$, we may assume that $\chi' \ge \chi'[i]$ for all $i$. 

If $\chi'<\chi_{n, m}$, then there exists $1<i<n$ such that $a_j=a'_j$ for $j<i$ and $a'_i<a_i$. So $\sum_{j=1}^i a'_i<\sum_{j=1}^i a_i=\lceil \frac{i m}{n} \rceil$ and $$\sum_{j=i+1}^n a'_j \ge m-\lceil \frac{i m}{n} \rceil+1 \ge \lceil \frac{(n-i) m}{n} \rceil=\sum_{j=1}^{n-i} a_i.$$

By induction hypothesis, there exists $k \ge i$, such that $$(a'_{k+1}, a'_{k+2}, \cdots, a'_n, a'_{i+1}, a'_{i+2}, \cdots, a'_k) \ge (a_1, \cdots, a_{n-i})$$ in $\ZZ^{n-i}$. In particular, $(a'_{k+1}, \cdots, a'_n) \ge (a_1, \cdots, a_{n-k}) \ge (a'_1, \cdots, a'_{n-k})$. 

By our assumption on $\chi'$, $(a'_1, \cdots, a'_k)$ is an upper bound of the length-$k$ subsequence of $(\chi', \chi', \cdots) \in \ZZ^\infty$. Hence $(a'_{k+1}, \cdots, a'_n)$ is also an upper bound of the length-$k$ subsequence of $(\chi', \chi', \cdots)$. By $\S$\ref{se1} (4), $(\chi'[k], \chi'[k], \cdots)=(a'_{k+1}, \cdots, a'_n, \chi', \chi', \cdots) \ge (\chi', \chi', \cdots)$. In particular, $\chi'[k] \ge \chi'$. Therefore $\chi'=\chi'[k]$.  Let $d=\gcd(k, n)$. Then there exists $\g \in \ZZ^d$, such that $\chi'=(\g, \g, \cdots, \g)$. Since $|\chi'|=m$ and $\chi'<\chi_{n, m}$, we have that $d \mid m$ and that $\g \le (a_1, \cdots, a_d)$. On the other hand, $(a'_{k+1}, \cdots, a'_n) \ge (a_1, \cdots, a_{n-k})$. So $\g \ge (a_1, \cdots, a_d)$. Hence $\g=(a_1, \cdots, a_d)$. By the definition of $\chi_{n, m}$, we have that $\chi_{n, m}=(\g, \g, \cdots, \g)=\chi'$, which is a contradiction. Part (2) is proved. 

\begin{lem}\label{asdf} Let $0<r<n$ and $c=\gcd(r, n)$. Let $a_i=\begin{cases} 1, & \text{ if } 1 \le i \le r \\ 0, & \text{ if } r<i \le n \end{cases}$ and let $\Om_r=t^{[a_1, \cdots, a_n]} w_{S-\{r\}} w_S$ be a length $0$-element in $\tW^!_A$. Then there exists $\t \in S_n$ such that $\Om_r (1 2 \cdots c)=\t \tw_{n, r} \t \i$. 
\end{lem}

We consider the action of $w_{I-\{r\}} w_I$ on $\{1, \cdots, n\}$. By definition, $w_{I-\{r\}} w_I(i)=j$ if and only if $j \equiv i+r \mod n$. Hence the orbits of $w_{I-\{r\}} w_I$ on $\{1, \cdots, n\}$ are given by the $\mod c$-equivalence classes. Therefore $w_{I-\{r\}} w_I (1 2 \cdots c)$ acts transitively on $\{1, 2, \cdots, n\}$. 

Assume that $\ua_1(\Om_r (1 2 \cdots c))=[b_1, b_2, \cdots, b_n, b_1, b_2, \cdots, b_n, \cdots]$. Then $b_i \in \{0, 1\}$ and $b_i=1$ if and only if $1-(i-1) r \equiv j \mod n$ for some $1 \le j \le r$. In other words, $b_i=\lceil \frac{1-(i-1) r}{n} \rceil-\lceil \frac{1-i r}{n} \rceil$. 

We show that 

(a) $\lceil \frac{1-i r}{n} \rceil+\lceil \frac{i r}{n} \rceil=1$ for any $i \in \ZZ$.

Assume that $i r=x n+y$ for some $x \in \ZZ$ and $0 \le y<n$. Then $\lceil \frac{1-i r}{n} \rceil=\lceil \frac{1-y}{n} \rceil-x$ and $\lceil \frac{i r}{n} \rceil=\lceil \frac{y}{n} \rceil+x$. Notice that if $y=0$, then $\lceil \frac{1-y}{n} \rceil=1$ and $\lceil \frac{y}{n} \rceil=0$. If $1<y<n$, then $\lceil \frac{1-y}{n} \rceil=0$ and $\lceil \frac{y}{n} \rceil=1$. Therefore, $\lceil \frac{1-i r}{n} \rceil+\lceil \frac{i r}{n} \rceil=1$. (a) is proved. 

Hence $b_i=\lceil \frac{i r}{n} \rceil-\lceil \frac{(i-1) r}{n} \rceil$ for all $i$. So $$\ua_1(\Om_r (1 2 \cdots c))=\ua_n(\tw_{n, r})=[\chi_{n, r}, \chi_{n, r}, \cdots].$$ 

Now define $\t(i)=\bigl(w_{I-\{r\}} w_I (1 2 \cdots c) \bigr)^i(1)$. Then we have that $\Om_r (1 2 \cdots c)=\t \tw_{n, r} \t \i$.  



\subsection{Proof of Proposition \ref{stable}}\label{pstable} By $\S$\ref{ps} (2), it suffices to prove the case where $0 \le m<n$. If $m=0$, then $\tw_{n, 0}=(1 2 \cdots n)$ and $P_\th(\tw_{n, 0})=(1 2 \cdots n)$ for any $\th \in \ZZ^\infty$. Now suppose that $0<m<n$. Let $d=\gcd(n, m)$. Then $\chi_{n, m}=(\chi_{\frac{n}{d}, \frac{m}{d}}, \cdots, \chi_{\frac{n}{d}, \frac{m}{d}})$. It is easy to see that $P_\th(\tw_{n, m})=t^{[\chi, \chi, \cdots, \chi]} (1 2 \cdots n)$, where $\chi \in \ZZ^{\frac{n}{d}}$ with $P_\th(\tw_{\frac{n}{d}, \frac{m}{d}})=t^\chi (1 2 \cdots \frac{n}{d})$. Therefore it suffice to prove the case where $n$ and $m$ are relatively prime. 

We regard $\tw_{n, m}$ as an element in $\tW^!_A$. By Lemma \ref{asdf}, there exists $\t \in S_n$ with $\t \tw_{n, m} \t \i=\Om_m$. Since $\Om_m$ is of length $0$, it is a terminal element. By Proposition \ref{red}, there exists $\t' \in S_n$ such that $l(\t' P_\th(\t \tw_{n, m} \t \i) (\t') \i)=l(\t \tw_{n, m} \t \i)=0$. 

Hence $\t' P_\th(\t \tw_{n, m} \t \i) (\t') \i=\t' \t P_\th(\tw_{n, m}) \t \i (\t') \i$ is of length $0$. Therefore $\t' \t P_\th(\tw_{n, m}) \t \i (\t') \i=\Om_m$. Set $\s'=\t \i \t' \t \in S_n$. Then $P_\th(\tw_{n, m})=\s' \tw_{n, m} (\s') \i$. In particular, $(1 2 \cdots n)=\s' (1 2 \cdots n) (\s') \i$. So $\s'=(1 2 \cdots n)^i$ for some $i$. 

\subsection{} 
In the rest of this section, we show that given any quasi-positive element with $\s \in (1 2 \cdots n) (\ZZ/2 \ZZ)^n$, one can get $\tw_{n, m}$ or $\tw'_{n, m}$ by applying successive $P$-operator. 

We introduce the notation $\lto$. 

Let $\tw=t^{[a_1, \cdots, a_n]} \s$ be a quasi-positive element with $a_i \in \ZZ$ and $\s \in S'_n$. We say that $\th \in \ZZ^\infty$ is {\it extreme} for $\tw$ if its leading entry equals $\max\{a_1, \cdots, a_n, 1\}$. We write $\tw \lto \tw'$ if there exits finite sequences $\tw=\tw_0, \tw_1, \cdots, \tw_k=\tw'$ and $\g_1, \cdots, \g_k \in \ZZ^n$ such that $\g_i$ is extreme for $\tw_{i-1}$ and $\tw_i=P_{\g_i}(\tw_{i-1})$ for all $i$. 

\begin{prop}\label{redred}
Let $\tw=t^{[a_1, \cdots, a_n]} \s$ be a quasi-positive element with $a_i \in \ZZ$ and $\s \in (1 2 \cdots n) (\ZZ/2 \ZZ)^n \subset S'_n$. Then $\tw \lto \tw'_{n, 0}$ or $\tw \lto \tw'_{n, m}$ with $m \in \NN$ with $m \mid n$ or $\tw \lto (1 2 \cdots n)^i \tw_{n, m} (1 2 \cdots n)^{-i}$ for some $m \ge 0$ and $i \in \ZZ$. If moreover, $\s=(1 2 \cdots n)$, then $\tw \lto (1 2 \cdots n)^i \tw_{n, m} (1 2 \cdots n)^{-i}$ for some $m \ge 0$ and $i \in \ZZ$.
\end{prop}

The proof will be given in $\S$\ref{pfredred}. 

\begin{lem}\label{xyz}
For any quasi-positive element $\tw=t^{[a_1, \cdots, a_n]} \s$, set $$ev_0(\tw)=\sum_i a_i \text{ and } ev_1(\tw)=(x, y),$$ where $x=\max_i a_i$ and $y=\sharp\{i; a_i=x\}$. Then for any $\tw'$ with $\tw \lto \tw'$, $ev_0(\tw') \le ev_0(\tw)$ and $ev_1(\tw') \le ev_1(\tw)$ for the lexicographic order on $\ZZ \times \NN$. 
\end{lem}

It suffices to prove the case where $\tw'=P_\th(\tw)$ for some $\th$ extreme for $\tw$. If $x=0$, then $\tw'=\tw$ and the lemma is obvious. Now we assume that $x \ge 1$.

Since $\tw$ is quasi-positive, we have that $a_1, \cdots, a_n \ge 0$. Suppose that $e_\th(\tw)=[c_1, \cdots, c_n]$. Set $a'_i=a_i-c_i+c_{\s \i(i)}$. We show that 

(a) $|a'_i| \le a_i-c_i+c_{|\s \i(i)|}$.

If $c_i=1$, then $a_i \ge 1$. Thus $a_i-c_i \ge 0$ and $a'_i \ge -1$ for any $i$. Notice that $c_{\s \i(i)} \le |c_{\s \i(i)}|=c_{|\s \i(i)|}$. Thus if $a'_i \ge 0$, then $a'_i \le a_i-c_i+c_{|\s \i(i)|}$. If $a'_i<0$, then $a_i-c_i=0$ and $c_{\s \i(i)}=-1$. In this case, we still have that $|a'_i| \le a_i-c_i+c_{|\s \i(i)|}$. (a) is proved. 

Hence $ev_0(\tw')=\sum_i |a'_i| \le \sum_i (a_i-c_i+c_{|\s \i(i)|})=\sum_i a_i=ev_0(\tw)$. 

We show that 

(b) $|a'_i| \le x$ for all $i$. 

We have seen that $a'_i \ge -1$. Therefore $-a'_i \le x$ for all $i$. On the other hand, it is easy to see that $a'_i \le x+1$. If $a_i-c_i+c_{\s \i(i)}=x+1$, then $a_i=x$, $c_i=0$, $\s \i(i)>0$ and $c_{\s \i(i)}=1$. By $\S$\ref{se1} (4), $\ua_i(\tw)=(x, \ua_{\s \i(i)}(\tw)) \ge \ua_{\s \i(i)}(\tw) \ge \th$. Thus $c_i=1$, which is a contradiction. (b) is proved. 

If $a_i-c_i+c_{\s \i(i)}=x$, then $(a_i, c_i, c_{\s \i(i)})=(x, 1, 1)$ or $(x, 0, 0)$ or $(x-1, 0, 1)$. 

If $a_i-c_i+c_{\s \i(i)}=-x$, then $x=1$, $a_i-c_i=0$ and $\s \i(i)<0$ and $c_{|\s \i(i)|}=1$. Notice that $\ua_i(\tw)=(a_i, -\ua_{|\s \i(i)|}(\tw)) \ge (0, 0, \cdots)$. Thus $a_i=1$, $c_i=1$, $\s \i(i)<0$ and $c_{|\s \i(i)|}=1$. 

As a summary, \begin{gather*} \sharp\{i; a_i-c_i=x\}=y-\sharp\{i; c_i=1\}, \\ \tag{c} \sharp\{i; a_i-c_i \neq x, |a_i-c_i+c_{\s \i(i)}|=x\} \le \sharp\{i; c_{|\s \i(i)|}=1\}.\end{gather*} Hence $ev_1(\tw') \le (x, y)$ with equality holds if and only if the inequality (c) is an equality. 

\begin{lem}
Let $\tw=t^{[a_1, \cdots, a_n]} \s$ be a quasi-positive element with $\s \in (1 2 \cdots n) (\ZZ/2 \ZZ)^n \subset S'_n$. If $ev_1(\tw)=ev_1(\tw')$ for any $\tw'$ with $\tw \lto \tw'$, then 

(1) $\max_i a_i-\min_i a_i \le 1$. 

(2) If $\max_i a_i>1$, then $\s=(1 2 \cdots n)$. 
\end{lem}

Let $ev_1(\tw)=(a, b)$. If $a \le 1$, then the statement is obvious. Now assume that $a>1$. Set $\th=[a, -a-1, -a-1, \cdots] \in \ZZ^\infty$.  By assumption on $\tw$, $ev_1(P_\th(\tw))=ev_1(\tw)$. By (c) in the proof of Lemma \ref{xyz}, for any $i$ with $a_i=a$, we have that $\s(i)>0$ and $a_{\s(i)} \ge a-1$. Thus $\th$ is extreme for $P_\th(\tw)$. Since $P_\th^2(\tw)=(a, b)$, for any $i$ with $a_i=a$, we have that $\s^2(i)>0$ and $a_{\s^2(i)} \ge a-1$. One can show by induction on $j$ that for $j>0$, $ev_1((P_\th)^j(\tw))=(a, b)$ and for any $i$ with $a_i=a$, $(\s)^j(i)>0$ and $b_{(\s)^j(i)} \ge a-1$. Therefore $\s=(1 2 \cdots n)$ and $a-1 \le b_i \le a$ for all $i$.  

\begin{lem}\label{red1} Let $\tw=t^{[a_1, \cdots, a_n]} \s$ be a quasi-positive element with $\s \in (1 2 \cdots n) (\ZZ/2 \ZZ)^n \subset S'_n$. Then one of the following holds:

(1) $\tw \lto t^{[a'_1, \cdots, a'_n]} (1 2 \cdots n)$ for some $a'_1, \cdots, a'_n \ge 0$ with $\max_i a'_i-\min_i a'_i \le 1$;

(2) $\tw \lto \tw'_{n, m}$, where $m=0$ or $m \in \NN$ and $m \mid n$. 
\end{lem}

Notice that for any quasi-positive element $\tx$, $ev_1(\tx) \ge (0, n)$. Thus there exists an element $\tw'$ such that $\tw \lto \tw'$ and $ev_1(\tw') \le ev_1(\tw'')$ for any $\tw''$ with $\tw \lto \tw''$. By Lemma \ref{xyz}, $ev_1(\tw'')=ev_1(\tw')$ for any $\tw''$ with $\tw' \lto \tw''$. We may assume that $\tw'=t^{[b_1, \cdots, b_n]} \s'$ and $ev_1(\tw')=(a, b)$. By our assumption on $\s$, $\s' \in (1 2 \cdots n) (\ZZ/2 \ZZ)^n$. 

If $a=0$, then $\tw'=\s'$. Since $\tw'$ is quasi-positive, then $\s'=(1 2 \cdots n)$ or $ (n, -n) (1 2 \cdots n)$. If $a>1$, then by the previous lemma, $\max_i a'_i-\min_i a'_i \le 1$ and $\s'=(1 2 \cdots n)$. 

Now we consider the case that $a=1$ and $\s' \notin S_n$. 

For any element of the form $\tx=t^{[c_1, \cdots, c_n]} \t$ with $0 \le c_1, \cdots, c_n \le 1$ and $\t \in S'_n$, let $m(\tx)$ be the smallest positive integer $j$ such that there exists $i$ with $c_i=1$ and $c_{\t^{-j}(i)}=-1$. We show that 

(a) $m(\tw')$ is finite. 

Since $\s' \in (1 2 \cdots n) (\ZZ/2 \ZZ)^n$ and $a \neq 0$, $\ua_i(\tw') \neq (0, 0, \cdots)$ for all $i$. Since $\tw'$ is quasi-positive, $\ua_i(\tw')>(0, 0, \cdots)$ for all $i$. Let $1 \le i \le n$ with $(\s') \i(i)<0$. Notice that $\ua_i(\tw')=(b_i, -\ua_{|(\s') \i(i)|}(\tw'))>(0, 0, \cdots)$ and $\ua_{|(\s') \i(i)|}(\tw')>(0, 0, \cdots)$. Thus $b_i=1$. Let $j$ be the smallest positive integer such that $b_{(\s')^{-j}(i)} \neq 0$. For any $0<j'<j$, $b_{(\s')^{-j'}(i)}=0$. Hence $(\s') \i |(\s')^{-j'}(i)|>0$ for $0<j'<j$. Therefore $(\s')^{-j}(i)<0$ and $b_{(\s')^{-j}(i)}=-1$.  (a) is proved. 

Now assume that $\tw' \lto \tw''$ and $m(\tw'') \le m(\tw''')$ for any $\tw'' \lto \tw'''$. Suppose that $\tw''=t^{[c_1, \cdots, c_n]} \t$ and $m(\tw'')=j$. Let $k \in \{1, \cdots, n\}$ such that $c_k=1$ and for any $k'$ with $c_{k'}=1$, we have that $\ua_k(\tw'') \le \ua_{k'}(\tw'')$. Then $c_{\t^{-j}(k)}=-1$ and $c_{\t^{-j'}(k)}=0$ for any $0<j'<j$. 

Suppose that $\ua_{-\t^{-j}(k)}(\tw'') \neq \ua_k(\tw'')$. Then $ev_1(P_{\ua_{-\t^{-j}(k)}(\tw'')}(\tw''))<ev_1(\tw'')$ if $j=1$ and $m(P_{\ua_{-\t^{-j}(k)}(\tw'')}(\tw'')) \le j-1$ if $j>1$. This is a contradiction. Hence $\ua_{-\t^{-j}(k)}(\tw'')=\ua_k(\tw'')$. Now it is easy to see that $\ua_k(\tw'')=(\g, -\g, \g, -\g, \cdots)$, here $\g=(1, 0, 0, \cdots, 0) \in \ZZ^j$. Hence $j \mid n$ and $\tw''=(1 2 \cdots n)^k \tw'_{n, j} (1 2 \cdots n)^{-k}$. Let $\th=\ua_k(\tw'')$. Then $(P_\th)^{n-k}(\tw'')=\tw'_{n, n/j}$ and $\tw \lto \tw' \lto \tw'' \lto \tw'_{n, n/j}$. 

\subsection{}\label{exp} By the above Lemma, to prove Proposition \ref{redred}, we only need to consider the elements of the form $t^{[a_1, \cdots, a_n]} (1 2 \cdots n)$, where $a_i \ge 0$ and $\max_i a_i-\min_i a_i \le 1$. By $\S$\ref{ps} (2), it suffices to consider the case where $a_i \in \{0, 1\}$ for all $i$. We will relate an element of this form to an element of the form $t^{[b_1, \cdots, b_p]} (1 2 \cdots p)$ with $b_i \ge 0$ for some $p<n$. 

\subsection{} Set \begin{gather*} \aa^+_i=(0, 1, 1, \cdots, 1, 0) \in \ZZ^{i+3} \qquad \text{ if } i \ge 0; \\ \aa^-_i=(1, 0, 0, \cdots, 0, 1) \in \ZZ^{3-i} \qquad \text{ if } i \le 0. \end{gather*} 

We define a sign function on $\ZZ$ which is different from the standard one in order to simplify some notation we're going to use below. For $i \in \ZZ$, set $$\sgn'(i)=\sgn(i-1/2)=\begin{cases} +, & \text{ if } i>0; \\ -, & \text{ if } i \le 0. \end{cases}$$ We define the sign function in this way so that $\aa^{\sgn'(i)}_{i-1}$ is defined for all $i \in \ZZ$. 

Let $\tw=t^{[a_1, \cdots, a_n]} (1 2 \cdots n)$ with $a_i \in \{0, 1\}$ and both $0$ and $1$ appear in $\{a_1, \cdots, a_n\}$. This is equivalent to say that $a_i \in \{0, 1\}$ and $\sum_i a_i \notin \{0, n\}$. Define $ev_2(\tw)=(c, d)$, where \begin{gather*} c=\max\{l; \exists i, \text{ with } \ua_i(\tw)=(\aa^{\sgn'(l)}_l, *) \text{ for some } * \}, \\ d=\sharp\{1 \le i \le n; \ua_i(\tw) \text{ is of the form } (\aa^{\sgn'(c)}_c, *)\}.\end{gather*}

Notice that if $\th$ is extreme for $\tw$, then $P_\th(\tw)=t^{[b_1, \cdots, b_n]} (1 2 \cdots n)$ for some $b_i \ge 0$. By Lemma \ref{xyz}, $0 \le b_i \le 1$. By Lemma \ref{sn}, $\sum_i b_i=\sum_i a_i$. Therefore $ev_2(P_\th(\tw))$ is also defined.

\begin{lem}
We keep the notation as in the previous subsection. Let $\s=(1 2 \cdots n)$. Then 

(1) If $c>0$ and $\th \in \ZZ^\infty$ with the first $c+1$ entries are all $1$, then $ev_2(P_\th(\tw)) \le ev_2(\tw)$. Moreover, the equality holds if and only if for any $i$ with $\ua_i(\tw) \ge \th$, either $\ua_{\s^{c+2}(i)}(\tw)$ is of the form $(\aa^+_c, *)$ or $\ua_{\s^{c+1}(i)}(\tw)$ is of the form $(\aa^+_{c-1}, *)$. 

(2) If $c \le 0$ and $\th \in \ZZ^\infty$ is of the form $[\aa^-_c, *]$, then $ev_2(P_\th(\tw)) \le ev_2(\tw)$. Moreover, the equality holds if and only if for any $i$ with $\ua_i(\tw) \ge \th$, either $\ua_{\s^{1-c}(i)}(\tw)$ is of the form $(\aa^-_c, *)$ or $\ua_{\s^{2-c}(i)}(\tw)$ is of the form $(\aa^-_{c-1}, *)$. 
\end{lem} 

We only prove part (1). Part (2) can be proved in the same way. Let $e_\th(\tw)=[c_1, \cdots, c_n]$. We show that 

(a) $(c_i, c_{\s \i(i)}) \neq (1, 1)$ for any $i$.

Otherwise, $a_i=1$ and $\ua_{\s \i(i)}(\tw) \ge \th$. Hence $\ua_i(\tw)=(1, \ua_{\s \i(i)}(\tw))$ has the first $c+2$ entries all equal $1$. That contradicts the assumption that $ev_2(\tw)=(c, d)$. (a) is proved. 

Let $\g=(c_n, c_{n-1}, \cdots, c_1)$. Then $$\ua_i(P_\th(\tw))=\ua_i(\tw)-(\g[n-i], \g[n-i], \cdots)+(\g[n-i+1], \g[n-i+1], \cdots).$$ We show that 

(b) for any $i$, if $\ua_i(P_\th(\tw))$ is of the form $(\aa^+_{c'}, *)$, then $c' \le c$.

By the proof of Lemma \ref{xyz}, if $a_j-c_j+c_{\s \i(j)}=1$, then $(a_j, c_j, c_{\s \i(j)})=(1, 0, 0)$ or $(0, 0, 1)$. Therefore if $a_j-c_j+c_{\s \i(j)}=1$ and $a_j=0$, then $a_{\s \i(j)}-c_{\s \i(j)}+c_{\s^{-2}(j)}=0$. Now suppose that $\ua_i(P_\th(\tw))$ is of the form $[\aa^+_{c'}, *]$ with $c'>c$. Then $a_{\s^{-k}(i)}-c_{\s^{-k}(i)}+c_{\s^{-k-1}(i)}=1$ for all $0<k \le c'+1$. Therefore $(a_{\s^{-k}(i)}, c_{\s^{-k}(i)}, c_{\s^{-k-1}(i)})=(1, 0, 0)$ for all $0<k \le c'$. Since $ev_2(\tw)=(c, d)$, we must have that $c'=c+1$ and $a_{\s \i(i)}=\cdots=a_{\s^{-c'}(i)}=1$. Moreover, $a_{\s^{-c'-1}(i)}=0$ and $c_{\s^{-c'-2}(i)}=1$. Let $\g'=(1, 1, \cdots, 1, 0) \in \ZZ^{c+2}$. Thus $\ua_{\s \i(i)}(\tw)=(\g', \ua_{\s^{-c'-2}(i)}(\tw))$. By $\S$\ref{se1} (4), we have that $\ua_{\s \i(i)}(\tw) \ge \ua_{\s^{-c'-2}(i)}(\tw) \ge \th$. So $c_{\s \i(i)}=1$ and $a_{\s \i(i)}-c_{\s \i(i)}+c_{\s^{-2}(i)}=0$. That is a contradiction. (b) is proved.  

Let $m=\sharp\{1 \le i \le n; \ua_i(\tw) \ge \th\}$. We have that $$\sharp\{i; \ua_i(\tw)-(\g[n-i], \g[n-i], \cdots) \text{ is of the form } (\aa^+_c, *)\}=d-m.$$ By the proof of (b), if $ \ua_i(\tw)-(\g[n-i], \g[n-i], \cdots)$ is not of the form $(\aa^+_c, *)$ and $\ua_i(P_\th(\tw))$ is of the form $(\aa^+_c, *)$, then $\ua_i(\tw)-(\g[n-i], \g[n-i], \cdots)$ is of the form $(\aa^+_{c-1}, *)$ and $\ua_{\s^{-c-1}(i)}(\tw) \ge \th$. Therefore $$\sharp\{i; \ua_i(P_\th(\tw)) \text{ is of the form } (\aa^+_c, *)\} \le d-m+m=d$$ and the equality holds if and only if for any $i$ with $\ua_i(\tw) \ge \th$, either $\ua_{\s^{c+2}(i)}(\tw)$ is of the form $(\aa^+_c, *)$ or $\ua_{\s^{c+1}(i)}(\tw)$ is of the form $(\aa^+_{c-1}, *)$. 

\subsection{} Let $\sum_{j=1}^k i_j=n$. Let $a_0, a_1, \cdots, a_n \in \ZZ$ with $a_0=a_n$. Set $$\g_j=[a_{\sum_{l=1}^{j-1} i_l}, a_{\sum_{l=1}^{j-1} i_l+1}, \cdots, a_{\sum_{l=1}^{j-1} i_l+i_j}] \in \ZZ^{i_j+1}.$$ Then define $$\wedge(\g_1, \cdots, \g_k)=[a_0, \cdots, a_{n-1}].$$ In other words, $\wedge(\g_1, \cdots, \g_k)$ is defined if and only if the tail entry of $\g_j$ is the head entry of $\g_{j+1}$ for all $0<j<k$ and the tail entry of $\g_k$ is the head entry of $\g_1$.

Given $\g_j \in \ZZ^{i_j}$ for all $j \in \NN$ such that the tail entry of $\g_j$ is the head entry of $\g_{j+1}$ for all $j \in \NN$, we define $\wedge(\g_1, \g_2, \cdots) \in \ZZ^\infty$ in the same way. 

\begin{lem}\label{red22}
Let $\tw=t^{[a_1, \cdots, a_n]} (1 2 \cdots n)$ with $a_i \in \{0, 1\}$ for all $i$ and both $0$ and $1$ appear in $\{a_1, \cdots, a_n\}$. Then there exists $m, c \in \ZZ$ and $\chi \in \ZZ^n$ of the form $\wedge(\g_1, \cdots, \g_p)$, where $\g_i \in \{\aa^{\sgn'(c)}_c, \aa^{\sgn'(c)}_{c-1}\}$ such that $\tw \lto t^{\chi[m]}(1 2 \cdots n)$. 
\end{lem}

By definition, $ev_2(\tw')>(-n, 1)$ for any $\tw \lto \tw'$. We assume that $\tw \lto \tw'$ and $ev_2(\tw') \le ev_2(\tw'')$ for any $\tw' \lto \tw''$. By the previous lemma, $ev_2(\tw')=ev_2(\tw'')$ for any $\tw' \lto \tw''$. Assume that $ev_2(\tw')=(c, d)$. We only prove the case where $c>0$. The case where $c \le 0$ can be proved in the same way. 

Set $\th=(\aa^+_c[1], 0, 0, \cdots)$. Then $ev_2(P_\th(\tw'))=(c, d)$. Suppose that $\tw'=t^{\chi'} (1 2 \cdots n)$. By the previous lemma, for any $i$ with $\ua_i(\tw') \ge \th$, $\ua_{\s(i)}(\tw')$ is of the form $\wedge(\aa^+_c, *, \aa^+_c)$ or $\wedge(\aa^+_c, *, \aa^+_{c-1})$. Applying the same argument to $ev_2(P_\th^2(\tw'))$, we have that for any $i$ with $\ua_i(\tw') \ge \th$, $\ua_{\s(i)}(\tw')$ is of the form $\wedge(\aa^+_c, *, \aa^+_{l'}, \aa^+_l)$ for $l' \in \{c, c-1\}$. The lemma follows by repeating the same argument several times. 

\begin{lem}\label{uni}
Let $n \in \NN$ and $m \in \ZZ$. Then there exists $\chi \in \ZZ^n$ with $|\chi|=m$ such that for any $\chi' \in \ZZ^n$ with $|\chi'|=m$, $$t^{\chi'} (1 2 \cdots n) \lto t^{\chi[i]} (1 2 \cdots n)=(1 2 \cdots n)^i \bigl(t^\chi (1 2 \cdots n)\bigr) (1 2 \cdots n)^{-i}$$ for some $i$. 
\end{lem}

We prove the lemma by induction on $n$. The case where $n=1$ is trivial. Now let us assume that the statement is true when $n$ is replaced by any integer $n'<n$.  By $\S$\ref{ps} (2), it suffices to consider the case where $0 \le m<n$. If $m=0$, the lemma follows from Lemma \ref{red1}. Now let us assume that $0<m<n$. 

Notice that if $\chi' \in \ZZ^n$ is of the form $\wedge(\g_1, \cdots, \g_p)$, where $\g_i \in \{\aa^{\sgn'(c)}_c, \aa^{\sgn'(c)}_{c-1}\}$ and $|\chi'|=m$. Then we must have that $$p=\min\{m, n-m\},  \qquad c=\begin{cases} \lceil \frac{n}{p} \rceil-2, & \text{ if } m>\frac{n}{2} \\ \lceil -\frac{n}{p} \rceil+2, & \text{ if } m \le \frac{n}{2} \end{cases}$$ and $\sharp\{i; g_i=\aa^{\sgn'(c)}_c\}=d$. Here $d=\begin{cases} n-(c+1)p, & \text{ if } m>\frac{n}{2} \\ (-c+3)p-n, & \text{ if } m \le \frac{n}{2} \end{cases}$. In particular, $p, c, d$ are uniquely determined by $n$ and $m$ and $p<n$. Set \begin{gather*} \Xi=\{\chi \in \ZZ^n; \chi \text{ is of the form } \wedge(\g_1, \cdots, \g_p), \g_i \in \{\aa^{\sgn'(c)}_c, \aa^{\sgn'(c)}_{c-1}\}\}; \\ \Xi'=\{\chi \in \ZZ^\infty; \chi \text{ is of the form } \wedge(\g_1, \g_2, \cdots), \g_i \in \{\aa^{\sgn'(c)}_c, \aa^{\sgn'(c)}_{c-1}\}\}; \\ \Pi=\{\pi=[a_1, \cdots, a_p] \in \ZZ^p; a_i \in \{c, c-1\}, |\pi|=(c-1)p+d\}; \\ \Pi'=\{\pi=[a_1, a_2, \cdots] \in \ZZ^\infty; a_i \in \{c, c-1\}\}. \end{gather*} 

Define $f: \Pi \to \Xi$ by $[a_1, \cdots, a_p] \mapsto \wedge(\aa^{\sgn'(c)}_{a_1}, \cdots, \aa^{\sgn'(c)}_{a_p})$ and $f': \Pi' \to \Xi'$ by $[a_1, a_2] \mapsto \wedge(\aa^{\sgn'(c)}_{a_1}, \aa^{\sgn'(c)}_{a_2}, \cdots)$. Then $f$ and $f'$ are well-defined bijective maps. For $\pi \in \Pi$ and $\th \in \Pi'$, $\pi \ge \th$ if and only if $f(\pi) \ge f'(\th)$. 

Let $\th \in \Pi'$ and $\pi \in \Pi$. Set $\th_1=\begin{cases} f'(\th)[1], & \text{ if } c>0 \\ f'(\th), & \text{ if } c \le 0 \end{cases}$. Then $\th_1$ is extreme for $f(\pi)$. Now assume that $P_\th(t^\pi (1 2 \cdots p))=t^{\pi'} (1 2 \cdots p)$ and $P_{\th_1}(t^{f(\pi)} (1 2 \cdots n))=t^{\p''} (1 2 \cdots n)$. Then it is easy to see that $\pi''=f(\pi')$ or $f(\pi')[1]$. Therefore for $\pi, \pi' \in \Pi$, if $t^{\pi'} (1 2 \cdots p) \lto t^{\pi} (1 2 \cdots p)$, then $t^{f(\pi')} (1 2 \cdots n)\lto t^{f(\pi)[i]} (1 2 \cdots n)$ for some $i$. 

By induction hypothesis on $\Pi$, there exists $\pi \in \Pi$ such that for any $\pi' \in \Pi$, $t^{\pi'} (1 2 \cdots p) \lto t^{\pi[i]} (1 2 \cdots p)$ for some $i$. Thus for any $\chi' \in \Xi$, $t^{\chi'} (1 2 \cdots n) \lto t^{f(\pi)[j]} (1 2 \cdots n)$ for some $j$. 

Let $\chi' \in \ZZ^n$ with $|\chi'|=m$. By Lemma \ref{red1}, there exists $\chi'_1 \in \ZZ^n$ with $|\chi'_1|=m$ and the entries of $\chi'_1$ are $0$ or $1$, such that $t^{\chi'} (1 2 \cdots n) \lto t^{\chi'_1} (1 2 \cdots n)$. By Lemma \ref{red22}, there exists $\chi'_2 \in \Xi$ and $i \in \ZZ$ with $t^{\chi'_1} (1 2 \cdots n) \lto t^{\chi'_2[i]}(1 2 \cdots n)$.  Therefore $$t^{\chi'} (1 2 \cdots n) \lto t^{f(\pi)[i][j]} (1 2 \cdots n)=t^{f(\pi)[i+j]} (1 2 \cdots n)$$ for some $i, j \in \ZZ$. 

\subsection{Proof of Proposition \ref{redred}}\label{pfredred} By Lemma \ref{red1}, one only needs to consider the element of the form  $t^{[a_1, \cdots, a_n]} (1 2 \cdots n)$ for some $a_1, \cdots, a_n \ge 0$ with $\max_i a_i-\min_i a_i \le 1$. By $\S$\ref{exp}, it suffices to consider the case where $a_i \in \{0, 1\}$. Let $m=\sum_i a_i$. By Lemma \ref{uni}, there exists $\chi \in \ZZ^n$ with $|\chi|=m$ such that for any $\chi' \in \ZZ^n$ with $|\chi'|=m$, $t^{\chi'} (1 2 \cdots n) \lto t^{\chi[i]} (1 2 \cdots n)$ for some $i$. Taking $\chi'=\chi'_{n, m}$, then by Proposition \ref{stable}, $\chi=\chi'_{n, m}[i]$ for some $i \in \ZZ$. Taking $\chi'=[a_1, \cdots, a_n]$, then there exists $j \in \ZZ$ such that \begin{align*} t^{\chi'} (1 2 \cdots n) & \lto t^{\chi'_{n, m}[i][j]} (1 2 \cdots n)=t^{\chi'_{n, m}[i+j]} (1 2 \cdots n) \\ &=(1 2 \cdots n)^{i+j} \tw_{n, m} (1 2 \cdots n)^{-i-j}.\end{align*}

This finishes the proof. 

\section{The main theorem}

\subsection{}\label{dd} Let $\tilde \l=[(b_1, c_1), \cdots, (b_k, c_k)]$ and $\tilde \mu=[(b_{k+1}, c_{k+1}), \cdots, (b_l, c_l)]$ with $(\tilde \l, \tilde \mu) \in \cd \cp$. Set $$\tw^{st}_{(\tilde \l, \tilde \mu)}=t^{[\chi'_{b_1, c_1}, \cdots, \chi'_{b_l, c_l}]} w_{(\l, \mu)}.$$ Here $\chi'$ is defined in $\S$\ref{chi}. We call $\tw^{st}_{(\tilde \l, \tilde \mu)}$ the {\it standard} element associated to $(\tilde \l, \tilde \mu)$. 

We show that 

\begin{thm}\label{61}
(1) Let $\tW^!=\tW^!_A$ and $\tw \in [(\tilde \l, \tilde \O)]$ for $(\tilde \l, \tilde \O) \in \cd \cp$. Then there exists $\tw' \in S_n \cdot \tw^{st}_{(\tilde \l, \tilde \mu)}$ such that $\tw \tilde \to \tw'$. 

(2) Let $\tW^!$ be of type B, C or D and $\tw \in [(\tilde \l, \tilde \mu)]'$ for $(\tilde \l, \tilde \mu) \in \cd \cp_{\ge 0}$. Then there exists $\tw' \in S_n \cdot \tw^{st}_{(\tilde \l, \tilde \mu)}$ such that $\tw \tilde \to \tw'$. 
\end{thm}

The proof here is long and technical. We will first prove part (1), which is relatively easier. 

\subsection{} Notice that $t^{[a, a, \cdots, a]} \in \tW^!$ is of length $0$ and commutes with any elements in $\tW^!$. Hence for any $\tx, \tx' \in \tW^!$, $\tx \tilde \to \tx'$ if and only if $t^{[a, \cdots, a]} \tx \tilde \to t^{[a, \cdots, a]} \tx'$. Also for $\tilde \l=[(b_1, c_1), \cdots, (b_k, c_k)]$, we have that $t^{[a, \cdots, a]} \tw^{st}_{(\tilde \l, \O)}=\tw^{st}_{(\tilde \l', \O)}$, where $\tilde \l'=[(b_1, c_1+a b_1), \cdots, (b_k, c_k+a b_k)]$. Replacing $\tw$ by $t^{[a, \cdots, a]} \tw$ for some $a \in \NN$ if necessary, we may assume that $\tw$ is quasi-positive. 

By definition, it suffices to prove the case where $\tw$ is a terminal element in $\tW^!$. We assume that $\tilde \l=[(b_1, c_1), \cdots, (b_k, c_k)]$. Since $\tw$ is quasi-positive, $c_i \ge 0$. By definition, $\tw \in S_n \cdot t^{[a_1, \cdots, a_n]} w_{(\l, \O)}$ for $a_i \ge 0$ with $\sum_{j=1}^{b_{l'}} a_{b_1+\cdots+b_{l'-1}+j}=c_{l'}$ for $1 \le l' \le k$. 

By Proposition \ref{redred}, there exists a finite sequence $\th_1, \cdots, \th_m \in \ZZ^\infty$ such that $$P_{\th_1} \cdots P_{\th_m}(t^{[a_1, \cdots, a_n]} w_{(\l, \O)})=t^{[a'_1, \cdots, a'_n]} w_{(\l, \O)},$$ where $[a'_1, \cdots, a'_{b_1}]=\chi'_{b_1, c_1}[i]$ for some $i$ and there exists a finite sequence $\th_m, \cdots, \th_{m'}$ such that $$P_{\th_m} \cdots P_{\th_{m'}} (t^{[a'_1, \cdots, a'_n]} w_{(\l, \O)})=t^{[a''_1, \cdots, a''_n]} w_{(\l, \O)},$$ where $[a''_{b_1+1}, \cdots, a''_{b_1+b_2}]=\chi'_{b_2, c_2}[j]$ for some $j$. Moreover, by Proposition \ref{stable}, $[a''_1, \cdots, a''_{b_1}]=\chi'_{b_1, c_1}[i']$ for some $i'$.  Repeating the procedure for each pair $(b_k, c_k)$, one can show that there exists a finite sequence $\th_1, \cdots, \th_p \in \ZZ^\infty$ such that $$P_{\th_1} \cdots P_{\th_p} (t^{[a_1, \cdots, a_n]} w_{(\l, \O)})=t^{[f_1, \cdots, f_n]} w_{(\l, \O)},$$ where for $1 \le l \le k$, $[f_{b_1+\cdots+b_{l-1}+1}, \cdots, f_{b_1+\cdots+b_{l-1}+b_l}]=\chi'_{b_l, c_l}[i_l]$ for some $i_l$. Now set $\s=(1 2 \cdots b_1)^{-i_1} \cdots (n-b_k+1, n-b_k+2, \cdots, n)^{-i_k}$. Then $t^{[f_1, \cdots, f_n]} w_{(\l, \O)}=\s \tw^{st}_{(\tilde \l, \tilde \O)} \s \i$. 

By Proposition \ref{red}, there exists $\tw' \in S_n \cdot P_{\th_1} \cdots P_{\th_p} (t^{[a_1, \cdots, a_n]} w_{(\l, \mu)})=S_n \cdot \tw^{st}_{(\tilde \l, \O)}$ such that $\tw \tilde \to \tw'$. This finishes the proof of part (1). 

\subsection{} Now we consider the case where $\tW^!$ is of type BCD. Again it suffices to prove the case where $\tw$ is a terminal element in $\tW^!$. By definition, $ev_0(\tx) \ge 0$ for any $\tx \in \tW^!_{int}$. Then there exists $\tw'$ with $\tw \tilde \to \tw'$ and $ev_0(\tw') \le ev_0(\tw'')$ for all $\tw''$ with $\tw \tilde \to \tw''$. By Lemma \ref{xyz}, $ev_0(\tw'')=ev_0(\tw')$ for all $\tw''$ with $\tw' \tilde \to \tw''$. 

By Corollary \ref{xxxxx}, there exists $\t \in S'_n$ such that $\tw \tilde \to \t \tw' \t \i$ and $\t \tw' \t \i$ is quasi-positive. Since $\tw$ is terminal, $\tw'$ and $\t \tw' \t \i$ are also terminal. Replacing $\tw$ by $\t \tw' \t \i$, we only need to consider the case that $\tw$ is quasi-positive and $ev_0(\tx)=ev_0(\tw)$ for all $\tx$ with $\tw \tilde \to \tx$. 

Let $(\tilde \l', \tilde \mu', \tilde \g')$ such that $\tilde \l'=[(b_1, c_1), \cdots, (b_k, c_k)]$ is positive, $\tilde \mu'=[(b_{k+1}, c_{k+1}), \cdots, (b_l, c_l)]$ is special, $\tilde \g'=[(b_{l+1}, c_{l+1}), \cdots, (b_m, c_m)]$ with $\sum_{i=1}^m b_i=n$ and for $l<i \le m$, $c_i \ge 2$ and $c_i \mid b_i$. We associate an element $$\tw_{(\tilde \l', \tilde \mu', \tilde \g')}=t^{[\chi'_{b_1, c_1}, \cdots, \chi'_{b_m, c_m}]} \s (1 2 \cdots b_1) \cdots (n-b_m+1 \cdots n),$$ where $\s \in (\ZZ/2 \ZZ)^n \subset S'_n$ is defined by $\s(j)=-j$ if and only if $j=b_1+\cdots+b_i$ for some $i>k$ or $j=b_1+\cdots+b_{i-1}+r \frac{b_i}{c_i}$ for some $i>l$ and $0<r<c_i$. It is easy to see that $\tw_{(\tilde \l', \tilde \mu', \tilde \g')}$ is quasi-positive. 

By Proposition \ref{redred}, one can show in the same way as we did for type A above that there exist a triple $(\tilde \l', \tilde \mu', \tilde \g')$ and $\t' \in S'_n$ such that $\tx \tilde \approx \t' \tw_{(\tilde \l', \tilde \mu', \tilde \g')} (\t') \i$. By Corollary \ref{xxxxx}, there exists $\t \in S_n$ such that $\t' \tw_{(\tilde \l', \tilde \mu', \tilde \g')} (\t') \i \tilde \approx \t \tw_{(\tilde \l', \tilde \mu', \tilde \g')} \t \i$, $\t \tw_{(\tilde \l', \tilde \mu', \tilde \g')} \t \i$ is quasi-positive and \break $\ua_1(\t \tw_{(\tilde \l', \tilde \mu', \tilde \g')} \t \i) \ge \cdots \ge \ua_n(\t \tw_{(\tilde \l', \tilde \mu', \tilde \g')} \t \i) \ge (0, 0, \cdots)$. By our assumption on $\tx$, $$ev_0(\tx)=ev_0(\t \tw_{(\tilde \l', \tilde \mu', \tilde \g')} \t \i)=ev_0(\tw_{(\tilde \l', \tilde \mu', \tilde \g')})=\sum_{i=1}^m c_i.$$ 

It remains to show that $\tilde \g'$ is empty. We will first illustrate the idea through an example and then give the  rigorous proof. 

\subsection{} Let $\tilde \l'=[(1, 2)], \tilde \mu'=[(2, 1), (2, 1)], \tilde \g'=[(4, 2)]$. Then $\tw_{(\tilde \l', \tilde \mu', \tilde \g')}=t^{[2, 0, 1, 0, 1, 0, 1, 0, 1]} (2, -3, -2, 3) (4, -5, -4, 5) (6, -7, -8, 9)$ and $ev_0(\tw_{(\tilde \l', \tilde \mu', \tilde \g')})=6$. Let $\t=(4 9) (2 6 7 3)$ and $\tx'=\t \tw_{(\tilde \l', \tilde \mu', \tilde \g')} \t \i$. Then $$\tx'=t^{[2, 1, 1, 1, 1, 0, 0, 0, 0]}(2, 6, -2, -6) (3, 7, -4, -8) (5, 9, -5, -9)$$ is quasi-positive and lies in $\cw$.

Let $\o=[1, 1, 1, 0, 0, 0, 0, 0, 0]$. Then $\o$ is minuscule for $\tx'$ and by the same argument as we did in the proof of Proposition \ref{red}, there exists $\tx'' \in S'_n \cdot t^{-\o} \tx' t^\o$ such that $\tx' \tilde \approx \tx''$. By definition, $$t^{-\o} \tx' t^\o=t^{[2, 0, 0, 1, 1, 1, 1, 0, 0]} (2, 6, -2, -6) (3, 7, -4, -8) (5, 9, -5, -9)$$ and the quasi-positive element in the $(\ZZ/2 \ZZ)^n$-conjugacy class of $t^{-\o} \tx' t^\o$ is $$\tx_1=t^{[2, 0, 0, 1, 1, 1, 1, 0, 0]} (2, -6, -2, 6) (3, -7, 4, 8) (5, 9, -5, -9).$$

Let $\th=\ua_7(\tx_1)$. Then $$P_\th(\tx_1)=t^{[2, 0, 0, 0, 1, 1, 0, 0, 0]} (2, -6, -2, 6) (3, 7, 4, 8) (5, 9, -5, -9)$$ and $ev_0(P_\th(\tx_1))=4<ev_0(\tx')$. 

\subsection{}\label{diffi} Now we continue our proof. Suppose that $\tilde \g$ is not empty. Let $i=\min\{\t(j); j>b_1+\cdots+b_l\}$ and $\th=\ua_i(\t \tw_{(\tilde \l', \tilde \mu', \tilde \g')} \t \i)$. Then $\th=\max_{j>b_1+\cdots+b_l} \ua_j(\tw_{(\tilde \l', \tilde \mu', \tilde \g')})$. In particular, the leading entry of $\th$ is positive. 

Let $\o=[d_1, \cdots, d_n]$, where $d_j=\begin{cases} 1, & \text{ if } j \le i \\ 0, & \text{ if } j>i \end{cases}$. By the same argument as we did in the proof of Proposition \ref{red}, $\o$ is minuscule for $\t \tw_{(\tilde \l', \tilde \mu', \tilde \g')} \t \i$ and there exists $x \in S'_n$ such that $$\t \tw_{(\tilde \l', \tilde \mu', \tilde \g')} \t \i \tilde \approx x t^{-\o} \t \tw_{(\tilde \l', \tilde \mu', \tilde \g')} \t \i t^\o x \i \in S'_n \cdot (t^{-\t \i \o} \tw_{(\tilde \l', \tilde \mu', \tilde \g')} t^{\t \i \o}).$$ 

We will calculate the element $t^{-\t \i \o} \tw_{(\tilde \l', \tilde \mu', \tilde \g')} t^{\t \i \o}$.

Choose $\th' \in \ZZ^\infty$ with $\th'>\th$ and $\th'<\ua_j(\tw_{(\tilde \l', \tilde \mu', \tilde \g')})$ for all $j$ such that $\ua_j(\tw_{(\tilde \l', \tilde \mu', \tilde \g')})>\th$. Let $l \le m'<m$ with $b_1+\cdots+b_{m'}<\t \i(i) \le b_1+\cdots+b_{m'+1}$. 

Assume that $e_\th(\tw_{(\tilde \l', \tilde \mu', \tilde \g')})=[e_1, \cdots, e_n]$ and $e_{\th'}(\tw_{(\tilde \l', \tilde \mu', \tilde \g')})=[e'_1, \cdots, e'_n]$. Then it is easy to see that for $j \le b_1+\cdots+b_l$ with $\ua_j(\tw_{(\tilde \l', \tilde \mu', \tilde \g')})=\th$, then $j=b_1+\cdots+b_{l'}$ for $k<l' \le l$. Therefore, 

(i) for $0 \le l' \le k-1$, $[d_{\t(b_1+\cdots+b_{l'}+1)}, d_{\t(b_1+\cdots+b_{l'}+2)}, \cdots, d_{\t(b_1+\cdots+b_{l'+1})}]=[c_{b_1+\cdots+b_{l'}+1}, \cdots, c_{b_1+\cdots+b_{l'+1}}]=[c'_{b_1+\cdots+b_{l'}+1}, \cdots, c'_{b_1+\cdots+b_{l'+1}}]$.

(ii) for $k \le l' \le l-1$ with $\ua_{b_1+\cdots+b_{l'+1}}(\tw_{(\tilde \l', \tilde \mu', \tilde \g')}) \neq \th$ or $\t(b_1+\cdots+b_{l'+1})>i$, $[d_{\t(b_1+\cdots+b_{l'}+1)}, \cdots, d_{\t(b_1+\cdots+b_{l'+1})}]=[c'_{b_1+\cdots+b_{l'}+1}, \cdots, c'_{b_1+\cdots+b_{l'+1}}]$.

(iii) for $k \le l' \le l-1$ with $\ua_{b_1+\cdots+b_{l'+1}}(\tw_{(\tilde \l', \tilde \mu', \tilde \g')})=\th$ and $\t(b_1+\cdots+b_{l'+1})<i$, $[d_{\t(b_1+\cdots+b_{l'}+1)}, \cdots, d_{\t(b_1+\cdots+b_{l'+1})}]=[c_{b_1+\cdots+b_{l'}+1}, \cdots, c_{b_1+\cdots+b_{l'+1}}]$.

(iv) for $l' \ge l$ with $l' \neq m'$, $[d_{\t(b_1+\cdots+b_{l'}+1)}, \cdots, d_{\t(b_1+\cdots+b_{l'+1})}]=[c'_{b_1+\cdots+b_{l'}+1}, \cdots, c'_{b_1+\cdots+b_{l'+1}}]$.

(v) for $l'=m'$ and $1 \le j \le b_{l'+1}$ with $\t(b_1+\cdots+b_{l'}+j) \neq i$, we have that $d_{\t(b_1+\cdots+b_{l'}+j)}=c'_{b_1+\cdots+b_{l'}+j}$. We have that $d_i=1$ and $c'_{\t \i(i)}=0$. 

Assume that \begin{gather*} t^{-e_\th(\tw_{(\tilde \l', \tilde \mu', \tilde \g')})} \tw_{(\tilde \l', \tilde \mu', \tilde \g')} t^{e_\th(\tw_{(\tilde \l', \tilde \mu', \tilde \g')})}=t^{[a_1, \cdots, a_n]} \s (1 2 \cdots b_1) \cdots (n-b_m+1 \cdots n), \\ t^{-e_{\th'}(\tw_{(\tilde \l', \tilde \mu', \tilde \g')})} \tw_{(\tilde \l', \tilde \mu', \tilde \g')} t^{e_{\th'}(\tw_{(\tilde \l', \tilde \mu', \tilde \g')})}=t^{[a'_1, \cdots, a'_n]} \s (1 2 \cdots b_1) \cdots (n-b_m+1 \cdots n), \\ t^{-\t \i \o} \tw_{(\tilde \l', \tilde \mu', \tilde \g')} t^{\t \i \o}=t^{[a''_1, \cdots, a''_n]} \s (1 2 \cdots b_1) \cdots (n-b_m+1 \cdots n).\end{gather*} By (i)--(iv) above, for any $l' \neq m'$, $[a''_{b_1+\cdots+b_{l'}+1}, \cdots, a''_{b_1+\cdots+b_{l'+1}}]$ equals $[a_{b_1+\cdots+b_{l'}+1}, \cdots, a_{b_1+\cdots+b_{l'+1}}]$ or $[a'_{b_1+\cdots+b_{l'}+1}, \cdots, a'_{b_1+\cdots+b_{l'+1}}]$. 
By Proposition \ref{stable} and Proposition \ref{stable1}, $t^{-\t \i \o} \tw_{(\tilde \l', \tilde \mu', \tilde \g')} t^{\t \i \o}$ is conjugated by $S'_n$ to the element $t^{[a'''_1, \cdots, a'''_n]} \s (1 2 \cdots b_1) \cdots (n-b_m+1 \cdots n)$, where $$[a'''_{b_1+\cdots+b_{l'}+1}, \cdots, a'''_{b_1+\cdots+b_{l'+1}}]=\begin{cases} \chi'_{b_{l'+1}, c_{l'+1}}, & \text{ if } l' \neq m'; \\ [a''_{b_1+\cdots+b_{l'}+1}, \cdots, a''_{b_1+\cdots+b_{l'+1}}], & \text{ if } l'=m'. \end{cases}$$

If $b_{m'+1}=c_{m'+1}$, then $\th=\ua_i(\t \tw_{(\tilde \l', \tilde \mu', \tilde \g')} \t \i)=(1, -1, 1, -1, \cdots)$ and by direct calcuation $|a''_{b_1+\cdots+b_{m'}+1}|+\cdots+|a''_{b_1+\cdots+b_{m'+1}}|=c_{m'+1}-2$. Thus $$ev_0(t^{-\t \i \o} \tw_{(\tilde \l', \tilde \mu', \tilde \g')} t^{\t \i \o})=\sum_{j \neq m'} c_j+(c_{m'}-2)=ev_0(\tw_{(\tilde \l', \tilde \mu', \tilde \g')})-2.$$ This contradicts our assumption on $\tx$. Therefore $\frac{b_{m'}}{c_{m'}} \ge 2$. By direct calculation, $$t^{[a'''_1, \cdots, a'''_n]} \s (1 2 \cdots b_1) \cdots (n-b_m+1 \cdots n) \in S'_n \cdot \tx',$$ where $\tx'=t^{[\chi'_{b_1, c_1}, \cdots, \chi''_{b_{m'+1}, c_{m'+1}}, \cdots, \chi'_{b_m, c_m}]}  \s' (1 2 \cdots b_1) \cdots (n-b_m+1 \cdots n)$. Here $\chi''_{b_{m'+1}, c_{m'+1}}=[\chi'_{\frac{b_{m'+1}}{c_{m'+1}}-1, 1}, \chi'_{\frac{b_{m'+1}}{c_{m'+1}}, 1}, \chi'_{\frac{b_{m'+1}}{c_{m'+1}}, 1}, \cdots, \chi'_{\frac{b_{m'+1}}{c_{m'+1}}, 1}, \chi'_{\frac{b_{m'+1}}{c_{m'+1}}+1, 1}]$ and $\s' \in (\ZZ/2 \ZZ)^n \subset W(C_n)$ is defined by $\s'(j)=-j$ if and only if one of the following holds:

(i) $j=b_1+\cdots+b_p$ for some $p>k$;

(ii) $j=b_1+\cdots+b_{p}+r \frac{b_{p+1}}{c_{p+1}}$ for some $p \ge l$ and $p \neq m'$ and $0<r<c_{p+1}$;

(iii) $j=b_1+\cdots+b_{m'}+r \frac{b_{m'+1}}{c_{m'+1}}-1$ for some $0<r<c_{m'+1}$. 

Then $\tx'$ is quasi-positive. By Corollary \ref{xxxxx}, there exists $\t'' \in S_n$ such that $\tx \tilde \approx \t'' \tx' (\t'') \i$. Let $\th''=\ua_{b_1+\cdots+b_{m'+1}}(\tx')$. Assume that $(P'_{\th''})^{\frac{b_{m'+1}}{c_{m'+1}}-1}(\tx')=t^{[f_1, \cdots, f_n]} w$ for $f_i \in \ZZ$ and $w \in S'_n$. By Proposition \ref{stable} and Proposition \ref{stable1}, $\sum_{j \le b_1+\cdots+b_{m'}} f_j=\sum_{j \le m'} c_j$ and $\sum_{j>b_1+\cdots+b_{m'+1}} f_j=\sum_{j>m'+1} c_j$. 

By direct calculation, $\sum_{b_1+\cdots+b_{m'}<j \le b_1+\cdots+b_{m'+1}} f_j=c_{m'+1}-2$. Therefore, $ev_0((P'_{\th''})^{\frac{b_{m'+1}}{c_{m'+1}}-1}(\tx'))<ev_0(\tx)$. This contradicts our assumption on $\tx$. This finishes the proof for part (2).

\

Now we state our main theorem. 

\begin{thm}\label{thA}
Let $\tW^!$ be of classical type and $\co$ be an integral conjugacy class of $\tW^!$. Then 

(1) For any $\tw \in \co$, there exists $\tw' \in \co_{\min}$ with $\tw \tilde \to \tw'$.

(2) For any $\tw, \tw' \in \co_{\min}$, $\tw \tilde \sim \tw'$. 
\end{thm}

The proof will be given in $\S$\ref{pfthA}. 

\begin{lem}\label{cd1}
Let $(\tilde \l, \tilde \mu) \in \cd \cp_{\ge 0}$ and $\o=[\frac{1}{2}, \cdots, \frac{1}{2}]$. Then for any $\tw \in S_n \cdot \tw^{st}_{(\tilde \l, \tilde \mu)}$, $t^{-\o} \tw t^\o \in S'_n \cdot \tw^{st}_{(\tilde \l, \underline{\tilde \mu})}$. 
\end{lem}

Let $\s \in S_n$ with $\tw=\s \tw^{st}_{(\tilde \l, \tilde \mu)} \s \i$. Then $t^{-\o} \tw t^\o=\s t^{-\o} \tw^{st}_{(\tilde \l, \tilde \mu)} t^\o \s \i \in S_n \cdot t^{-\o} \tw^{st}_{(\tilde \l, \tilde \mu)} t^\o$. By direct calculation, \begin{align*} t^{-\o_n} \tw^{st}_{(\tilde \l, \tilde \mu)} t^{\o_n} &=t^{[\chi'_{b_1, c_1}, \cdots, \chi'_{b_k, c_k}, -\chi'_{b_{k+1}, 1-c_{k+1}}, \cdots, -\chi'_{b_l, 1-c_l}]} w_{(\l, \mu)} \\ &=\t t^{[\chi'_{b_1, c_1}, \cdots, \chi'_{b_k, c_k}, \chi'_{b_{k+1}, 1-c_{k+1}}, \cdots, \chi'_{b_l, 1-c_l}]} w_{(\l, \mu)} \t \i, \end{align*} where $\t \in (\ZZ/2 \ZZ)^n \subset W(C_n)$ is defined by $\t(j)=-j$ if and only if $j>b_1+\cdots+b_k$. 

\begin{lem}\label{cd2}
Let $\tW^!$ be of type C or D and $(\tilde \l, \tilde \mu), (\tilde \l, \underline{\tilde \mu}) \in \cd \cp_{\ge 0}$. Then 

(1) for any $\tw \in S'_n \cdot \tw^{st}_{(\tilde \l, \tilde \mu)}$, there exists $\tw' \in S_n \cdot \tw^{st}_{(\tilde \l, \underline{\tilde \mu})}$ with $\tw \tilde \to \tw'$. 

(2) for any minimal length element $\tw$ in $S'_n \cdot \tw^{st}_{(\tilde \l, \tilde \mu)}$, there exists a minimal length element $\tw'$ in $S_n \cdot \tw^{st}_{(\tilde \l, \underline{\tilde \mu})}$ with $\tw \tilde \approx \tw'$. 
\end{lem}

By Corollary \ref{xxxxx}, $\tw \tilde \to \tw_1$ for some quasi-positive element $\tw_1 \in S_n \cdot \tw^{st}_{(\tilde \l, \tilde \mu)}$. Let $\o=[\frac{1}{2}, \cdots, \frac{1}{2}]$. Then $\o$ is minuscule and there exists $w \in W$ such that $l(t^\o w)=0$. By Lemma \ref{cd1}, $\tw_1 \tilde \approx (t^\o w) \i \tw_1 (t^\o w)=w \i t^{-\o} \tw_1 t^\o w \in S'_n \cdot \tw^{st}_{(\tilde \l, \underline{\tilde \mu})}$. So $\tw \tilde \to \tw'$ for some $\tw' \in S'_n \cdot \tw^{st}_{(\tilde \l, \underline{\tilde \mu})}$. By Corollary \ref{xxxxx}, we may choose $\tw'$ in $S_n \cdot \tw^{st}_{(\tilde \l, \underline{\tilde \mu})}$. Part (1) is proved. 

Now assume that $\tw$ is of minimal length in $S'_n \cdot \tw^{st}_{(\tilde \l, \tilde \mu)}$. By (1), $\tw \tilde \to \tw'$ for some $\tw' \in S_n \cdot \tw^{st}_{(\tilde \l, \underline{\tilde \mu})}$. Again by (1), there exits $\tw'' \in S_n \cdot \tw^{st}_{(\tilde \l, \tilde \mu)}$ with $\tw' \tilde \to \tw''$. Hence $l(\tw'') \le l(\tw') \le l(\tw)$. Since $\tw$ is of minimal length in $S'_n \cdot \tw^{st}_{(\tilde \l, \tilde \mu)}$, $l(\tw'') \ge l(\tw)=l(\tw')$. Thus $l(\tw)=l(\tw')=l(\tw'')$ and $\tw \tilde \approx \tw'$.  

For any $\tx \in S'_n \cdot \tw^{st}_{(\tilde \l, \underline{\tilde \mu})}$, there exists $\tx' \in S'_n \cdot \tw^{st}_{(\tilde \l, \tilde \mu)}$ with $\tx \tilde \to \tx'$. Hence $l(\tx) \ge l(\tx') \ge l(\tw)$. Therefore $\tw'$ is of minimal length in $S'_n \cdot \tw^{st}_{(\tilde \l, \underline{\tilde \mu})}$. Part (2) is proved. 

\subsection{Proof of Theorem \ref{thA}}\label{pfthA} Let $\co_{ter}$ the set of terminal elements in $\co$. We show that 

(a) for any $\tw, \tw' \in \co_{ter}$, $\tw \tilde \sim \tw'$. 

If $\tW$ is of type $A$ or $B$, then there exists a unique element $(\tilde \l, \tilde \mu) \in \cd \cp$ that represents $\co$ in the sense of $\S$\ref{psc}. By Theorem \ref{61}, there exists $\tw_1, \tw'_1 \in S_n \cdot \tw^{st}_{(\tilde \l, \tilde \mu)}$ such that $\tw \tilde \to \tw_1$ and $\tw' \tilde \to \tw'_1$. 

If $\tW^!$ is of type $\tilde C_n$ or $\tilde D_n$, then $\co=[(\tilde \l, \tilde \mu)]' \cup [(\tilde \l, \underline{\tilde \mu})]'$ for some $(\tilde \l, \tilde \mu) \sim (\tilde \l, \underline{\tilde \mu}) \in \cd \cp_{\ge 0}$. By Theorem \ref{61}, there exists $\tw_2 \in S_n \cdot \tw^{st}_{(\tilde \l, \tilde \mu)} \cup S_n \cdot \tw^{st}_{(\tilde \l, \underline{\tilde \mu})}$ with $\tw \tilde \to \tw_2$. If $\tw_2 \in S_n \cdot \tw^{st}_{(\tilde \l, \tilde \mu)}$, then we set $\tw_1=\tw_2$. If $\tw_2 \in S_n \cdot \tw^{st}_{(\tilde \l, \underline{\tilde \mu})}$, then by Lemma \ref{cd2}, there exists $\tw_1 \in S_n \cdot \tw^{st}_{(\tilde \l, \tilde \mu)}$ such that $\tw \tilde \to \tw_2 \tilde \to \tw_1$. Similarly, there exists $\tw'_1 \in S_n \cdot \tw^{st}_{(\tilde \l, \tilde \mu)}$ such that $\tw' \tilde \to \tw'_1$.

Since $\tw$ and $\tw'$ are terminal elements, $\tw_1, \tw'_1$ are also terminal elements and $\tw \tilde \approx \tw_1$ and $\tw' \tilde \approx \tw'_1$. By Proposition \ref{par}, $\tw_1$ and $\tw'_1$ are minimal length elements in $W \cdot \tw^{st}_{(\tilde \l, \tilde \mu)}$ and $\tw' \tilde \sim \tw'_1$. Hence $\tw \tilde \sim \tw'$. 

(a) is proved. 

It is easy to see that minimal length element in $\co$ are terminal. In particular, part (2) follows from (a). 

By (a), all the terminal elements have the same length. Hence an element in $\co$ is a terminal element if and only if it is a minimal length element.

For any $\tw \in \co$, by definition there exists a terminal element $\tw'$ with $\tw \tilde \to \tw'$. We have that $\tw'$ is a minimal length element in $\co$. Part (1) is proved. 

\begin{cor}\label{Bruhat1}
Let $\tW$ be a classical extended affine Weyl group and $\co$ be an integral conjugacy class in $\tW$. Then $$\co_{\min}=\{\tw \in \co; \tw \text{ is a minimal element for the Bruhat order in } \co\}.$$
\end{cor}

Let $\tw \in \co$. If $\tw$ is a minimal length element, then $\tw$ is a minimal element for the Bruhat order in $\co$. 
If $\tw$ is not a minimal length element, then by \cite[Lemma 4.4]{He072} and Theorem \ref{thA}, there exists $\tw' \in \co_{\min}$ with $\tw'<\tw$ for the Bruhat order on $\co$. 

\

In the rest of this section, we consider $\tW$ instead of $\tW^!$. 

\begin{thm}\label{wa}
Let $\tW$ be a classical extended affine Weyl group and $\co$ be a conjugacy class of $\tW$. We assume furthermore that $\co \subset W_a$. Then 

(1) For any $\tw \in \co$, there exists $\tw' \in \co_{\min}$ such that $\tw \tilde \to \tw'$. 

(2) For any $\tw, \tw' \in \co_{\min}$, $\tw \tilde \approx \tw'$.
\end{thm}

If $\tW$ is of type A, then we have a natural surjective group homomorphism $p: \tW^!_A \to \tW$ defined in $\S$\ref{tww}. This map is $W_a$-equivalent for the conjugation action of $W_a$ and is length-preserving. For conjugacy class $\co$ of $\tW$, there exists a conjugacy class $\co^!$ of $\tW^!$ such that $p: \co^! \to \co$ is a bijection. Now part (1) and (2) for $\co$ follows from Theorem \ref{thA} for $\co^!$. 

If $\tW$ is of type B or C, then $\tW=\tW^!$ and $\co \subset \tW^!_{int}$. This case has been proved in Theorem \ref{thA}.

If $\tW$ is of type D, then $\co, \iota \co \iota \i \subset \tW^!_{int}$. Set $\co^!=\co \cup \iota \co \iota \i$. Then $\co^!$ is an integral conjugacy class of $\tW^!$. 

We first consider the case that $\co \cap \iota \co \iota \i=\emptyset$. 

By Theorem \ref{thA}, for any $\tw \in \co$, there exists $\tw_1 \in \co^!_{\min}$ such that $\tw \tilde \to \tw_1$ with resepct to $\tW^!$. In other words, $\tw \tilde \to \tw_1$ or $\tw \tilde \to \iota \tw_1 \iota \i$ with respect to $\tW$. Here the difference is that in $\tW^!$ we allow conjugation by $\iota$, but in $\tW$ we don't.  Notice that only one of $\tw_1, \iota \tw_1 \iota \i$ is contained in $\co$. This element is a minimal length element in $\co$. Part (1) is proved in this case. 

If $\tw, \tw' \in \co_{\min}$, then by Theorem \ref{thA}, $\tw \tilde \approx \tw'$ with respect to $\tW^!$. Hence $\tw \tilde \approx \tw'$ or $\tw \tilde \approx \iota \tw' \iota \i$ with respect to $\tW$. However, the latter case can't happen since $\tw, \iota \tw' \iota \i$ lies in different conjugacy classes of $\tW$. Part (2) is proved in this case.

Now we assume that $\co=\iota \co \iota \i$. Let $(\tilde \l, \tilde \mu) \in \cd \cp_{\ge 0}$ with $\co=[(\tilde \l, \tilde \mu)]' \cup [(\tilde \l, \underline{\tilde \mu})]'$. By $\S$\ref{psc}, $\tilde \mu \neq \tilde \O$ or some entry of $\tilde \l$ is of the form $(b, 0)$ with $b$ odd. Replacing $\tilde \mu$ by $\underline{\tilde \mu}$ if necessary, we may assume that some entry of $\tilde \mu$ is of the form $(*, 0)$ if $\tilde \mu \neq \tilde \O$. Hence for any $\tx \in (S_n \cdot \tw^{st}_{(\tilde \l, \tilde \mu)}) \cap \cw$, $\ua_n(\tx)=(0, 0, \cdots)$ and $\iota \tx \iota \i \in W \cdot \tx$. 

By Theorem \ref{61} and Proposition \ref{par}, for any $\tw \in \co$, there exists $\tw_1 \in (S_n \cdot \tw^{st}_{(\tilde \l, \tilde \mu)}) \cap \cw$ such that $\tw \tilde \to \tw_1$ or $\tw \tilde \to \iota \tw_1 \iota \i$ with respect to $\tW$. In other words, $\tw \tilde \to \tw'$ for some $\tw' \in W \cdot \tw^{st}_{(\tilde \l, \tilde \mu)}$. Now part (1) and (2) follows from Proposition \ref{par}. 

\

The following special case will be used in \cite{GH}. 

\begin{prop}
Let $\tW$ be a classical extended affine Weyl group and $\tw=t^\chi w$ with $w$ a Coxeter element in $W$.  

(1) If $\tw \in W_a$, then $\tw \tilde \to w$. 

(2) If $\tW=\tW(A_{n-1})$, $0<r<n$ and $\t_r=t^{\o_r} w_{S-\{r\}} w_S$ is a length $0$ element in $\tW$ with $\tw \in \t_r W_a$, then $\tw \tilde \to \t_r (1 2 \cdots \gcd(n, r))$. 
\end{prop}

We first consider the case that $\tW=\tW(A_{n-1})$. Then there exists $\tw' \in [(\tilde \l, \tilde \O)]$ with $\tilde \l=(n, r)$ for some $0 \le r<n$ such that $p(\tw')=\tw$. Here $p: \tW^! \to \tW$ is defined in $\S$\ref{tww}. By Theorem \ref{61}, there exists $\tw'_1 \in S_n \cdot \tw^{st}_{(\tilde \l, \tilde \O)}$ such that $\tw' \tilde \to \tw'_1$ with respect to $\tW^!$. 

If $\tw \in W_a$, then $r=0$. In this case, $\tw^{st}_{(\tilde \l, \tilde \O)}$ is a Coxeter element in $S_n$. Hence $\tw' \tilde \to w'$ with respect to $\tW$ for some $w'$ in the conjugacy class of $W$ that contains $w$. By \cite[Proposition 3.1.6 \& Theorem 3.2.7]{GP00}, $\tw' \tilde \to w' \to w$. 

If $\tw \notin W_a$, then $r \neq 0$ and $\tw \in \t_r W_a$. Let $c=\gcd(n, r)$. By Lemma \ref{asdf}, $\t_r (1 2 \cdots c) \in S_n \cdot \tw^{st}_{(\tilde \l, \tilde \O)}$. Notice that $\t_r \in {}^S \tW$ and $I(\t_r)=S-\{c, 2c, \cdots, n-c\}=\{1, 2, \cdots, c-1\} \sqcup \{c+1, c+2, \cdots, 2 c-1\} \sqcup \cdots \sqcup \{n-c+1, n-c+2, \cdots, n-1\}$. By Proposition \ref{par}, there exists $x \in W_{I(\t_r)}$ such that $\tw'_1 \to \t_r x$. Since the action of $\t_r$ on $I(\t_r)$ sends $i$ to $i+r$, each orbit is of the form $\{i, i+c, \cdots, i+n-c\}$ for some $0<i<c$. By \cite[Lemma 2.7]{GKP}, there exists $x_1 \in W_{\{1, 2, \cdots, c\}}$ such that $\t_r x \to \t_r x_1$ and $x_1$ is conjugated to $(1 2 \cdots c)$ by $W_{\{1, 2, \cdots, c\}}$. By \cite[Proposition 3.1.6 \& Theorem 3.2.7]{GP00}, $\t_r x_1 \approx \t_r (1  2 \cdots c)$ and $\tw' \tilde \to \tw'_1 \to \t_r (1 2 \cdots c)$. 

Now we consider the case where $\tW$ is of type B, C or D. Set $\tilde \mu=(n, 0)$ if $\tW$ is of type B or C and $\tilde \mu=[(n-1, 0), (1, 0)]$ if $\tW$ is of type D. Then $\tw \in [(\tilde \O, \tilde \mu)]'$ if $\tW$ is of type B and $\tw \in [(\tilde \O, \tilde \mu)]' \cup [(\tilde \O, \underline{\tilde \mu})]'$ if $\tW$ is of type C or D. By Theorem \ref{61}, $\tw \tilde \to w'$ with respect to $\tW^!$ for some $w' \in S_n \cdot \tw^{st}_{(\tilde \O, \tilde \mu)}$.  

Notice that in either case, $\tw^{st}_{(\tilde \O, \tilde \mu)}=w_{(\O, \mu)} \in W \cdot w$ and $\iota (S_n \cdot \tw^{st}_{(\tilde \O, \tilde \mu)}) \iota \i \subset W \cdot w$ if $W$ is of type D. Thus $\tw \tilde \to w''$ with respect to $\tW$ for some $w'' \in W \cdot w$. By \cite[Proposition 3.1.6 \& Theorem 3.2.7]{GP00}, $w'' \to w$. The Proposition is proved. 

\section{Distinguished conjugacy class}

\subsection{}\label{zzzz} Let $\tilde \l=[(b_1, c_1), \cdots, (b_k, c_k)]$ and $\tilde \mu=[(b_{k+1}, c_{k+1}), \cdots, (b_l, c_l)]$ with $(\tilde \l, \tilde \mu) \in \cd \cp$. We say that $(\tilde \l, \tilde \mu)$ is {\it distinguished} if $b_i$ and $c_i$ are coprime for all $i \le k$ and $c_i=1$ for $i>k$. 

Let $(\tilde \l, \tilde \mu) \in \cd \cp$ be distinguished. Then for $1 \le l' \le l$ and $1 \le j<j' \le b_{l'}$, $\ua_{b_1+\cdots+b_{l'-1}+j}(\tw^{st}_{(\tilde \l, \tilde \mu)}) \neq \ua_{b_1+\cdots+b_{l'-1}+j'}(\tw^{st}_{(\tilde \l, \tilde \mu)})$. Moreover, if $\ua_j(\tw^{st}_{(\tilde \l, \tilde \mu)})=\ua_{j'}(\tw^{st}_{(\tilde \l, \tilde \mu)})$ for some $1 \le j<j' \le n$, then there exists $1 \le m<m' \le l$ such that 

(1) either $m, m' \le k$ or $m, m'>k$.

(2) $(b_m, c_m)=(b_{m'}, c_{m'})$ and $1 \le j'' \le b_m$ such that $j=b_1+\cdots+b_{m-1}+j''$ and $j'=b_1+\cdots+b_{m'-1}+j''$. 

Define $\t \in S_n$ in such a way that $\t(i)<\t(j)$ if $\ua_i(\tw^{st}_{(\tilde \l, \tilde \mu)})>\ua_j(\tw^{st}_{(\tilde \l, \tilde \mu)})$ or $\ua_i(\tw^{st}_{(\tilde \l, \tilde \mu)})=\ua_j(\tw^{st}_{(\tilde \l, \tilde \mu)})$ and $i<j$. Set $\tw^f_{(\tilde \l, \tilde \mu)}=\t \tw^{st}_{(\tilde \l, \tilde \mu)} \t \i$. We call $\tw^f_{(\tilde \l, \tilde \mu)}$ the {\it fundamental element} corresponding to $(\tilde \l, \tilde \mu)$. 

It is easy to see that 

(1) If $(\tilde \l, \tilde \O) \in \cd \cp$ is distinguished, then $\tw^f_{(\tilde \l, \tilde \O)} \in {}^S \tW^!_A$.

(2) If $\tW^!$ is of type B, C or D and $(\tilde \l, \tilde \mu) \in \cd \cp_{\ge 0}$ is distinguished, then $\tw^f_{(\tilde \l, \tilde \mu)} \in {}^S \tW^!$. 

\subsection{} 
Let $\tW^!$ be of classical type. We call an integral conjugacy class in $\tW^!$ {\it distinguished} if for any $(\tilde \l, \tilde \mu) \in \cd \cp$ that represents this conjugacy class in the sense of $\S$\ref{psc}, $(\tilde \l, \tilde \mu)$ is distinguished. In other words, an integral conjugacy class in $\tW^!_A$ or $\tW^!_B$ is distinguished if it corresponds to distinguished $(\tilde \l, \tilde \mu)$. An integral conjugacy class in $\tW^!_C$ or $\tW^!_D$ is distinguished if it corresponds to $(\tilde \l, \tilde \mu) \sim (\tilde \l, \underline{\tilde \mu})$ such that $(\tilde \l, \tilde \mu)$ and $(\tilde \l, \underline{\tilde \mu})$ are distinguished. In this case, we must have that $\tilde \mu=\tilde \O$. 

The following two Theorems sharpen Theorem \ref{thA} for distinguished integral conjugacy classes. 

\begin{thm}\label{thB}
Let $\tW^!$ be of type A, B or C and $\co$ be a distinguished integral conjugacy class of $\tW^!$. We assume that $\co$ is represented by $(\tilde \l, \tilde \mu)$. Then 

(1) $\tw \to \tw^f_{(\tilde \l, \tilde \mu)}$ for any $\tw \in \co$. 

(2) $\tw \approx \tw^f_{(\tilde \l, \tilde \mu)}$ for any $\tw \in \co_{\min}$. 
\end{thm}

\begin{thm}\label{thC}
Let $\tW^!=\tW^!_D$ and $\co$ be a distinguished integral conjugacy class of $\tW^!$. We assume that $\co$ is represented by the double partition $(\tilde \l, \tilde \O)$. Then 

(1) $\tw \to \tw^f_{(\tilde \l, \tilde \O)}$ or $\tw \to \iota(\tw^f_{(\tilde \l, \tilde \O)})$ for any $\tw \in \co$. 

(2) $\tw \approx \tw^f_{(\tilde \l, \tilde \O)}$ or $\tw \approx \iota(\tw^f_{(\tilde \l, \tilde \O)})$ for any $\tw \in \co_{\min}$. 
\end{thm}

The proofs of these two theorems will be given in $\S$\ref{pfthB}. 

\begin{lem}\label{cd3}
Let $(\tilde \l, \tilde \mu) \in \cd \cp_{\ge 0}$ be distinguished and $\o \in \ZZ^n$ be of the form $[1, 1, \cdots, 1, 0, 0, \cdots, 0]$. Then $t^{-\o} \tw^f_{(\tilde \l, \tilde \mu)} t^\o \in S'_n \cdot \tw^f_{(\tilde \l, \tilde \mu)}$. If moreover, $\tilde \mu=\O$, then $t^{-\o} \tw^f_{(\tilde \l, \tilde \O)} t^\o \in S_n \cdot \tw^f_{(\tilde \l, \tilde \O)}$.
\end{lem}

Assume that the first $i$ entries of $\o$ are $1$ and the rest entries are $0$. Let $\th=\ua_i(\tw^f_{(\tilde \l, \tilde \mu)})$ and $\th' \in \ZZ^\infty$ such that $\th'>\th$ and $\th'<\ua_j(\tw^f_{(\tilde \l, \tilde \mu)})$ for all $j$ with $\ua_j(\tw^f_{(\tilde \l, \tilde \mu)})>\th$. Let $\s \in S_n$ with $\tw^f_{(\tilde \l, \tilde \mu)}=\s \tw^{st}_{(\tilde \l, \tilde \mu)} \s \i$. Assume that \begin{gather*} t^{-e_\th(\tw^{st}_{(\tilde \l, \tilde \mu)})} \tw^{st}_{(\tilde \l, \tilde \mu)} t^{e_\th(\tw^{st}_{(\tilde \l, \tilde \mu)})}=t^{[a_1, \cdots, a_n]} w_{\l, \mu}, \\ t^{-e_{\th'}(\tw^{st}_{(\tilde \l, \tilde \mu)})} \tw^{st}_{(\tilde \l, \tilde \mu)} t^{e_{\th'}(\tw^{st}_{(\tilde \l, \tilde \mu)})}=t^{[a'_1, \cdots, a'_n]} w_{\l, \mu}, \\ t^{-\s \i \o} \tw^{st}_{(\tilde \l, \tilde \mu)} t^{\s \i \o}=t^{[a''_1, \cdots, a''_n]} w_{\l, \mu}. \end{gather*} By the same argument as in $\S$\ref{diffi}, $[a''_{b_1+\cdots+b_{j-1}+1}, \cdots, a''_{b_1+\cdots+b_j}]$ equals $[a_{b_1+\cdots+b_{j-1}+1}, \cdots, a_{b_1+\cdots+b_j}]$ or $[a'_{b_1+\cdots+b_{j-1}+1}, \cdots, a'_{b_1+\cdots+b_j}]$ for any $1 \le j \le l$. By Proposition \ref{stable1} and \ref{stable}, $t^{-\s \i \o} \tw^{st}_{(\tilde \l, \tilde \mu)} t^{\s \i \o} \in S'_n \cdot \tw^{st}_{(\tilde \l, \tilde \mu)}$. Therefore $$t^{-\o} \tw^f_{(\tilde \l, \tilde \mu)} t^{\o}=\s t^{-\s \i \o} \tw^{st}_{(\tilde \l, \tilde \mu)} t^{\s \i \o} \s \i \in S'_n \cdot \tw^{st}_{(\tilde \l, \tilde \mu)}=S'_n \cdot \tw^f_{(\tilde \l, \tilde \mu)}.$$ 

If moreover, $\tilde \mu=\tilde \O$, then $t^{-\s \i \o} \tw^{st}_{(\tilde \l, \tilde \mu)} t^{\s \i \o} \in S_n \cdot \tw^{st}_{(\tilde \l, \tilde \mu)}$ and $t^{-\o} \tw^f_{(\tilde \l, \tilde \O)} t^\o \in S_n \cdot \tw^f_{(\tilde \l, \tilde \O)}$. 

\begin{lem}\label{ttt}
Let $\tW^!$ be of classical type and $(\tilde \l, \tilde \mu) \in \cd \cp$ be distinguished.  Let $\t \in \tW^!$ with $l(\t)=0$.

(1) If $\tW^!=\tW^!_A$ and $\tilde \mu=\tilde \O$, then $\t \i \tw^f_{(\tilde \l, \tilde \mu)} \t \in S_n \cdot \tw^f_{(\tilde \l, \tilde \mu)}$. 

(2) If $\tW^!=\tW^!_B$ and $(\tilde \l, \tilde \mu) \in \cd \cp_{\ge 0}$, then $\t \i \tw^f_{(\tilde \l, \tilde \mu)} \t \in S'_n \cdot \tw^f_{(\tilde \l, \tilde \mu)}$. 

(3) If $\tW^!$ is of type C or D, $\tilde \l$ is positive, $\tilde \mu=\tilde \O$ and $\t \neq \iota$, then $\t \i \tw^f_{(\tilde \l, \tilde \mu)} \t \in S'_n \cdot \tw^f_{(\tilde \l, \tilde \mu)}$. 
\end{lem}

If $\tW^!=\tW^!_A$ and $\tilde \mu=\tilde \O$, then $\t=t^{\o} w$ for some $\o \in \ZZ^n$ of the form $[1, 1, \cdots, 1, 0, 0, \cdots, 0]$ and $w \in S_n$. Replacing $\tw^f_{(\tilde \l, \tilde \O)}$ by $t^{[a, a, \cdots, a]} \tw^f_{(\tilde \l, \tilde \O)}$ for some $a \in \NN$ if necessary, we may assume that $(\tilde \l, \tilde \O) \in \cd \cp_{\ge 0}$. By Lemma \ref{cd3}, $\t \i \tw^f_{(\tilde \l, \tilde \mu)} \t \in S_n \cdot \tw^f_{(\tilde \l, \tilde \mu)}$. 

If $\tW^!=\tW^!_B$ and $(\tilde \l, \tilde \mu) \in \cd \cp_{\ge 0}$, then $\t=t^{\o} w$ for $\o=[1, 0, 0, \cdots, 0]$ and some $w \in S'_n$. By Lemma \ref{cd3}, $\t \i \tw^f_{(\tilde \l, \tilde \mu)} \t \in S'_n \cdot \tw^f_{(\tilde \l, \tilde \mu)}$. 

If $\tW^!=\tW^!_C$, $\tilde \l$ is positive and $\tilde \mu=\tilde \O$, then $\t=t^{\o} w$ for $\o=[\frac{1}{2}, \frac{1}{2}, \cdots, \frac{1}{2}]$ and some $w \in S'_n$. Since $\tilde \mu=\tilde \O$, $w_{(\l, \mu)} \in S_n$ and $\tw^f_{(\tilde \l, \tilde \mu)}=t^\chi \s$ for some $\chi \in \ZZ^n$ and $\s \in S_n$. By direct calculation, $t^{-\o} \tw^f_{(\tilde \l, \tilde \mu)} t^\o=\tw^f_{(\tilde \l, \tilde \mu)}$ and $\t \i \tw^f_{(\tilde \l, \tilde \mu)} \t \in S'_n \cdot \tw^f_{(\tilde \l, \tilde \mu)}$. 

If $\tW^!=\tW^!_D$, $\tilde \l$ is positive and $\tilde \mu=\tilde \O$, then $\t=t^{\o} w$ for $\o=[1, 0, 0, \cdots, 0]$ or $[\frac{1}{2}, \frac{1}{2}, \cdots, \frac{1}{2}, \pm \frac{1}{2}]$ and some $w \in S'_n$. If $\o=[1, 0, 0, \cdots, 0]$, then by Lemma \ref{cd3}, $\t \i \tw^f_{(\tilde \l, \tilde \mu)} \t \in S'_n \cdot \tw^f_{(\tilde \l, \tilde \mu)}$. If $\o=[\frac{1}{2}, \frac{1}{2}, \cdots, \frac{1}{2}]$, then $t^{-\o} \tw^f_{(\tilde \l, \tilde \mu)} t^\o=\tw^f_{(\tilde \l, \tilde \mu)}$ and $\t \i \tw^f_{(\tilde \l, \tilde \mu)} \t \in S'_n \cdot \tw^f_{(\tilde \l, \tilde \mu)}$. Set $\g=[1, 1, \cdots, 1, 0] \in \ZZ^n$ and $\o'=[\frac{1}{2}, \frac{1}{2}, \cdots, \frac{1}{2}]$. If $\o=[\frac{1}{2}, \frac{1}{2}, \cdots, \frac{1}{2}, -\frac{1}{2}]=\g-\o'$, then by Lemma \ref{cd3}, $$t^{-\o} \tw^f_{(\tilde \l, \tilde \mu)} t^{\o}=t^{-\g} t^{\o'} \tw^f_{(\tilde \l, \tilde \mu)} t^{-\o'} t^\g=t^{-\g} \tw^f_{(\tilde \l, \tilde \mu)} t^{\g} \in S'_n \cdot \tw^f_{(\tilde \l, \tilde \mu)}.$$ Hence $\t \i \tw^f_{(\tilde \l, \tilde \mu)} \t \in S'_n \cdot \tw^f_{(\tilde \l, \tilde \mu)}$. 

\subsection{Proof of Theorems \ref{thB} and \ref{thC}}\label{pfthB} By Theorem \ref{61}, for any $\tw \in \co$, $\tw \tilde \to \tw'$ for some $\tw' \in S_n \cdot \tw^{st}_{(\tilde \l, \tilde \mu)}=S_n \cdot \tw^f_{(\tilde \l, \tilde \mu)}$. Since $\tw^f_{(\tilde \l, \tilde \mu)} \in {}^S \tW$, by Proposition \ref{par}, $\tw' \tilde \to \tw^f_{(\tilde \l, \tilde \mu)}$. By definition, $\tw \to \t \tw^f_{(\tilde \l, \tilde \mu)} \t \i$ for some $\t \in \tW^!$ with $l(\t)=0$. 

By Lemma \ref{ttt} and Proposition \ref{par}, if $\tW^!$ is of type ABC or $\tW^!=\tW^!_D$ and $\t \neq \iota$, then $\t \tw^f_{(\tilde \l, \tilde \mu)} \t \i \to \tw^f_{(\tilde \l, \tilde \mu)}$. This proves part (1) of Theorems \ref{thB} and \ref{thC}. Part (2) follows from part (1) and the definition of $\approx$. 

\begin{cor}\label{por}
Let $\tW^!$ be of classical type and $\co$, $\co'$ be two distinguished integral conjugacy classes in $\tW$. Then the following conditions are equivalent:

(1) For some $w' \in \co'_{\min}$, there exists $w \in \co_{\min}$ such that $w \le w'$.

(2) For any $w' \in \co'_{\min}$, there exists $w \in \co_{\min}$ such that $w \le w'$.
\end{cor}

It is obvious that (2) implies (1). Now we rove that (1) implies (2). Let $w', w'_1 \in \co'_{\min}$ and $w \in \co_{\min}$ with $w \le w'$. By Theorems \ref{thB} and \ref{thC}, $w'_1 \tilde \approx w'$. Then by \cite[Lemma 4.4]{He072}, there exists $w_1 \in \co_{\min}$ with $w_1 \le w'_1$. 

\subsection{}\label{pw} Let $\tilde \l=[(b_1, c_1), \cdots, (b_k, c_k)]$ and $\tilde \mu=[(b_{k+1}, c_{k+1}), \cdots, (b_l, c_l)]$ with $(\tilde \l, \tilde \mu) \in \cd \cp$. Set $d_i=\gcd(b_i, c_i)$ and $m=\sum_{l'>k, c_{l'}=0} b_{l'}$. Let $\tilde \l'$ be the double partition whose entries are $(\frac{b_1}{d_1}, \frac{c_1}{d_1}), \cdots, (\frac{b_1}{d_1}, \frac{c_1}{d_1})$ (with $d_1$-times), $\cdots$, $(\frac{b_k}{d_k}, \frac{c_1}{d_k}), \cdots, (\frac{b_k}{d_k}, \frac{c_k}{d_k})$ (with $d_k$-times) and $(1, 0), (1, 0), \cdots, (1, 0)$ (with $m$-times) and $\tilde \mu'$ be the double partition whose entries are $(b_{l'}, c_{l'})$ with $l'>k$ and $c_{l'}=1$. Set $d(\tilde \l, \tilde \mu)=(\tilde \l', \tilde \mu')$. Then $d(\tilde \l, \tilde \mu)$ is distinguished and $w^f_{d(\tilde \l, \tilde \mu)}$ is defined. Set $d'(\tilde \l, \tilde \mu)=d(\tilde \l', \underline{\tilde \mu'})$. Then it is easy to see that  $d'(\tilde \l, \tilde \mu)$ is distinguished and is of the form $(*, \tilde \O)$. 

The following theorem relates an integral conjugacy class to a conjugacy class corresponding to a distinguished pair of double partitions. 

\begin{thm}\label{62}
(1) Let $\tW^!=\tW^!_A$ and $\tw \in [(\tilde \l, \tilde \O)]$ for some $(\tilde \l, \tilde \O) \in \cd \cp$. Then there exists $x \in  W_{I(\tw^f_{d(\tilde \l, \tilde \O)})}$ with $\tw \tilde \to x \tw^f_{d(\tilde \l, \tilde \O)}$. 

(2) Let $\tW^!=\tW^!_B$ and $\tw \in [(\tilde \l, \tilde \mu)]'$ for some $(\tilde \l, \tilde \mu) \in \cd \cp_{\ge 0}$. Then there exists $x \in W_{I(\tw^f_{d(\tilde \l, \tilde \mu)})}$ with $\tw \tilde \to x \tw^f_{d(\tilde \l, \tilde \mu)}$. 

(3) Let $\tW^!=\tW^!_C$ and $\tw \in [(\tilde \l, \tilde \mu)]' \cup [(\tilde \l, \underline{\tilde \mu})]'$ for some $(\tilde \l, \tilde \mu) \in \cd \cp_{\ge 0}$. Then there exists $x \in W_{I(\tw^f_{d(\tilde \l, \tilde \mu)})}$ with $\tw \tilde \to x \tw^f_{d(\tilde \l, \tilde \mu)}$. 

(4) Let $\tW^!=\tW^!_D$ and $\tw \in [(\tilde \l, \tilde \mu)]' \cup [(\tilde \l, \underline{\tilde \mu})]'$ for some $(\tilde \l, \tilde \mu) \in \cd \cp_{\ge 0}$. Then there exists $x \in W_{I(\tw^f_{d(\tilde \l, \tilde \mu)})}$ with $\tw \tilde \to x \tw^f_{d(\tilde \l, \tilde \mu)}$ or $\ua_n(\tw^f_{d(\tilde \l, \tilde \mu)})=(0, 0, \cdots)$ and $\tw \tilde \to x \tw^f_{d(\tilde \l, \tilde \mu)} (n, -n)$.
\end{thm}

If $\tW^!=\tW^!_A$, then replacing $\tw$ by $t^{[a, a, \cdots, a]} \tw$ for some $a \in \NN$, we may assume that $\tw \in [(\tilde \l, \tilde \O)]$ with $(\tilde \l, \tilde \O) \in \cd \cp_{\ge 0}$. If $\tW^!$ is of type C or D and $\tw \in [(\tilde \l, \underline{\tilde \mu})]'$, then there exists $w \in W$ such that $l(t^\o w)=0$, where $\o=[\frac{1}{2}, \cdots, \frac{1}{2}]$. Conjugation by $t^\o w$ preserve the length and sends $[(\tilde \l, \underline{\tilde \mu})]'$ to $[(\tilde \l, {\tilde \mu})]'$. Thus replacing $\tw$ by $t^\o w \tw (t^\o w) \i$ if necessary, we may assume that $\tw \in [(\tilde \l, \tilde \mu)]'$.  

As a summary, it suffices to consider the case where $\tw \in [(\tilde \l, \tilde \mu)]' \cap \tW^!$ with $(\tilde \l, \tilde \mu) \in \cd \cp_{\ge 0}$. 

Let $(\tilde \l', \tilde \mu')=d(\tilde \l, \tilde \mu)$ be the double partition defined in $\S$\ref{pw}. Then the multi-set $\{\ua_1(\tw^{st}_{(\tilde \l, \tilde \mu)}), \cdots, \ua_n(\tw^{st}_{(\tilde \l, \tilde \mu)})\}$ is the same as the multi-set $$\{\ua_1(\tw^{st}_{(\tilde \l', \tilde \mu')}), \cdots, \ua_n(\tw^{st}_{(\tilde \l', \tilde \mu')})\}=\{\ua_1(\tw^f_{(\tilde \l', \tilde \mu')}), \cdots, \ua_n(\tw^f_{(\tilde \l', \tilde \mu')})\}.$$ 

By definition $(\tilde \l', \tilde \mu')$ is distinguished. By Theorem \ref{61}, $\tw \tilde \to \tw'$ for some $\tw' \in S_n \cdot \tw^{st}_{(\tilde \l, \tilde \mu)}$ with $\ua_1(\tw') \ge \ua_2(\tw') \cdots \ge \ua_n(\tw) \ge (0, 0, \cdots)$. Since $\ua_1(\tw^f_{(\tilde \l', \tilde \mu')}) \ge \cdots \ge \ua_n(\tw^f_{(\tilde \l', \tilde \mu')})$, $\ua_i(\tw')=\ua_i(\tw^f_{(\tilde \l', \tilde \mu')})$ for all $i$. Hence we may assume that $\tw'=t^\chi w'$ and $\tw^f_{(\tilde \l', \tilde \mu')}=t^\chi w$ for some $\chi=[a_1, \cdots, a_n] \in \ZZ^n$ and $w, w' \in S'_n$. By definition, $\ua_i(\tw')=(a_i, \ua_{(w') \i(i)}(\tw'))$ and $\ua_i(\tw^f_{(\tilde \l', \tilde \mu')})=(a_i, \ua_{w \i (i)}(\tw^f_{(\tilde \l', \tilde \mu')}))$. Thus $$\ua_{w \i(i)}(\tw^f_{(\tilde \l', \tilde \mu')})=\ua_{w \i(i)}(\tw')=\ua_{(w') \i(i)}(\tw^f_{(\tilde \l', \tilde \mu')})$$ for any $i$. In other words, for any $i \in \{\pm 1, \cdots, \pm n\}$ and $m \in \NN$, $\ua_i(\tw^f_{(\tilde \l', \tilde \mu')})=\ua_{(w' w \i)^m(i)}(\tw^f_{(\tilde \l', \tilde \mu')})$. 

Notice that $w \i w'$ is generated by $(i, (w' w \i)(i))$. 

If $\tW^!=\tW^!_A$ and $\tilde \mu=\tilde \O$, then $w, w' \in S_n$ and $\tilde \mu'=\tilde \O$. By Lemma \ref{qqq}, $e_i-e_{(w \i w')(i)} \in \Phi_{I(\tw^f_{(\tilde \l', \tilde \mu')})}$ for all $i$ and $w' w \i \in W_{I(\tw^f_{(\tilde \l', \tilde \mu')})}$. So $\tw'=w' w \i \tw^f_{(\tilde \l', \tilde \mu')} \in W_{I(\tw^f_{(\tilde \l', \tilde \mu')})} \tw^f_{(\tilde \l', \tilde \mu')}$. 

If $\tW^!$ is of type BCD and $w \i w' \in W$, then again by Lemma \ref{qqq}, $w' w \i \in W_{I(\tw^f_{(\tilde \l', \tilde \mu')})}$. So $\tw'=w' w \i \tw^f_{(\tilde \l', \tilde \mu')} \in W_{I(\tw^f_{(\tilde \l', \tilde \mu')})} \tw^f_{(\tilde \l', \tilde \mu')}$. 

Now we consider the case that $\tW^!=\tW^!_D$ and $w' w \i \in S'_n-W$. Then there exists $i>0$ with $w' w \i(i)<0$. Since $\ua_i(\tw^f_{(\tilde \l', \tilde \mu')})=\ua_{w' w \i(i)}(\tw^f_{(\tilde \l', \tilde \mu')})$ and $\ua_i(\tw^f_{(\tilde \l', \tilde \mu')}), \ua_{-w' w \i(i)}(\tw^f_{(\tilde \l', \tilde \mu')}) \ge (0, 0, \cdots)$, we have that $\ua_i(\tw^f_{(\tilde \l', \tilde \mu')})=\ua_{w' w \i(i)}(\tw^f_{(\tilde \l', \tilde \mu')})=(0, 0, \cdots)$. In particular, $\ua_n(\tw^f_{(\tilde \l, \tilde \mu')})=(0, 0, \cdots)$ and $I(\tw^f_{(\tilde \l', \tilde \mu')})=\iota(I(\tw^f_{(\tilde \l', \tilde \mu')})$ with $\iota=(n, -n)$ be an outer automorphism of $W$. Notice that $\ua_i(\tw^f_{(\tilde \l', \tilde \mu')})=\ua_{(w' (n, -n) w \i)^m(i)}(\tw^f_{(\tilde \l', \tilde \mu')})$. By Lemma \ref{qqq}, $w' (n, -n) w \i  \in W_{I(\tw^f_{(\tilde \l', \tilde \mu')})}$. So $\tw'=w' (n, -n) w \i \tw^f_{(\tilde \l', \tilde \mu')} (n, -n) \in W_{I(\tw^f_{(\tilde \l', \tilde \mu')})} \tw^f_{(\tilde \l', \tilde \mu')} (n, -n)$.




\section{Good elements} 

\subsection{}\label{ne} Notice that each conjugacy class of $\tW$ lies in a coset of $W_a$. Let $\eta: \tW \to \tW/W_a$ be the natural projection. Then $\eta$ is constant on each conjugacy class of $\tW$. 

Let $P^\vee_{\QQ}=P^\vee \otimes_\ZZ \QQ$. Then the action of $W$ on $P^\vee$ extends in a natural way to an action on $P^\vee_\QQ$ and the quotient $P^\vee_\QQ/W$ can be identified with $$(P^\vee_\QQ)_+=\{\chi \in P^\vee_\QQ; <\chi, \a> \ge 0, \forall \a \in R^+\}.$$ 

For each element $\tw \in \tW$, there exists $n \in \NN$ such that $\tw^n=t^\chi$ for some $\chi \in P^\vee$. Let $v_{\tw}=\chi/n \in P^\vee_\QQ$ and $[v_{\tw}]$ the corresponding element in $P^\vee_\QQ/W$. It is easy to see that $v_{\tw}$ is independent of the choice of $n$. Moreover, if $\tw, \tw'$ are in the same conjugacy class of $\tW$, then $\tw^n$ is conjugated to $(\tw')^n$. Hence $(\tw')^n=t^{\chi'}$ for some $\chi' \in W \cdot \chi$. Therefore $v_{\tw}$ and $v_{\tw'}$ are in the same $W$-orbit. We call the map $\tw \mapsto [v_{\tw}]$ the {\it Newton map}. It is constant on each conjugacy class of $\tW$. 

Define $f: \tW \to (P^\vee_\QQ)_+ \times \tW/W_a$ by $\tw \mapsto ([v_{\tw}], \eta(\tw))$. Then the map $f$ is constant on each conjugacy class of $\tW$. This map is the restriction to $\tW$ of the map $G(L) \to (P^\vee_\QQ)_+ \times \tW/W_a$ in \cite[4.13]{Ko97}. We denote the image of the map $f$ by $B(\tW)$. 

\begin{lem}
Let $\tw \in \tW$. Then the following conditions are equivalent:

(1) For any $n \in \NN$, $l(\tw^n)=n l(\tw)$. 

(2) Let $v$ be the unique element in $W \cdot v_{\tw} \cap (P^\vee_\QQ)_+$. Then $$l(\tw)=<v, 2 \rho>,$$ where $\rho$ is the half sum of positive roots. 
\end{lem}

Assume that $\tw^m=t^\chi$ for some $\chi \in P^\vee$. 

If $l(\tw^n)=n l(\tw)$ for any $n \in \NN$, then $l(\tw)=\frac{1}{m} l(t^\chi)$. Let $\chi' \in W \cdot \chi \cap P^\vee_+$, then $l(t^\chi)=l(t^{\chi'})=<\chi', 2 \rho>$ and $v=\chi'/m$. Hence $l(\tw)=\frac{1}{m}<\chi', 2 \rho>=<v, 2 \rho>$. 

On the other hand, if $l(\tw)=<v, 2 \rho>$, then $$l(\tw^m)=l(t^\chi)=<\chi', 2 \rho>=m l(\tw).$$ Let $n \in \NN$, then there exists $k \in \NN$ such that $n \le k m$. We have that $l(\tw^{m k})=l(t^{k \chi})=<k \chi', 2 \rho>=m k l(\tw)$ and $$m k l(\tw)=l(\tw^{m k}) \le l(\tw^n)+l(\tw^{m k-n}) \le n l(\tw)+(m k-n) l(\tw)=m k l(\tw).$$ Therefore, both inequalities above are actually equalities. In particular, $l(\tw^n)=n l(\tw)$. 

\subsection{} We call an element $\tw \in \tW$ a {\it good element} if it satisfies the conditions in the previous lemma. The following result characterizes the good elements.

\begin{prop}\label{good}
Let $C$ be a fiber of $f: \tW \to B(\tW)$ and $\tw \in C$. Then $\tw$ is a good element if and only if $\tw \in C_{\min}$. 
\end{prop}

Notice that for any $\tx \in \tW$ and $n \in \NN$, $l(\tx) \ge \frac{1}{n} l(\tx^n)$. In particular, $l(\tx) \ge <v, 2 \rho>$, where $v$ is the unique element in $W \cdot v_{\tx} \cap (P^\vee_\QQ)_+$. Hence if $\tw$ is a good element, then $\tw$ is a minimal length element in $C$.

Let $O$ be the $\d_F$-conjugacy class on $G(L)$ whose image under the map $G(L) \to (P^\vee_\QQ)_+ \times \tW/W_a$ equals $f(C)$. By \cite[Proposition 13.1.3 \& Corollary 13.2.4]{GHKR}, there exists a good element $\tx$ such that $I \dot \tx I \subset O$. Since $\dot x \in O$, $\tx \in C$. Therefore $\tx$ is a minimal length element and all the minimal length elements in $C$ are good elements.  

\begin{cor}
Let $\tw$ be a Coxeter element of $W_a$. Then $w$ is a good element. 
\end{cor}

\begin{rmk}
This result was first proved by Speyer in \cite{Spe}. Here we give a different proof. 
\end{rmk}

By \cite{Ho82}, $\tw$ has infinite order. Hence $[v_{\tw}] \neq 0$. Let $C$ be the fiber of the map $f: \tW \to B(\tW)$ that contains $\tw$. If $w' \in C$ and $l(\tw')<l(\tw)$, then $\tw'$ lies in some $W_J$ with $J \neq \tS$. In particular, $\tw'$ lies in some finite Weyl group and $[v_{\tw'}]=0$. That is a contradiction. So $\tw$ is a minimal length element in $C$. By Proposition \ref{good}, $\tw$ is a good element. 

\

In the rest of this section, we discuss good elements for $\tW^!_A$ and $\tW^!_C$. 

\begin{lem}
(1) Let $\tW^!=\tW^!_A$. Then each fiber of the map $f: \tW^! \to B(\tW^!)$ is of the form $\sqcup_{d(\tilde \l', \tilde \O)=(\tilde \l, \tilde \O)} [(\tilde \l', \tilde \O)]$ for some distinguished $(\tilde \l, \O) \in \cd \cp$. 

(2) Let $\tW^!=\tW^!_C$. Then each fiber of the map $f: W_a \to B(\tW^!)$ is of the form $\sqcup_{d'(\tilde \l', \tilde \mu)=(\tilde \l, \tilde \O)} [(\tilde \l', \tilde \mu)]'$ for some distinguished $(\tilde \l, \O) \in \cd \cp_{\ge 0}$. 
\end{lem}

(1) For any $\tilde \l'=[(b_1, c_1), \cdots, (b_k, c_k)]$ and $\tw' \in [(\tilde \l', \tilde \O)]$, $[v_{\tw'}]$ is the $S_n$-orbit of $(e_1, \cdots, e_n)$, where $e_i=\frac{c_j}{b_j}$ for $b_1+\cdots+b_{j-1}<i \le b_1+\cdots+b_j$, i.e, the multiset whose elements are $(\frac{c_1}{b_1}, \cdots, \frac{c_1}{b_1})$ (with $b_1$-times), $\cdots$, $(\frac{c_k}{b_k}, \cdots, \frac{c_k}{b_k})$ (with $b_k$-times). Let $(\tilde \l, \tilde \O)=d(\tilde \l', \tilde \O)$ and $\tw \in [(\tilde \l, \tilde \O)]$, then $[v_{\tw}]$ is also the $S_n$-orbit of $(e_1, \cdots, e_n)$ and $\eta(\tw)=\eta(\tw')$. Thus $f(\tw)=f(\tw')$. So for any distinguished $(\tilde \l, \tilde \O) \in \cd \cp$, $\sqcup_{d(\tilde \l', \tilde \O)=(\tilde \l, \tilde \O)} [(\tilde \l', \tilde \O)]$ lies in a single fiber of $f: \tW^! \to B(\tW^!)$. On the other hand, it is easy to see that different distinguished double partitions correspond to different multisets defined above. Thus $f([(\tilde \l, \tilde \O)]) \neq f([\tilde \l_1, \tilde \O)])$ for distinguished double partitions $(\tilde \l, \tilde \O)$ and $(\tilde \l_1, \tilde \O)$. Part (1) is proved. 

(2) For any $\tilde \l'=[(b_1, c_1), \cdots, (b_k, c_k)]$ and $\tilde \mu=[(b_{k+1}, c_{k+1}), \cdots, (b_l, c_l)]$ with $(\tilde \l' \tilde \mu) \in \cd \cp_{\ge 0}$ and $\tw \in [(\tilde \l', \tilde \mu)]', \tw' \in [d(\tilde \l', \tilde \mu)]', \tw'' \in [d'(\tilde \l', \tilde \mu)]'$, we have that $[v_{\tw}]=[v_{\tw'}]=[v_{\tw''}]$ is the $S'_n$-orbit of $(e_1, \cdots, e_n)$, where $e_i=\frac{c_j}{b_j}$ for $b_1+\cdots+b_{j-1}<i \le b_1+\cdots+b_j$ and $e_i=0$ for $i>b_1+\cdots+b_k$. Hence for any distinguished $(\tilde \l, \tilde \O) \in \cd \cp_{\ge 0}$, $\sqcup_{d'(\tilde \l', \tilde \mu)=(\tilde \l, \tilde \O)} [(\tilde \l', \tilde \mu)]'$ lies in a single fiber of $f: W_a \to B(\tW^!)$. On the other hand, $f([(\tilde \l, \tilde \O)]') \neq f([\tilde \l_1, \tilde \O)]')$ for distinguished double partitions $(\tilde \l, \tilde \O)$ and $(\tilde \l_1, \tilde \O)$ in $\cd \cp_{\ge 0}$. Part (2) is proved. 

\begin{prop}\label{Bruhat2}
(1) Let $(\tilde \l, \tilde \O) \neq (\tilde \l', \tilde \O) \in \cd \cp$ with $d(\tilde \l', \tilde \O)=(\tilde \l, \tilde \O)$. Then for any $\tw' \in [(\tilde \l', \tilde \O)]$, there exists $\tw \in [(\tilde \l, \tilde \O)]_{\min}$ with $\tw<\tw'$ for the Bruhat order of $\tW^!_A$. 

(2) Let $(\tilde \l, \tilde \O) \neq (\tilde \l', \tilde \mu) \in \cd \cp_{\ge 0}$ with $d'(\tilde \l', \tilde \mu)=(\tilde \l, \tilde \O)$. Then for any $\tw' \in [(\tilde \l', \tilde \mu)]'$, there exists $\tw \in [(\tilde \l, \tilde \O)]'_{\min}$ with $\tw<\tw'$ for the Bruhat order of $\tW^!_C$. 
\end{prop}

(1) By Theorem \ref{62}, there exists $x \in W$ such that $\tw' \tilde \to x \tw^f_{(\tilde \l, \tilde \O)}$ and $l(x \tw^f_{(\tilde \l, \tilde \O)})=l(x)+l(\tw^f_{(\tilde \l, \tilde \O)})$. Since $(\tilde \l, \tilde \O) \neq (\tilde \l', \tilde \O)$, $x>1$ and $x \tw^f_{(\tilde \l, \tilde \O)}>\tw^f_{(\tilde \l, \tilde \O)}$. By Theorem \ref{thA} and \cite[Lemma 4.4]{He072}, there exists $\tw \in [(\tilde \l, \tilde \O)]_{\min}$ with $\tw<\tw'$. Part (1) is proved. 

(2) By the same argument as we did in part (1), there exists $\tw_1 \in [d(\tilde \l', \tilde \mu)]'_{\min}$ with $\tw_1<\tw'$. If $d(\tilde \l', \tilde \mu)=(\tilde \l, \tilde \O)$, then the statement is proved. Now suppose that $d(\tilde \l, \tilde \mu)=(\tilde \l'', \tilde \mu') \neq (\tilde \l, \tilde \O)$. Then $d(\tilde \l'', \underline{\tilde \mu'})=(\tilde \l, \tilde \O)$. Let $\t=t^{[\frac{1}{2}, \cdots, \frac{1}{2}]} w_{S-\{n\}} w_S \in \tW^!$. Then $l(\t)=0$ and $\t \i \tw_1 \t \in [(\tilde \l'', \underline{\tilde \mu'})]'$. Again by the same argument as we did in part (1), there exists $\tw \in [(\tilde \l, \tilde \O)]'_{\min}$ with $\tw<\t \i \tw_1 \t$. Hence $\t \tw \t \i<\tw_1<\tw'$. Since $\t [(\tilde \l, \tilde \O)]' \t \i=[(\tilde \l, \tilde \O)]'$, $\t \tw \t \i \in [(\tilde \l, \tilde \O)]'_{\min}$. Part (2) is proved.  

\begin{cor}\label{xx}
Let $\tW^!$ be of type A or C, $\aa \in f(\tW^!_{int})$ and $\tw \in f \i(\aa)$. Then the following conditions are equivalent:

(1) $\tw$ is a good element; 

(2) $\tw$ is a minimal length element in the unique distinguished conjugacy class $\co$ in $f \i(\aa)$. 

(3) $\tw$ is a minimal element for the Bruhat order on $f \i(\aa)$. 
\end{cor}

By Proposition \ref{good}, if $\tw$ is a good element, then $\tw$ is a minimal length element in $f \i(\aa)$, hence 
a minimal element for the Bruhat order on $f \i(\aa)$. This proves that $(1) \Rightarrow (3)$. Also $(3) \Rightarrow (2)$ follows from Propositions \ref{Bruhat1} and \ref{Bruhat2}. Therefore $(f \i(\aa))_{\min}=\co_{\min}$. Hence $(1) \Leftrightarrow (2)$. 

\begin{cor}\label{good3}
Let $\tW^!$ be of type A or C and $\tw, \tw' \in \tW^!_{int}$ be good elements. Then $\tw$ and $\tw'$ are in the same fiber of the map $f: \tW_{int} \to B(\tW)$ if and only if $\tw \approx \tw'$. 
\end{cor}

This follows from Corollary \ref{xx} and Theorem \ref{62}. 

\subsection{}\label{porr} Let $\tW^!$ be of type A or C. Set $B(\tW^!_{int})=f(\tW^!_{int})$. For $\aa \in B(\tW^!_{int})$ and $\tw \in \tW^!_{int}$, we write $\aa \preceq \tw$ if there exists $\tw' \in f \i(\aa)_{\min}$ with $\tw' \le \tw$ for the Bruhat order in $\tW$. By \cite[Lemma 4.4]{He072}, 

(a) if $\tw_1 \tilde \to \tw$ and $\aa \preceq \tw$, then $\aa \preceq \tw_1$. 

For $\aa, \aa' \in B(\tW^!_{int})$, we write $\aa' \preceq \aa$ if there exists $\tw \in f\i(\aa)_{\min}$ with $\aa' \preceq \tw$. By Corollary \ref{por}, 

(b) $\preceq$ is a partial order on $B(\tW^!_{int})$.

\section{Reductive method}

\subsection{}\label{Hecke} Let $H$ be the Hecke algebra associated to $\tW^!$, i.e., $H$ is the associated $\ZZ[v, v \i]$-algebra with basis $T_{\tw}$ for $\tw \in \tW^!$ and multiplication is given by 
\begin{gather*} T_{\tx} T_{\ty}=T_{\tx \ty}, \quad \text{ if } l(\tx)+l(\ty)=l(\tx \ty); \\ (T_s-v)(T_s+v \i)=0, \quad \text{ for } s \in \tS. \end{gather*}

Then $T_s \i=T_s-(v-v \i)$ and $T_{\tw}$ is invertible in $H$ for all $\tw \in \tW^!$. 

By \cite[Lemma 8.2.1]{GP00}, we have that 

(1) If $\tw \tilde \sim \tw'$, then $T_{\tw} \equiv T_{\tw'} \mod [H, H]$. 

It is also easy to check that 

(2) If $l(s \tw s)<l(\tw)$ for some $s \in \tS$, then $$T_{\tw} \equiv T_s^2 T_{s \tw s}=(v-v \i) T_{\tw s}+T_{s \tw s} \mod [H, H].$$ 

\subsection{} Let $\tW^!$ be of classical type, $\co$ be an integral conjugacy class in $\tW^!$ and $\tw \in \tW^!_{int}$. We construct some polynomials $f_{\tw, \co} \in \ZZ[v-v \i]$ as follows. 

If $\tw$ is a minimal element in the conjugacy class of $\tW^!$ that contains it, then we set $f_{\tw, \co}=\begin{cases} 1, & \text{ if } \tw \in \co \\ 0, & \text{ if } \tw \notin \co \end{cases}$. 

If $\tw$ is not a minimal element in the conjugacy class of $\tW^!$ that contains it and that $f_{\tw', \co}$ is already defined for all integral elements $\tw' \in \tW^!$ with $l(\tw')<l(\tw)$, then by Theorem \ref{thA}, there exists $\tw_1 \approx \tw$ and $s \in \tS$ such that $l(s \tw_1 s)<l(\tw_1)=l(\tw)$, we define $f_{\tw, \co}$ as $$f_{\tw, \co}=(v-v \i) f_{\tw_1 s, \co}+f_{s \tw_1 s, \co}.$$

This completes the definition of $f_{\tw, \co}$. One also sees from the definition that all the coefficients of $f_{w, \co}$ are nonnegative integer. 

By \ref{Hecke} (1) and Theorem \ref{thA}, 

(1) for any $w, w' \in \co_{\min}$, $T_w \equiv T_{w'} \mod [H, H]$. 

Now we choose a representative $\tw_\co \in \co_{\min}$ for each integral conjugacy class $\co$ in $\tW^!$. By \ref{Hecke} (1) \& (2) and Theorem \ref{thA}, for any integral element $\tw \in \tW^!$, 

\[\tag{2} T_{\tw} \equiv \sum_\co f_{\tw, \co} T_{w_\co} \mod [H, H].\] 

We call $f_{\tw, \co}$ the {\it class polynomials}. 

Notice that the definition of $f_{\tw, \co}$ depends on the choice of the sequence of elements in $\tS$ used to conjugate $\tw$ to a minimal length element in its conjugacy class. We expect that $f_{\tw, \co}$ is in fact, independent of such choice and is uniquely determined by the condition (2) above. This is true if one replaces $\tW^!$ by a finite Coxeter group and $H$ by the corresponding finite Hecke algebra (see \cite[Theorem 4.2]{GR}). 

\

In the rest of this paper, we discuss some applications on the loop groups. 

\subsection{} Similar to the definition of $\tW^!$, we define $G^!$ as follows.

Set $G^!_A=GL_n(L)$, $G^!_B=PSO_{2n+1}(L)$, $G^!_C=PSP_{2n}(L)$ and $G^!_D=PSO_{2n}(L) \rtimes <\iota>$, where $\iota$ is the outer diagonal automorphism on $G$ whose induces action on $W(D_n)$ is the conjugation by $(n, -n)$. We say that $G^!$ is {\it of classical type} if $G^!=G^!_*$ for $*$ is A, B, C or D. Then $\tW^!$ is the extended affine Weyl group of $G^!$. Moreover, the bijective group homomorphisms $\d_a$ and $\d_F$ on $G(L)$ defined in $\S$\ref{sa} extend in a natural way to group homomorphisms on $G^!$, which we still denote by the same symbol. Unless otherwise stated, we write $\d$ for $\d_a$ or $\d_F$. 

Recall that $\Om=\{\t \in \tW^!; l(\t)=0\}$. Let $G'=\sqcup_{\tw \in W_a} I \dot \tw I$ be the identity component of $G^!$. Then \begin{align*} G &=\sqcup_{\tw \in \tW^!} I \dot \tw I=\sqcup_{\t \in \Om} \sqcup_{\tw \in W_a} I \dot \tw \dot \t I \\&=\sqcup_{\t \in \Om} \sqcup_{\tw \in W_a} I \dot \tw I \dot \t=\sqcup_{\t \in \Om} G' \dot \t. \end{align*}

\begin{lem}\label{red11}
Let $\tw, \tw' \in \tW^!$. 

(1) If $\tw \to \tw'$, then $$G' \cdot_\d I \dot \tw I \subset G' \cdot_\d I \dot \tw' I \cup \cup_{\tx \in \tw W_a, l(\tx)<l(\tw)} G' \cdot_\d I \dot \tx I.$$ 

(2) If $\tw \approx \tw'$, then $$G' \cdot_\d I \dot \tw I=G' \cdot_\d I \dot \tw' I.$$ 
\end{lem}

By definition, there exists a finite sequence $\tw=\tw_0 \xrightarrow {i_1} \tw_1 \xrightarrow {i_2} \cdots \xrightarrow {i_m} \tw_n=\tw'$, where $i_j \in \tS$ for all $j$. We prove the lemma by induction on $m$. 

The statements are true for $m=0$. Now assume that the statements hold for $m-1$. By \cite[Lemma 1.6.4]{DL}, we have that $\tw=\tw_1$ or $s_{i_1} \tw<\tw$ or $\tw s_{i_1}<\tw$. If $\tw=\tw_1$, then the statements follow from induction hypothesis. Now we prove the case where $s_{i_1} \tw<\tw$. The case $\tw s_{i_1}<\tw$ can be proved in the same way. 

Since $s_{i_1} \tw<\tw$, then $G' \cdot_\d I \dot \tw I=G' \cdot_\d I \dot s_{i_1} I \dot s_{i_1} \dot \tw I=G' \cdot_\d I \dot s_{i_1} \dot \tw I \dot s_{i_1} I$. Moreover, \[I \dot s_{i_1} \dot w I \dot s_{i_1} I=\begin{cases} I \dot \tw_1 I, & \text{ if }  l(\tw_1)=l(s_{i_1} \tw)+1=l(\tw); \\ I \dot s_{i_1} \dot \tw I \sqcup I \dot \tw_1 I, & \text{ if } l(\tw_1)=l(s_i \tw)-1=l(\tw)-2. \end{cases}\] In either case, $$G' \cdot_\d I \dot \tw I \subset G' \cdot_\d I \dot \tw_1 I \cup \cup_{\tx \in \tw W_a, l(\tx)<l(\tw)} G' \cdot_\d I \dot \tx I.$$ Notice that $l(\tw_1) \le l(\tw)$. By induction hypothesis, $G' \cdot_\d I \dot \tw_1 I \subset G' \cdot_\d I \dot \tw' I \cup \cup_{\tx \in w W_a, l(\tx)<l(\tw_1)} G' \cdot_\d I \dot \tx I$. Hence $G' \cdot_\d I \dot \tw I \subset G' \cdot_\d I \dot \tw' I \cup \cup_{\tx \in w W_a, l(\tx)<l(\tw)} G' \cdot_\d I \dot \tx I$. 

If moreover, $\tw \approx \tw'$, then $l(\tw_1)=l(\tw)$ and $\tw_1 \approx \tw'$. By induction hypothesis, $G' \cdot_\d I \dot \tw I=G' \cdot_\d I \dot \tw_1 I=G' \cdot_\d I \dot \tw' I$. 

\begin{lem}\label{red2}
Let $\tw, \tw' \in \tW^!$. 

(1) If $\tw \tilde \to \tw'$, then $$G^! \cdot_\d I \dot \tw I \subset G^! \cdot_\d I \dot \tw' I \cup \cup_{\tx \in \tw W_a, l(\tx)<l(\tw)} G^! \cdot_\d I \dot \tx I.$$ 

(2) If $\tw \tilde \approx \tw'$, then $$G^! \cdot_\d I \dot \tw I=G^! \cdot_\d I \dot \tw' I.$$ 
\end{lem}

For any $\tx \in \tW$ and $\t \in \Om$, we have that $$G^! \cdot_\d I \dot \tx I=G^! \cdot_\d \dot \t I \dot \tx I \s(\dot \t) \i=G^! \cdot_\d I \dot \t \dot \tx \dot \t \i I.$$ Now if $\tw \tilde \to \tw'$, then there exists $\t \in \tW$ with $l(\t)=0$ such that $\tw \to \t \tw' \t \i$. By Lemma \ref{red11},  \begin{align*} G^! \cdot_\d I \dot \tw I &=G^! \cdot_\d (G' \cdot_\d I \dot \tw I) \\ &\subset G^! \cdot_\d \bigl(G' \cdot_\d I \dot \t \dot \tw' \dot \t \i I \cup \cup_{\tx \in w W_a, l(\tx)<l(\tw)} G' \cdot_\d I \dot \tx I \bigr) \\ &=G^! \cdot_\d I \dot \t \dot \tw' \t \dot \i I \cup \cup_{\tx \in \tw W_a, l(\tx)<l(\tw)} G^! \cdot_\d I \dot \tx I \\ &=G^! \cdot_\d I \dot \tw' I \cup \cup_{\tx \in \tw W_a, l(\tx)<l(\tw)} G^! \cdot_\d I \dot x I. \end{align*} 

If $\tw \tilde \approx \tw'$, then there exists $\t \in \tW$ with $l(\t)=0$ such that $\tw \approx \t \tw' \t \i$. By Lemma \ref{red11},  \begin{align*} G^! \cdot_\d I \dot \tw I &=G^! \cdot_\d (G' \cdot_\d I \dot \tw I) =G^! \cdot_\d (G' \cdot_\d I \dot \t \dot \tw' \dot \t \i I) \\ &=G^! \cdot_\d I \dot \t \dot \tw' \dot \t \i I=G^! \cdot_\d I \dot \tw' I. \end{align*} The Lemma is proved. 

\

The following result is proved in \cite[Lemma 2.4 \& 2.8]{H7}. 

\begin{lem}\label{red3} Let $\tw \in {}^S \tW^!$ and $\tw'=x \tw$ for some $x \in W_{I(w)}$, then $G' \cdot_\d I \dot \tw' I \subset G' \cdot_\d I \dot \tw I$. 
\end{lem}

\subsection{} For any $b \in G^!$ and $\tw \in \tW^!$, define $$X_{\tw, \d}(b)=\{g I \in G^!/I; g \i b \d(g) \in I \dot \tw I\}.$$ Now we discuss some reductive method for $X_{\tw, \d}(b)$. The first two Lemmas below are proved in \cite{GH}. 

\begin{lem}\label{indd}
Let $\tw \in \tW^!$, and let $s\in\tS$. Then 

(1) If $l(s\tw s)=l(\tw)$, then $X_{\tw, \d}(b)$, $X_{s \tw s, \d}(b)$ are universally homeomorphic. 

(2) If $l(s \tw s) =l(\tw)-2$, then $X_{\tw, \d}(b)$ can be written as a disjoint union $X_{\tw, \d}(b) = X_1 \cup X_2$ where $X_1$ is closed and $X_2$ is open, and such that $X_1$ admits a morphism to $X_{s \tw s, \d}(b)$, all of whose fibers are isomorphic to $\mathbb A^1$, and such that $X_2$ admits a morphism to $X_{s \tw, \d}(b)$, all of whose fibers are isomorphic to $\mathbb A^1\setminus\{ 0\}$. 
\end{lem}

\begin{lem}\label{tt}
Let $\tx, \t \in \tW$ with $l(\t)=0$. Then $X_{\tx, \d}(b)$ is isomorphic to $X_{\t \tx \t \i, \d}(b)$. 
\end{lem}

As a consequence of Lemma \ref{indd} and \ref{tt}, we have that 

\begin{cor}\label{tttt}
Let $\tw, \tw' \in \tW^!$ with $\tw \tilde \approx \tw'$. Then $X_{\tw, \d}(b)$ and $X_{\tw', \d}(b)$ are universally homeomorphic. 
\end{cor}

\begin{lem}\label{lang}
Let $\tw \in {}^S \tW^!$ and $x \in W_{I(w)}$. Then $\dim(X_{x \tw, \d_F}(b))=\dim(X_{\tw, \d_F}(b))+l(x)$ and $\dim(X_{x \tw, \d_a}(b)) \le \dim(X_{\tw, \d_a}(b))+l(w_{I(w)})$.
\end{lem}

\begin{rmk}
By convention, we set $\dim(\O)=-\infty$. 
\end{rmk}

Let $J=I(\tw)$, $L_J$ be the standard Levi subgroup of $G$ corresponding to $J$ and $P=L_J I$ be the standard parahoric subgroup of $G^!$. Set $$X=\{g P \in G^!/P; g \i b \d(g) \in P \tw P\}.$$ Notice that for any $u \in W_J$, $I \dot u \dot \tw I=I \dot u I \dot \tw I \subset P \dot \tw P$. Thus the map $g I \mapsto g P$ sends $X_{u \tw, \d}(b)$ to $X$. 

Let $U_J$ be the subgroup of $G$ generated by $u_\a(\kk)$ for $\a \in \Phi^-\setminus\Phi^-_J$ and $I_J$ be the inverse image of $U_J$ under the map $G(o) \to G$. Then $I_J$ is normal in $P_J$ and $P_J=I_J L_J$. Thus $P \dot \tw P=I_J L_J \dot \tw L_J I_J=I_J L_J \dot \tw I_J$. 

Let $g \in G^!$ with $g P \in X$. Then $g \i b \d(g) \in I_J l \dot \tw I_J$ for some $l \in L_J$. Now for $p \in L_J$, $p \i g \i b \d(g) \d(p) \in p \i I_J l \dot \tw I_J \d(p) \i=I_J p \i l \dot \tw \d(p) \i I_J$. Notice that $I=I_J (B^- \cap L_J)=(B^- \cap L_J) I_J$ and $I \dot u \dot \tw I=I_J (B^- \cap L_J) \dot u (B^- \cap L_J) \dot \tw I_J$. Thus for $p \in L_J$, $g p I \in X_{u \tw, \d}(b)$ if and only if $p \i l \dot \tw \d(p) \dot \tw \i \in (B^- \cap L_J) \dot u (B^- \cap L_J)$. 

We also have that $P/I \cong L_J/(B^- \cap L_J)$. Define $\d': L_J \to L_J$ by $\d'(l)=\dot \tw \d(l) \dot \tw \i$. Set $$Y_g=\{p (B^- \cap L_J) \in L_J/(B^- \cap L_J); p \i l \d'(p) \in (B^- \cap L_J) \dot u (B^- \cap L_J)\}.$$ Then $\{p I \in P/I; g p I \in X_{u \tw, \d} (b)\} \cong Y_g$ is a subvariety of the flag variety of $L_J$. In particular, each fiber of the map $X_{u \tw, \d}(b) \to X$ is at most of dimension $l(w_J)$. 

If $\d=\d_F$, $Y_g$ is a Deligne-Lusztig variety in $L_J/(B^- \cap L_J)$ and it is known that it is of dimension $l(u)$. So $\dim(X_{u \tw, \d_F}(b))=\dim(X)+l(u)$ for any $u \in W_J$. In particular, $\dim(X_{\tw, \d_F}(b))=\dim(X)$ and $\dim(X_{x \tw, \d_F}(b))=\dim(X_{\tw, \d_F}(b))+l(x)$. 

If $\d=\d_a$, then by \cite[Lemma 7.3]{St68}, $Y_g$ is nonempty for $u=1$. Hence $\dim(X_{w, \d_a}(b)) \ge \dim(X)$. Therefore $\dim(X_{x \tw, \d_a}(b)) \le \dim(X)+l(w_J) \le \dim(X_{\tw, \d_a}(b))+l(w_J)$. 

\begin{cor}
Let $G^!$ be of classical type and $b \in G^!$. Let $\co$ be an integral conjugacy class of $\tW^!$ and $\tw, \tw' \in \co_{\min}$. Then $\dim(X_{\tw, \d_F}(b))=\dim(X_{\tw', \d_F}(b))$. 
\end{cor}

By Theorem \ref{61}, there exists $\tx \in {}^S \tW^!$ and $v, v' \in W_{I(\tx)}$ such that $\tw \tilde \approx v \tx$ and $\tw' \tilde \approx v' \tx$. Since $l(\tw)=l(\tw')$, $l(v)=l(v')$. By Corollary \ref{tttt} and Lemma \ref{lang}, \begin{gather*} \dim(X_{\tw, \d_F}(b))=\dim(X_{v \tx, \d_F}(b))=\dim(X_{\tx, \d_F}(b))+l(v), \\ \dim(X_{w', \d_F}(b))=\dim(X_{v' \tx, \d_F}(b))=\dim(X_{\tx, \d_F}(b))+l(v'). \end{gather*} Therefore $\dim(X_{\tw, \d_F}(b))=\dim(X_{\tw', \d_F}(b))$. 

\begin{prop}\label{class}
Let $G^!$ be of classical type, $b \in G^!$ and $\tw \in \tW^!_{int}$. Then 

$$\dim(X_{\tw, \d_F}(b))=\max_{\co} \frac{1}{2}(l(\tw)-l(\tw_\co)+\deg(f_{\tw, \co}))+\dim(X_{\tw_\co, \d_F}(b)),$$ where $\co$ runs over integral conjugacy classes of $\tW^!$ and $\tw_\co$ is a minimal length element in $\co$.
\end{prop}

Let $\co'$ be the integral conjugacy class that contains $\tw$. 

If $\tw \in \co'_{\min}$, then $\frac{1}{2} \deg(f_{\tw, \co})+\dim(X_{\tw_\co, \d_F}(b)) \neq -\infty$ if and only $f_{\tw, \co} \neq 0$ and $X_{\tw_\co, \d_F}(b) \neq \emptyset$, i.e., $\co=\co'$ and $b \in G^! \cdot I \dot \tw_\co I$. In this case, $f_{\tw, \co}=1$ and $\deg(f_{\tw, \co})=0$. By the previous Corollary, $\dim(X_{\tw, \d_F}(b))=\dim(X_{
\tw_\co, \d_F}(b))$. The Proposition holds in this case. 

If $\tw \notin \co'_{\min}$, we use the same sequence of elements in $\tS$ to conjugate $\tw$ to a minimal element in $\co'$ as we did in the definition of $f_{\tw, \co}$. Then there exists $\tw_1 \approx \tw$ and $s \in \tS$ such that $l(s \tw_1 s)<l(\tw_1)=l(\tw)$ and $f_{\tw, \co}=(v-v \i) f_{\tw_1 s, \co}+f_{s \tw_1 s, \co}$. Hence $\deg(f_{\tw, \co})=\max\{\deg(f_{\tw_1 s, \co})+1, \deg(f_{s \tw_1 s, \co})\}$ and $\frac{1}{2}(l(\tw)+\deg(f_{\tw, \co}))=\max\{\frac{1}{2}(l(\tw_1 s)+\deg(f_{\tw_1 s, \co}))+1, \frac{1}{2}(l(s \tw_1 s)+\deg(f_{s \tw_1 s, \co}))+1\}$.

On the other hand, by Lemma \ref{indd}, $\dim(X_{\tw, \d_F}(b))=\dim(X_{\tw_1, \d_F}(b))=\max\{\dim(X_{\tw_1 s, \d_F}(b))+1, \dim(X_{s \tw_1 s, \d_F}(b))+1\}$. Now the Proposition follows from induction on $l(\tw)$. 

\

By the same argument, we have the following result for $\d_a$.

\begin{prop}\label{class2}
Let $G^!$ be of classical type, $b \in G(L)$ and $w \in \tW^!_{int}$. If $X_{\tx, \d_a}(b)$ is finite dimensional for any $\tx \in \tW^!$ that is of minimal length in its conjugacy class of $\tW^!$, then $X_{\tw, \d_a}(b)$ is also finite dimensional. 
\end{prop}

\section{More on good elements}

We first consider $G(L) \cdot_\d I \dot \tw I$ for some good element $\tw$. 

\begin{lem}\label{dg}
Let $\tw, \tw' \in \tW$ be good elements. If $f(\tw) \neq f(\tw')$, then $$G(L) \cdot_\d I \dot \tw I \cap G(L) \cdot_\d I \dot \tw' I=\emptyset.$$ 
\end{lem}

It is easy to see that for $\tx,\ty \in \tW$, $I \dot \tx I \dot \ty I \subset \cup_{\tilde z \in x \ty W_a} I \dot {\tilde z} I$. Hence \begin{gather*} G(L) \cdot_\d I \dot \tw I \subset \cup_{\tx \in \tW} I \dot \tx I \dot \tw I \dot \tx \i I \subset \cup_{\tilde z \in \tx \tw \tx \i W_a=\tw W_a} I \dot {\tilde z} I, \\ G(L) \cdot_\d I \dot \tw' I \subset \cup_{\tx \in \tW} I \dot \tx I \dot \tw' I \dot \tx \i I \subset \cup_{\tilde z \in \tx \tw' \tx \i W_a=\tw W_a} I \dot {\tilde z} I. \end{gather*} Therefore if $\tw W_a \neq \tw' W_a$, then $G(L) \cdot_\d I \dot \tw I \cap G(L) \cdot_\d I \dot \tw' I=\emptyset$. 

Now assume that $\eta(\tw)=\eta(\tw')$. By our assumption, $[v_{\tw}] \neq [v_{\tw'}]$. If $G(L) \cdot_\d I \dot \tw I \cap G(L) \cdot_\d I \dot \tw' I \neq \emptyset$, then there exist $g \in G(L)$ such that $I \dot \tw I \cap g I \dot \tw' I \s(g) \i \neq \emptyset$. Let $z \in I \dot \tw I \cap g I \dot \tw' I \s(g) \i$. Since $\tw$ and $\tw'$ are good elements, then for any $n \in \NN$, \begin{align*} z \s(z) \cdots \s^{n-1}(z) & \in (I \dot \tw I) (I \dot \tw I) \cdots (I \dot \tw I)=I \dot \tw^n I, \\ z \s(z) \cdots \s^{n-1}(z) & \in (g I \dot \tw' I \s(g) \i) (\s(g) I \dot \tw' I \s^2(g) \i) \cdots (\s^{n-1}(g) I \dot \tw' I \s^n(g) \i) \\ &=g I (\dot \tw')^n I \s^n(g). \end{align*} In particular, for any $n \in \NN$, $I \dot \tw^n I \cap g I (\dot \tw')^n I \s^n(g) \i \neq \emptyset$. 

There exists $m \in \NN$ such that $\tw^m=t^\mu$ and $(\tw')^m=t^{\mu'}$ for some $\mu, \mu' \in P^\vee$. Since $[v_{\tw}] \neq [v_{\tw'}]$, $\mu \notin W \cdot \mu'$. Assume that $g \in I \tx I$ for some $\tx \in \tW$. Then $I \e^{k \mu} I \cap I \dot \tx I \e^{k \mu'} I \dot \tx \i I \neq \emptyset$ for all $k \in \NN$. 

Notice that $I \dot \tx I \e^{k \mu'} I \dot \tx \i I \subset \cup_{\ty, \ty' \le \tx} I \ty \e^{k \mu'} (\dot \ty') \i I$. Thus $t^{k \mu}=\ty t^{k \mu'} (\ty') \i$ for some $\ty, \ty' \le \tx$. Assume that $\ty=y t^\chi$ and $\ty'=y' t^{\chi'}$ with $\chi, \chi' \in P^\vee$ and $y, y' \in W$. Then $\ty t^{k \mu'} (\ty') \i=t^{y (k \mu'+\chi-\chi')} y (y')\i$. Hence $y=y'$ and $k \mu'+\chi-\chi'=k y \i \mu$. By definition, $l(t^\chi) \le l(\ty)+l(y) \le l(\tx)+l(w_S)$. Similarly, $l(t^{\chi'}) \le l(\tx)+l(w_S)$. Since $\mu \notin W \cdot \mu'$, then $l(t^{\mu'-y \i \mu}) \ge 1$. Now $$k \le l(t^{k(\mu'-y \i \mu)})=l(t^{\chi'-\chi}) \le l(t^{\chi'})+l(t^\chi) \le 2 l(\tx)+2 l(w_S).$$ That is a contradiction.  

\begin{cor}\label{single}
Let $\tw, \tw' \in \tW$ be good elements with $f(\tw)=f(\tw')$. Then $G(L) \cdot_{\d_F} I \dot \tw I=G(L) \cdot_{\d_F} I \dot \tw' I$ is a single $\d_F$-conjugacy class. 
\end{cor}

As in the proof of Proposition \ref{good}, for any $\d_F$-conjugacy class $O$ of $G(L)$, there exists a good element $\tx_O \in \tW$ such that $\dot \tx_O \in O$. If $O \subset G(L) \cdot_{\d_F} I \dot w I$, then we must have that $\dot \tx_O \in G(L) \cdot_{\d_F} I \dot \tw I$. By the previous Lemma, $f(\tx_O)=f(\tw)$. By \cite[4.13]{Ko97}, $O$ is uniquely determined by $f(\tx_O)$. Therefore $G(L) \cdot_{\d_F} I \dot \tw I$ is the single $\d_F$-conjugacy class that contains $\dot \tx_O$. In particular,  $G(L) \cdot_{\d_F} I \dot \tw I=G(L) \cdot_{\d_F} I \dot \tw' I$. 

\subsection{}\label{reK} Now we reformulate Kottwitz's classification of $\d_F$-conjugacy classes as follows.

Let $\tW_{good}$ be the set of good elements in $\tW$. For $\tw, \tw' \in \tW_{good}$, we write $\tw \asymp \tw'$ if $f(\tw)=f(\tw')$. Then we have that \[\tag{a} G(L)=\sqcup_{[\tw] \in \tW_{good}/\asymp} G(L) \cdot_{\d_F} I \dot \tw I.\]

However, if $\d=\d_a$, then in general, $G(L) \cdot_\d I \dot \tw I$ contains infinitely many $\d$-conjugacy classes and we don't know if $G(L) \cdot_\d I \dot \tw I=G(L) \cdot_\d I \dot \tw' I$ for $\tw \asymp \tw'$. We'll see that it is true for type A and C. 

\

In the rest of this section, we discuss the dimension of $X_{\tw, \d}(b)$ for good element $\tw$. 

\begin{prop}\label{00}
Let $\tw \in \tW$ be a good element. Let $n \in \NN$ and $\mu \in P^\vee$ with $\tw^n=t^\mu$. Let $M$ be the Levi subgroup of $G$ generated by $T$ and $u_\a(\kk)$ for $\a \in R$ with $<\mu, \a>=0$. Then $X_{\tw, \d}(\dot \tw) \subset M(L) I/I \cong M(L)/(M(L) \cap I)$. Moreover, for any $b \in G(L) \cdot_\d I \dot \tw I$, $\dim(X_{\tw, \d_F}(b))=0$ and $\dim(X_{\tw, \d_a}(b))<\infty$ if $a$ is not a root of unity. 
\end{prop}

The proof will be given in $\S$\ref{pf00}. 

\subsection{}\label{comm} Let $\tilde \Phi=\{\a+n \d; \a \in \Phi, n \in \ZZ\}$ be the set of real affine roots and $\tilde \Phi^+=\{\a+n \d; \a \in \Phi^+, n>0\}\sqcup \{\a+n \d; \a \in \Phi^-, n \ge 0\}$ be the set of positive real affine roots. The affine simple roots are $-\a_i$ for $i \in S$ and $\a_0=\th+\d$. Then any positive real affine root $\a$ can be written in a unique way as $\sum_{i \in S} -a_i \a_i+a_0 \a_0$ with $a_i \in \NN \cup \{0\}$ for $i \in \tS$. We set $ht(\a)=\sum_{i \in \tS} a_i$. For any real root $\a+n \d$, we define $x_{\a+n \d}: \kk \to G(L)$ by $x_{\a+n \d}(k)=u_\a(k \e^n)$ for $k \in \kk$. 

Let $I_1$ be the inverse image of $U^-$ under the projection $G(o) \to G$, where $U^-$ is the unipotent radical of $B^-$. Then $I_1$ is generated by $x_{\a}(\kk)$ for $\a \in \tilde \Phi^+$ and $I_1 \cap T(L)$. For $n>1$, let $I_n$ be the subgroup of $I$ generated by $x_{\a}(\kk)$ with $ht(\a) \ge n$ and $I_1 \cap T(L)$. Then it is easy to see that 

(1) for $\a, \b \in \tilde \Phi^+$ and $a, b \in \kk$, $x_\a(a) x_\b (b) \in x_\b(b) x_\a(a) I_{ht(\a)+ht(\b)}$.

(2) for any $t \in I_1 \cap T(L)$, $\a \in \tilde \Phi^+$ and $a \in \kk$, $t x_\a(a) \in x_\a(a) t I_{ht(\a)+1}$. 

Thus $I_n$ is a normal subgroup of $I$ for any $n \in \NN$.

\begin{lem}\label{0}
Let $b \in I$. Then $\dim(X_{1, \d_F}(b))=0$ and $\dim(X_{1, \d_a}(b))<\infty$ if $a$ is not a root of unity. 
\end{lem} 

For $\tx \in \tW$, define $Y_{\tx}=I \dot \tx I/I \cap X_{1, \d}(b)$.  Assume that $\{\a \in \tilde \Phi^+; \tx \i \a<0\}=\{\a_{i_1}, \cdots, \a_{i_k}\}$. We arrange the roots such that $ht(\a_{i_1}) \le ht(\a_{i_2}) \le \cdots \le ht(\a_{i_k})$. Let $g \in G(L)$ with $g I \in Y_{\tx}$, then $g=x_{\a_{i_1}}(a_1) \cdots x_{\a_{i_k}}(a_k) \dot \tx I$ for some $a_1, \cdots, a_k \in \kk$. Since $g I \in X_{1, \d}(b)$, $g I \cap b \d(g) I \neq \emptyset$.

We first consider the case where $\d=\d_F$. By \cite[Prop 6.3.1]{GHKR}, any element in $I$ is $\d_F$ conjugate to $1$. So we may take $b=1$. Now $$x_{\a_{i_1}}(a_1) \cdots x_{\a_{i_k}}(a_k)=x_{\a_{i_1}}(F(a_1)) \cdots x_{\a_{i_k}}(F(a_k)).$$ Therefore $a_1=F(a_1), \cdots, a_k=F(a_k)$. So there are only finite many elements in $Y_{\tx}$ for each $\tx \in \tW$ and $\dim(X_{1, \d_F}(b))=0$. 

We then consider the case that $\d=\d_a$ with $a$ not a root of unity.  Notice that $\dim Y_{\tx}$ equals the dimension of the variety consists of $(a_1, \cdots, a_k)$ satisfying \[\tag{*} x_{\a_{i_1}}(a_1) \cdots x_{\a_{i_k}}(a_k) \dot \tx I \cap b \d(x_{\a_{i_1}}(a_1) \cdots x_{\a_{i_k}}(a_k)) \dot \tx I \neq \emptyset.\]

We may assume that $b \in t I_1$ for some $t \in T$. Let $n \in \NN$ with $t \d(x_\a(1)) t \i \neq x_\a(1)$ for all $\a$ with $ht(\a) \ge n$. Such $n$ exists since $a$ is not a root of unity. 

Let $m \ge n$ and $j \in \NN$ with $ht(\a_{i_{j-1}})<m$ and $ht(\a_{i_j}) \ge m$. We show that 

(a) for any $u \in I_m ({}^{\dot \tx} I \cap I)$, there is a unique $(a_j, \cdots, a_k)$ such that $x_{\a_{i_j}}(a_j) \cdots x_{\a_{i_k}}(a_k) \in u t \d(x_{\a_{i_j}}(a_j) \cdots x_{\a_{i_k}}(a_k)) ({}^{\dot \tx} I \cap I)$. 

Notice that (a) is obvious if $ht(\a_{i_k})<m$. We prove this statement by descending induction on $m$. 

Assume that $ht(\a_{i_j})=\cdots=ht(\a_{i_l})=m<ht(\a_{i_{l+1}})$. Then $$x_{\a_{i_j}}(a_j) \cdots x_{\a_{i_k}}(a_k) \in x_{\a_{i_j}}(a_j) \cdots x_{\a_{i_l}}(a_l) I_{m+1}.$$ 
It is easy to see that the map $\kk^{l-j} \times I_{m+1} ({}^{\dot \tx} I \cap I_m) \to I_m$ defined by $(b_j, \cdots, b_l, z) \mapsto x_{\a_{i_j}}(b_j) \cdots x_{\a_{i_l}}(b_l) z$ is a bijection. 

Let $u'=x_{\a_{i_j}}(b_j) \cdots x_{\a_{i_l}}(b_l) \in u I_{m+1} ({}^{\dot \tx} I \cap I_m)$. Since the commutator of $I_m$ and $I_m$ is contained in $I_{m+1}$, for any $v \in I_m$, $I_{m+1} ({}^{\dot \tx} I \cap I_m) t v \subset t v I_{m+1}$. Hence 
 \begin{align*} u t \d(x_{\a_{i_j}}(a_j) \cdots x_{\a_{i_k}}(a_k)) & \in u' I_{m+1} ({}^{\dot \tx} I \cap I_m) t \d(x_{\a_{i_j}}(a_j) \cdots x_{\a_{i_l}}(a_l)) I_{m+1} \\ & \subset u' t \d(x_{\a_{i_j}}(a_j) \cdots x_{\a_{i_l}}(a_l)) I_{m+1} ({}^{\dot \tx} I \cap I).\end{align*}  

Now we have that $$x_{\a_{i_j}}(a_j) \cdots x_{\a_{i_l}}(a_l) \in u' t \d(x_{\a_{i_j}}(a_j) \cdots x_{\a_{i_l}}(a_l)) t\i I_{m+1} ({}^{\dot \tx} I \cap I).$$ 

For $p \ge j$, $t \d(x_{\a_p}(a)) t \i=x_{\a_p}(c_p a)$ for some $c_p \neq 1$. Then $$u' t \d(x_{\a_{i_j}}(a_j) \cdots x_{\a_{i_l}}(a_l)) t\i \in x_{\a_{i_j}}(b_j+c_j a_j) \cdots x_{\a_{i_l}}(b_l+c_l a_l) I_{m+1}.$$ Thus we must have that $a_j=b_j+c_j a_j, \cdots, a_l=b_l+c_l a_l$. In particular, $a_j, \cdots, a_l$ are uniquely determined. 

Now set $u_1=(x_{\a_{i_j}}(a_j) \cdots x_{\a_{i_l}}(a_l)) \i u' t \d(x_{\a_{i_j}}(a_j) \cdots x_{\a_{i_l}}(a_l)) t \i \in I_{m+1} ({}^{\dot x} I \cap I)$, we have that $$x_{\a_{i_{l+1}}}(a_{l+1}) \cdots x_{\a_{i_k}}(a_k) \in u_1 t \d(x_{\a_{i_{l+1}}}(a_{l+1}) \cdots x_{\a_{i_k}}(a_k)) ({}^{\dot \tx} I \cap I).$$ By induction hypothesis on $m+1$, $a_{l+1}, \cdots, a_k$ are also uniquely determined. 

(a) is proved. 

Let $j' \in \NN$ with $ht(\a_{i_{j'-1}})<n$ and $ht(\a_{i_j'}) \ge n$. By (a) for $m=n$, given any $(a_1, \cdots, a_{j'}) \in \kk^{j'}$, there exists at most one $(a_{j'+1}, \cdots, a_k) \in \kk^{k-j'}$ such that $(a_1, \cdots, a_k)$ satisfies (*) above. In other words, $\dim(Y_{\tx}) \le \sharp\{\a \in \tilde \Phi^+; ht(\a)<n\}$ for any $\tx \in \tW$. Therefore $\dim(X_{1, \d}(b))=\max_{\tx \in \tW} Y_{\tx}<\infty$. 

\begin{lem}\label{finite}
Let $\mu \in P^\vee$ and $M$ be the Levi subgroup of $G$ generated by $T$ and $u_\a(\kk)$ for $\a \in R$ with $<\mu, \a>=0$. Then the map $$I \times_{I \cap M(L)} (I \cap M(L)) \e^\mu \to I \e^\mu I$$ defined by $(i, i') \mapsto i i' \d(i) \i$, is bijective. 
\end{lem}

\begin{rmk}
If $\s=\d_F$, then the lemma is a special case of \cite[Theorem 2.1.2]{GHKR}. The case where $\d=\d_a$ is essentially the same as in loc.cit. Here we give a proof to convince the readers that no problem occurs for $\d=\d_a$. 
\end{rmk}

Notice that $I \times_{I \cap M(L)} (I \cap M(L)) \e^\mu \cong I_1 \times_{I_1 \cap M(L)} (I \cap M(L)) \e^\mu$ and $I \e^\mu I=I_1 (I \cap M(L) \e^\mu I_1$. It suffices to prove that for any $n$, the map $I_n \times_{I_n \cap I_{n+1} (I_n \cap M(L))} I_{n+1} (I \cap M(L)) \e^\mu I_{n+1} \to I_n (I \cap M(L)) \e^\mu I_n$ defined by $(i, i') \mapsto i i' \d(i) \i$ is bijective. 

Let $P$ be the parabolic subgroup of $G$ generated by $T$ and $u_\a(\kk)$ for $\a \in R$ with $<\mu, \a> \ge 0$ and Let $P^-$ be the opposite parabolic subgroup of $G$ generated by $T$ and $u_\a(\kk)$ for $\a \in R$ with $<\mu, \a> \le 0$. Let $U_P$ be unipotent radical of $P$ and $U_{P^-}$ be the unipotent radical of $P^-$. Set $I_n'=I_n \cap U_P(L)$ and $I_n''=I_n \cap U_{P^-}(L)$. Since $I$ normalizes $I_n$ and $M(L)$ normalizes $U_{P}(L)$ and $U_{P^-}(L)$, $I \cap M(L)$ normalizes $I'_n$ and $I''_n$ for all $n$. Moreover, we have that $I_n=I'_n I''_n (I_n \cap M(L))$. Now \begin{align*} I_n (I \cap M(L)) \e^\mu I_n & \subset I_n \cdot_\d I_n (I \cap M(L)) \e^\mu=I_n \cdot_\d I'_n I''_n (I \cap M(L)) \e^\mu \\ &=I_n \cdot_\d I''_n (I \cap M(L)) \e^\mu I'_n=I_n \cdot_\d (I \cap M(L)) I''_n \e^\mu I'_n. \end{align*} By definition, $\e^\mu I'_n \e^{-\mu} \subset I_{n+1}$ and $\e^{-\mu} I''_n \e^\mu \subset I_{n+1}$. So $I''_n \e^\mu I'_n \subset I''_n I_{n+1} \e^\mu=I_{n+1} I''_n \e^\mu \subset I_{n+1} \e^\mu I_{n+1}$. This proves the surjectivity.  

On the other hand, for any $i_1 \in I'_n$ and $i_2 \in I''_n$, \begin{align*} & i_1 i_2 (I_n \cap M(L)) \cdot_\d I_{n+1} (I \cap M(L)) \e^\mu I_{n+1}=i_1 i_2 I_{n+1} (I \cap M(L)) \e^\mu I_{n+1} \d(i_1 i_2) \i  \\ & \subset i_1 I_{n+1} (I \cap M(L)) I''_n \e^\mu I'_n I_{n+1} \d(i_2) \i=i_1 I_{n+1} (I \cap M(L)) \e^\mu I_{n+1} \d(i_2) \i.\end{align*} So $i_1 i_2 (I_n \cap M(L)) \cdot_\d I_{n+1} (I \cap M(L)) \e^\mu I_{n+1} \cap I_{n+1} (I \cap M(L)) \e^\mu I_{n+1} \neq \emptyset$ if and only if $i_1, i_2 \in I_{n+1}$. This proves the injectivity. 

\subsection{Proof of Proposition \ref{00}}\label{pf00} Since $\tw$ is a good element, we have that $(I \dot \tw I) \d(I \dot \tw I) \cdots \d^{n-1}(I \dot \tw I) \subset I \e^\mu I$. By Lemma \ref{finite}, for any $b \in G(L) \cdot_\d I \dot \tw I$, there exists $h \in G(L)$ such that $$h \i b \d(b) \cdots \d^{n-1}(b) \d^n(h) \in (I \cap M(L)) \e^\mu.$$ If $b=\dot \tw$, then we may just take $h=1$. Now set $b'=h \i b \d(b) \cdots \d^{n-1}(b) \d^n(h)$. If $g I \in X_{\tw, \d}(b)$, then again by Lemma \ref{finite}, $g I=g' I$ for some $g'$ with $(g') \i b \d(b) \cdots \d^{n-1}(b) \d^n(g') \in (I \cap M(L)) \e^\mu$. Thus $$X_{\tw, \d}(b) \subset \{ h g I/I; g \i b' \d^n(g) \in (I \cap M(L)) \e^\mu\}.$$ 

Let $x \in W$ such that $x \mu$ is dominant. Set $J=I(t^{x \mu})$ and $M'={}^{\dot x} M$. Then $M'$ is a Levi factor of the standard parabolic subgroup $P$ of $G$ corresponding to $J$ and ${}^{\dot x} (I \cap M(L)) \subset \cap M'(o)$. It is easy to see that $$X_{\tw, \d}(b) \subset \{ h \dot x g \dot x \i I/I; g \i M'(o) \e^{x \mu} \d^n(g) \cap M'(o) \e^{x \mu} \neq \emptyset\}.$$

Let $g \in G(L)$ with $g \i M'(o) \e^{x \mu} \d^n(g) \cap M'(o) \e^{x \mu} \neq \emptyset$. Then we must have that $g=m k$ for $m \in M'(L)$ and $k \in G(o)$ with $m \i M'(o) \e^{x \mu} \d^n(m) \in M'(o) \e^{x \mu}$. The case where $\d=\d_F$ is proved in \cite[Theorem 1.1 (2)]{Ko2}. The case where $\d=\d_a$ can be proved in the same way. 

Assume that $k \in I \dot u I$ for $u=u_1 v_1$ with $u_1 \in W^J$ and $v_1 \in W_J$. Then $I \dot u I=I \dot u_1 I \dot v_1 I \subset I \dot u_1 I M'$. Thus $I \dot u_1 I M'(o) \e^{x \mu} \cap M'(o) \e^{x \mu} \d^n(I \dot u_1 I) \neq \emptyset$. Notice that \begin{gather*} M'(o) \e^{x \mu} \d^n(I \dot u_1 I) \subset \cup_{v \in W_J} I \dot v I \e^{x \mu} I \dot u_1 I=\cup_{v \in W_J} I \e^{x \mu} \dot v I \dot u_1 I \subset \cup_{y \in W} I \e^{x \mu} \dot y I, \\ I \dot u_1 I M'(o) \e^{x \mu} \subset \cup_{v \in W_J} I \dot u_1 I \dot v I \e^{x \mu}=\cup_{v \in W_J} I \dot u_1 I \e^{x \mu} I \dot v I=\cup_{v \in W_J} I \dot u_1 \e^{x \mu} \dot v I. \end{gather*} Hence there exists $v \in W_J$ and $y \in W$ such that $\e^{x \mu} y=u_1 \e^{x \mu} v=\e^{u_1 x \mu} u_1 v$. In particular, we have that $x \mu=u_1 x \mu$. By \cite[Lemma 3.5]{HT}, $u_1 \in W_J$. So $u_1=1$ and $k \in I M'(o)$. 

Now $k \i m \e^{x \mu} \d^n(k)=m' \e^{x \mu}$ for some $m, m' \in M'(o)$. By \cite[Lemma 7.3]{St68} for $\d=\d_a$ and Lang's theorem for $\d=\d_F$, $m' \in k_1 (I \cap M'(o)) \d^n(k_1) \i$ for some $k_1 \in M'(o)$. Set $k'=k k_1$. Then $(k') \i (M'(o) \e^{x \mu}) \d^n(k') \cap (I \cap M'(o)) \e^{x \mu} \neq \emptyset$. We assume that $k' \in M'(o) k''$ for some $k'' \in I$. Then $(k') \i (M'(o) \e^{x \mu}) \d^n(k')=(k'') \i (M'(o) \e^{x \mu}) \d^n(k'')$ and $k'' (I \cap M'(o)) \e^{x \mu} \d^n(k'') \i \cap (M'(o) \e^{x \mu}) \neq \emptyset$. 

By Lemma \ref{finite}, $k'' (I \cap M'(o)) \e^{x \mu} \d^n(k'') \i \cap (M'(o) \e^{x \mu}) \subset (I \cap M'(o)) \e^{x \mu}$ and $k'' \in I \cap M'(o)$. So $k \in M'(o)$. Therefore $X_{w, \d}(b) \subset \{h \dot x g \dot x \i I/I; g \in M'(L)\}=\{h g I/I; g \in M(L)\}$. The ``moreover'' part follows from Lemma \ref{0}. 

\begin{cor}\label{t0}
Let $\mu \in P^\vee$ be a regular coweight. Then $G(L) \cdot_\d I \e^\mu I$ is a single $\d$-conjugacy class of $G(L)$ and $\dim(X_{t^\mu, \d}(\e^\mu))=0$. 
\end{cor}

By Lemma \ref{finite}, $I \e^\mu I$ is a single $\d$-conjugacy class of $I$. Hence $G(L) \cdot_\d I \e^\mu I$ is a single $\d$-conjugacy class of $G(L)$. Now by Proposition \ref{00}, $X_{t^\mu, \d}(\e^\mu) \subset T(L) I/I \cong T(L)/(T(L) \cap I)$ is of dimension $0$.

\section{Applications on loop groups of type A and C}

We first discuss some stratification of $G(L)$, stable under the $\d$-conjugation action. 

\begin{lem}\label{eq1}
Let $G=GL_n$ or $PSP_{2 n}$. Let $\aa \in B(\tW_{int})$. Then 

(1) For any good elements $\tw, \tw' \in f \i(\aa)$, we have that $$G(L) \cdot_\d I \dot \tw' I=G(L) \cdot_\d I \dot \tw I=G' \cdot_\d I \dot \tw I.$$  Now we define $Z_{\aa, \d}=G(L) \cdot_\d I \dot \tw I$ for any good element $\tw \in f \i(\aa)$. 

(2) Let $\co \subset f \i(\aa)$ be a conjugacy class of $\tW$ and $\tw \in \co_{\min}$. Then $G(L) \cdot_\d I \dot \tw I \subset Z_{\aa, \d}$. 
\end{lem}

(1) By Corollary \ref{good3}, $\tw \approx \tw'$. Thus by Lemma \ref{red2}, $G(L) \cdot_\d I \dot \tw' I=G(L) \cdot_\d I \dot \tw I$. Notice that $G(L)=\sqcup_{\t \in \tW, l(\t)=0} G' \dot \t$. Thus \begin{align*} G(L) \cdot_\d I \dot \tw I &=\cup_{\t \in \tW, l(\t)=0} G' \cdot_\d \dot \t I \dot \tw I \d(\dot \t) \i \\ &=\cup_{\t \in \tW, l(\t)=0} G' \cdot_\d I \dot \t \tw \dot \t \i I.\end{align*} Since $l(\t \tw \t \i)=l(\tw)$, $\t \tw \t \i$ is also a good element in $f \i(\aa)$.  By Corollary \ref{good3} and Lemma \ref{red11}, $G' \cdot_\d I \dot \tw I=G' \cdot_\d I \dot \t \dot \tw \dot \t \i I$. Hence $G(L) \cdot_\d I \dot \tw I=G' \cdot_\d I \dot \tw I$. 

(2) If $G=GL_n$, then $f \i(\aa)=\sqcup_{d(\tilde \l', \tilde \O)=(\tilde \l, \tilde \O)} [(\tilde \l', \tilde \O)]$ for some distinguished $(\tilde \l, \tilde \O) \in \cd \cp$. By Theorem \ref{62}, $\tw \tilde \approx x \tw^f_{(\tilde \l, \tilde \O)}$ for some $x \in W_{I(\tw^f_{(\tilde \l, \tilde \O)})}$. By Lemma \ref{red2} and \ref{red3}, $$G(L) \cdot_\d I \dot \tw I=G(L) \cdot_\d I \dot \tx \dot \tw^f_{(\tilde \l, \tilde \O)} I \subset G(L) \cdot_\d I \dot \tw^f_{(\tilde \l, \tilde \O)} I.$$

If $G=PGL_{2n}$, then $f \i(\aa)=\sqcup_{d'(\tilde \l', \tilde \mu)=(\tilde \l, \tilde \O)} [(\tilde \l', \tilde \mu)]'$ for some distinguished $(\tilde \l, \tilde \O) \in \cd \cp_{\ge 0}$. Hence $\co=[(\tilde \l', \tilde \mu)]'$ for some $(\tilde \l', \tilde \mu) \in \cd \cp_{\ge 0}$ with $d'(\tilde \l', \tilde \mu)=(\tilde \l, \tilde \O)$. Assume that $d(\tilde \l', \tilde \mu)=(\tilde \l'', \tilde \mu')$. Then $d(\tilde \l'', \underline{\tilde \mu'})=(\tilde \l, \tilde \O)$. By Theorem \ref{62}, $\tw \tilde \approx x \tw^f_{(\tilde \l'', \tilde \mu')}$ for some $x \in W_{I(\tw^f_{(\tilde \l'', \tilde \mu')})}$. By Lemma \ref{red2} and \ref{red3}, $$G(L) \cdot_\d I \dot \tw I=G(L) \cdot_\d I \dot x \dot \tw^f_{(\tilde \l'', \tilde \mu')} I \subset G(L) \cdot_\d I \dot \tw^f_{(\tilde \l'', \tilde \mu')} I.$$ Again by Theorem \ref{62}, $\tw^f_{(\tilde \l'', \tilde \mu')} \tilde \approx y \tw^f_{(\tilde \l, \tilde \O)}$ for some $y \in W_{I(\tw^f_{(\tilde \l, \tilde \O)})}$. Then by Lemma \ref{red2} and \ref{red3}, $$G(L) \cdot_\d I \dot \tw I \subset G(L) \cdot_\d I \dot \tw^f_{(\tilde \l'', \tilde \mu')} I \subset G(L) \cdot_\d I \dot \tw^f_{(\tilde \l, \tilde \O)} I.$$  

\begin{prop}\label{par1}
Let $G=GL_n$ or $PSP_{2 n}$ and $\tw \in \tW_{int}$. Then $$\overline{G(L) \cdot_\d I \dot \tw I}=G(L) \cdot_\d \overline{I \dot \tw I}=\sqcup_{\aa \preceq \tw} Z_{\aa, \d}.$$
\end{prop}

By the proof of \cite[Prop 18]{V}, $\overline{G(L) \cdot_\d I \dot \tw I}=G(L) \cdot_\d \overline{I \dot \tw I}$. By Lemma \ref{dg}, the union is in fact a disjoint union. If $\aa \preceq \tw$, then there exists a minimal length element $\tw'$ in $f \i(\aa)$ such that $\tw' \le \tw$. Hence $Z_{\aa, \d}=G(L) \cdot_\d I \dot \tw' I \subset G(L) \cdot_\d \overline{I \dot \tw I}$. Now we prove that $G(L) \cdot_\d I \dot \tw I \subset \cup_{\aa \preceq \tw} Z_{\aa, \d}$ by induction on $l(\tw)$. 

If $\tw \in f \i(\aa)$ is a minimal length element in its conjugacy class, then by Corollary \ref{xx}, $\aa \preceq \tw$. By Lemma \ref{eq1}, $G(L) \cdot_\d I \dot \tw I \subset Z_{\aa, \d}$. 

If $\tw$ is not a minimal length element in the conjugacy class of $\tW$ that contains it, then by Theorem \ref{thA} there exists $\tw_1 \approx \tw$ and $\tw_1 \xrightarrow i \tw_2$ with $l(\tw_2)<l(\tw_1)$. In particular, $\tw_2<\tw_1$ and $s_i \tw_1<\tw_1$. By induction hypothesis, \begin{align*} G(L) \cdot_\d I \dot \tw I &=G(L) \cdot_\d I \dot \tw_1 I=G(L) \cdot_\d I \dot \tw_2 I \cup G(L) \cdot_\d I \dot s_i \dot \tw_1 I  \\ & \subset \cup_{\aa \preceq \tw_1} Z_{\aa, \d}.\end{align*} By $\S$ \ref{porr}, $\aa \preceq \tw$ if and only if $\aa \preceq \tw_1$. So $G(L) \cdot_\d I \dot \tw I \subset \cup_{\aa \preceq \tw} Z_{\aa, \d}$. 

\begin{cor}
Let $G=GL_n$ or $PSP_{2 n}$ and $\aa \in B(\tW_{int})$. Then $\overline{Z_{\aa, \d}}=\sqcup_{\aa' \preceq \aa} Z_{\aa', \s}$. 
\end{cor}

\begin{rmk}
If $\d=\d_F$ is a Frobenius morphism, then $Z_{\aa, \d}$ is a single $\d$-conjugacy class and any $\d$-conjugacy class is of this form. In this case, the closure of $Z_{\aa, \d}$ is a union of other $\d$-conjugacy classes and the explicit closure relation is obtained by Viehmann \cite[Prop 18]{V}. However, if $\d=\d_a$, then in general $Z_{\aa, \d}$ contains infinitely many $\d$-conjugacy classes. 
\end{rmk}

\begin{prop}\label{str}
(1) Let $G=GL_n$, then $G(L)=\sqcup_{\aa \in B(\tW)} Z_{\aa, \d}$ is a stratification of $G(L)$. 

(2) Let $G=PSP_{2n}$, then $G'=\sqcup_{\aa \in B(\tW_{int})} Z_{\aa, \d}$ is a stratification of $G'$. 
\end{prop}

\begin{rmk}
If $\d=\d_F$, then both parts follows from Kottwitz's classification of $\s$-conjugacy classes \cite{Ko97}. See subsection \ref{reK}. 
\end{rmk}

Notice that \begin{gather*} G(L)=\sqcup_{\tw \in \tW} I \dot \tw I=\cup_{\tw \in \tW} G(L) \cdot_\d I \dot \tw I=\cup_{\tw \in \tW} G(L) \cdot_\d \overline{I \dot \tw I}; \\ G'=\sqcup_{\tw \in W_a} I \dot \tw I=\cup_{\tw \in W_a} G(L) \cdot_\d I \dot \tw I=\cup_{\tw \in W_a} G(L) \cdot_\d \overline{I \dot \tw I}.\end{gather*} Now the Proposition follows from Proposition \ref{par1}.  

\

Now we discuss the dimension of $X_{\tw, \d}(b)$ for type A and C.

\begin{prop}\label{class3}
Let $G=GL_n$ or $PSP_{2n}$. Let $\aa \in B(\tW_{int})$. Then for any $\tw \in \tW$ and $b \in Z_{\aa, \d}$, we have that 

(1) $\dim (X_{\tw, \d_F}(b))=\max_{\co} \frac{1}{2}(l(\tw)+l(\co)+\deg(f_{\tw, \co}))-l(f \i(\aa))$, here $\co$ runs over conjugacy class of $\tW$ in $f \i(\aa)$.

(2) $\dim(X_{\tw, \d_a}(b))< \infty$ if $a$ is not a root of unity.
\end{prop}

Since $b \in Z_{\aa, \d}$ for $\aa \in B(\tW_{int})$, we have that $b \in G'$. Hence if $X_{\tw, \d}(b) \neq \emptyset$, then $\tw \in \tW_{int}$. Now let $\co$ be an integral conjugacy class of $\tW$ and $\tw_\co$ a minimal length element in $\co$. If $X_{\tw_\co, \d}(b) \neq \emptyset$, then $b \in G(L) \cdot_\d I \dot \tw_\co I$. By Lemma \ref{eq1} and Lemma \ref{dg}, we must have that 

(a) $\co \subset f \i(\aa)$. 

Case I: $G=GL_n$. Then $f \i(\aa)=\sqcup_{d(\tilde \l', \tilde \O)=(\tilde \l, \tilde \O)} [(\tilde \l', \tilde \O)]$ for some distinguished $(\tilde \l, \tilde \O) \in \cd \cp$. By Theorem \ref{62}, $\tw_\co \tilde \approx x \tw^f_{(\tilde \l, \tilde \O)}$ for some $x \in W_{I(\tw^f_{(\tilde \l, \tilde \O)})}$. By Corollary \ref{tttt}, $\dim(X_{\tw_\co, \d}(b))=\dim(X_{x \tw^f_{(\tilde \l, \tilde \O)}, \d}(b))$. 

If $\d=\d_F$, then by Lemma \ref{lang} and Prop \ref{00}, $\dim(X_{x  \tw^f_{(\tilde \l, \tilde \O)}, \d}(b))=\dim(X_{\tw^f_{(\tilde \l, \tilde \O)}, \d}(b))+l(x)=l(x)=l(\tw_\co)-l(\tw^f_{(\tilde \l, \tilde \O)})=l(\co)-l(f \i(\aa))$. Now by Proposition \ref{class}, $\dim (X_{\tw, \d}(b))=\max_{\co} \frac{1}{2}(l(\tw)+l(\co)+\deg(f_{\tw, \co}))-l(f \i(\aa))$. 

If $\d=\d_a$ for some $a$ not a root of unity, then by Lemma \ref{lang} and Prop \ref{00}, $\dim(X_{x \tw^f_{(\tilde \l, \tilde \O)}}(b)) \le \dim(X_{\tw^f_{(\tilde \l, \tilde \O)}, \d}(b))+l(w_S)<\infty$. Hence by Proposition \ref{class2}, $\dim(X_{\tw, \d}(b))<\infty$. 

Case II: $G=PGL_{2n}$. Then $f \i(\aa)=\sqcup_{d'(\tilde \l', \tilde \mu)=(\tilde \l, \tilde \O)} [(\tilde \l', \tilde \mu)]'$ for some distinguished $(\tilde \l, \tilde \O) \in \cd \cp_{\ge 0}$ and $\co=[(\tilde \l, \tilde \mu)]'$ for some $d'(\tilde \l', \tilde \mu)=(\tilde \l, \tilde \O)$. Let $d(\tilde \l', \tilde \mu)=(\tilde \l'', \tilde \mu')$. Then by Theorem \ref{62}, $\tw_\co \tilde \approx x \tw^f_{(\tilde \l'', \tilde \mu')}$ for some $x \in W_{I(\tw^f_{(\tilde \l'', \tilde \mu')})}$ and $\tw^f_{(\tilde \l'', \tilde \mu')} \tilde \approx y \tw^f_{(\tilde \l, \tilde \O)}$ for some $y \in W_{I(\tw^f_{(\tilde \l, \tilde \O)})}$. 

If $\d=\d_F$, then by Corollary \ref{tttt}, Lemma \ref{lang} and Prop \ref{00}, \begin{align*} & \dim(X_{\tw_\co, \d}(b))=\dim(X_{x \tw^f_{(\tilde \l'', \tilde \mu')}, \d}(b))=\dim(X_{\tw^f_{(\tilde \l'', \tilde \mu')}, \d}(b))+l(x) \\ &=\dim(X_{y \tw^f_{(\tilde \l, \tilde \O)}, \d}(b))+l(x)=\dim(X_{\tw^f_{(\tilde \l, \tilde \O)}, \d}(b))+l(x)+l(y). \end{align*} Therefore $\dim(X_{\tw_\co, \d}(b))=l(\co)-l(f \i(\aa))$. Now by Proposition \ref{class}, $\dim (X_{\tw, \d}(b))=\max_{\co} \frac{1}{2}(l(\tw)+l(\co)+\deg(f_{\tw, \co}))-l(f \i(\aa))$. 

If $\d=\d_a$ for some $a$ not a root of unity, then by Corollary \ref{tttt}, Lemma \ref{lang} and Prop \ref{00}, \begin{align*} & \dim(X_{\tw_\co, \d}(b))=\dim(X_{x \tw^f_{(\tilde \l'', \tilde \mu')}, \d}(b)) \le \dim(X_{\tw^f_{(\tilde \l'', \tilde \mu')}, \d}(b))+l(w_S) \\ &=\dim(X_{y \tw^f_{(\tilde \l, \tilde \O)}, \d}(b))+l(w_S) \le \dim(X_{\tw^f_{(\tilde \l, \tilde \O)}, \d}(b))+2 l(w_S)<\infty. \end{align*} Hence by Proposition \ref{class2}, $\dim(X_{w, \d}(b))<\infty$. 

\begin{prop}
We assume that either 

(1) $G=GL_n$ and $\mu$ is a regular coweight or 

(2) $G=PSP_{2n}$ and $\mu$ is a regular coweight that lies in the coroot lattice. 

Then for any $\tw \in \tW_{int}$ and $a \in \kk^\times$, $$\dim(X_{\tw, \d_a}(\e^\mu))=\frac{1}{2}(l(\tw)+\deg(f_{\tw, \co})-l(t^\mu)),$$ where $\co$ is the conjugacy class of $\tW$ that contains $t^\mu$. If moreover, $\kk$ is of positive characteristic, then $\dim(X_{\tw, \d_a}(\e^\mu))=\dim(X_{\tw, \d_F}(\e^\mu))$. 
\end{prop}

We prove the Proposition by induction on $l(\tw)$. Let $\co'$ be the integral conjugacy class that contains $\tw$. 

If $\co=\co'$. Then $l(\tw)=l(t^\mu)$ is a minimal length element in $\co$. By definition, $f_{\tw, \co}=1$. By Corollary \ref{t0}, $\dim(X_{\tw, \d}(\e^\mu))=0$. The Proposition holds in this case. 

If $\co \neq \co'$ and $\tw \in \co'_{\min}$. Then by definition, $f_{\tw, \co}=0$. Since $\mu$ is regular, $\co$ is a single fiber of $f: \tW \to B(\tW)$. By (a) in the proof of Proposition \ref{class3}, $X_{\tw_{\co'}, \d}(\e^\mu)=\emptyset$. The Proposition holds in this case.

If $\tw \notin \co'_{\min}$, the Proposition follows from induction on $l(\tw)$ by the same argument as we didi in the proof of Proposition \ref{class}. 

\subsection{} Let us come to the case where $\d=\d_1$ is the identity map. In this case, the $\d$-conjugacy classes in $G(L)$ are just the usual conjugacy classes. Let $b$ be a regular semisimple, integral element in $G(L)$. Here integral means the elements in $G(L) \cdot I$. It is shown by Kazhdan and Lusztig in \cite{KL} that $X_{(1, id)}(b)$ is finite dimensional and a conjectural dimension formula is also given there. If moreover, $b$ is elliptic (i.e., its centralizer is an anisotropic torus), then $X_{1, id}(b)$ has only finitely many irreducible components and is an algebraic variety. The conjecture is proved later by Bezrukavnikov in \cite{Be}. (Actually they considered only topologically unipotent elements, but the general can be reduced to that case using Jordan decomposition). Now by the same argument as we did in the proof of Proposition \ref{class3}, one can show the following result. This answers the question in \cite[Section 7]{L3} for $G=GL_n$ or $PSP_{2n}$. 

\begin{prop}\label{class5}
Let $G=GL_n$ or $PSP_{2n}$. Let $b$ be a regular semisimple integral element in $G(L)$. Then for any $\tw \in \tW$, $X_{\tw, id}(b)$ is finite dimensional. If moreover, $b$ is elliptic, then $X_{\tw, id}(b)$ is an algebraic variety. 
\end{prop}

\subsection{} From now on, we only consider the case that $\d=\d_F$. Similar results holds for $\d=\d_a$ with $a$ not a root of unity. Define $$J_{b, \d_F}=\{g \in G(L); g \i b \d_F(g)=b\}.$$ Then $J_{b, \d_F}$ acts on $X_{\tw, \d_F}(b)$ on the left for any $\tw \in \tW$. 

Notice that $X_{\tw, \d_F}(b)=\lim\limits_{\rightarrow} X_i$ for some closed subschemes $X_1 \subset X_2 \subset \cdots \subset X_{\tw, \d_F}(b)$ of finite type. Let $l$ be a prime with $l$ not equal to the characteristic of $\kk$. Then $H^j_c(X_i, \bar \Q_l)$ is defined for all $j$. Set $$H^{BM}_j(X_{\tw, \d_F}(b), \bar \Q_l)=\lim\limits_{\rightarrow} H^j_c(X_i, \bar \Q_l)^*.$$ Then $H^{BM}_j(X_{\tw, \d_F}(b), \bar \Q_l)$ is a smooth representation of $J_{b, \d_F}$. Hence it is a semisimple module for any open compact subgroup of $J_{b, \d_F}$. 

The following result can be proved along the line of \cite[Theorem 1.6]{DL}. 

\begin{lem}\label{inddd}
Let $b \in G(L)$ and $K$ be an open compact subgroup of $J_{b, \d_F}$. Let $\tw \in \tW$, and let $s\in\tS$ be a simple affine reflection. Then 

(1) If $l(s \tw s)=l(\tw)$, then for any $j \in \ZZ$, $$H^{BM}_j(X_{\tw, \d_F}(b), \bar \Q_l) \cong H^{BM}_j(X_{s \tw s, \d_F}(b), \bar \Q_l)$$ as $J_{b, \d_F}$-modules. 

(2) If $l(s \tw s) =l(\tw)-2$, then for any simple $K$-module $M$ that is a direct summand of $\oplus_j H^{BM}_j(X_{\tw, \d_F}(b), \bar \Q_l)$, $M$ is also a direct summand of $\oplus_j H^{BM}_j(X_{s \tw s, \d_F}(b), \bar \Q_l) \oplus \oplus_j H^{BM}_j(X_{s \tw, \d_F}(b), \bar \Q_l)$. 
\end{lem}

\begin{thm}\label{BM}
Let $G=GL_n$ or $PSP_{2n}$, $\aa \in B(\tW_{int})$ and $b \in Z_{\aa, \d_F}$. Let $K$ be an open compact subgroup of $J_{b, \d_F}$ and $M$ be a simple $K$-module. If $M$ is a direct summand of $\oplus_{\tw \in \tW} \oplus_j H^{BM}_j(X_{\tw, \d_F}(b), \bar \Q_l)$, then $M$ is a direct summand of $\oplus_{\tx \in \co_{\min}} \oplus_j H^{BM}_j(X_{\tx, \d_F}(b), \bar \Q_l)$, where $\co$ runs over conjugacy classes of $\tW$ on $f \i(\aa)$. 
\end{thm}

Notice that $X_{\tw, \d_F}(b) \neq \emptyset$ implies that $\tw \in \tW_{int}$. By Theorem \ref{thA} and Lemma \ref{inddd}, $M$ is a direct summand of $\oplus_j H^{BM}_j(X_{\tx, \d_F}(b), \bar \Q_l)$ for a minimal length element $\tx$ in some integral conjugacy class $\co$ of $\tW$. Then $X_{\tx, \d_F}(b) \neq \emptyset$. By Lemma \ref{eq1} and Lemma \ref{dg}, we must have that $\co \subset f \i(\aa)$. The theorem is proved. 





\subsection{} In the rest of this paper, we discuss in more details the special cases that $f \i(\aa)$ is a single conjugacy class. 

Let $\cd \cp_{>0}$ be the set of $(\tilde \l, \tilde \mu) \in \cd \cp$ such that all the entries of $\tilde \l$ and $\tilde \mu$ are of the form $(b, c)$ with $c>0$. We call $(\tilde \l, \tilde \O) \in \cd \cp$ {\it super distinguished} if it is distinguished and all the entries of $\tilde \l$ are distinct. Then it is easy to see that

(1) Let $(\tilde \l, \tilde \O) \in \cd \cp$. Then $\sharp\{(\tilde \l', \tilde \O) \in \cd \cp; d(\tilde \l', \tilde \O)=(\tilde \l, \tilde \O)\}=1$ if and only if $(\tilde \l, \tilde \O)$ is super distinguished. 

(2) Let $(\tilde \l, \tilde \O) \in \cd \cp_{\ge 0}$. Then $\sharp\{(\tilde \l', \tilde \mu) \in \cd \cp_{\ge 0}; d'(\tilde \l', \tilde \mu)=(\tilde \l, \tilde \O)\}=1$ if and only if $(\tilde \l, \tilde \O)$ is super distinguished and $(\tilde \l, \tilde \O) \in \cd \cp_{>0}$. 

In other words,

Let $G=GL_n$ and $\aa \in B(\tW)$. Then $f\i(\aa)$ is a single conjugacy class of $\tW$ if and only if $f \i(\aa)=[(\tilde \l, \tilde \O)]$ for some super distinguished $(\tilde \l, \tilde \O) \in \cd \cp$. 

Let $G=PSP_{2n}$ and $\aa \in B(\tW_{int})$. Then $f\i(\aa)$ is a single conjugacy class of $\tW$ if and only if $f \i(\aa)=[(\tilde \l, \tilde \O)]'$ for some super distinguished $(\tilde \l, \tilde \O) \in \cd \cp_{>0}$. 

In either case, for any $\tw, \tw' \in f \i(\aa)_{\min}$, $\tw \approx \tw'$ by Theorem \ref{62}. Hence by Lemma \ref{inddd}, $X_{\tw, \d_F}(b)$ is universally isomorphic to $X_{\tw', \d_F}(b)$. Hence by Theorem \ref{BM}, we have that 

\begin{cor}\label{unique}
Let $G=GL_n$ or $PSP_{2n}$, $\aa \in B(\tW_{int})$ and $b \in X_{\aa, \d_F}$. Let $K$ be an open compact subgroup of $J_{b, \d_F}$ and $M$ be a simple $K$-module. We assume furthermore that $f \i(\aa)$ is a single conjugacy class of $\tW$. If $M$ is a direct summand of $\oplus_{\tw \in \tW} \oplus_j H^{BM}_j(X_{\tw, \d_F}(b), \bar \Q_l)$, then $M$ is a direct summand of $\oplus_j H^{BM}_j(X_{\tx, \d_F}(b), \bar \Q_l)$ for any $\tx \in f \i (\aa)_{\min}$. 
\end{cor}

\subsection{} We mention two interesting cases of super distinguished pair of double partitions $(\tilde \l, \tilde \O)$. 

The first case is that all the entries of $\tilde \l$ are of the form $(1, *)$. In this case, $\tw^f_{(\tilde \l, \tilde \O)}=t^\chi$, where $\chi$ is a dominant regular coweight for type A and a dominant regular coweight that lies in the coroot lattice for type C. 

The second case is that $\tW$ is of type A and $\tilde \l=(n, r)$ for some $0<r<n$ such that $n$ and $r$ are coprime. By Lemma \ref{asdf}, $\tw^f_{(\tilde \l, \tilde \O)}=t^{\o_r} w_{{S-\{r\}}} w_S$. This is a superbasic element. 

\begin{cor}\label{ddd}
We assume that either 

(1) $G=GL_n$ and $\chi$ is a dominant regular coweight or 

(2) $G=PSP_{2n}$ and $\chi$ is a dominant regular coweight that lies in the coroot lattice. 

Then $J_{\e^\chi, \d_F}=T(L)^{\d_F}$ and $T(o)^{\d_F}$ acts trivially on $H^{BM}_j(X_{\tw, \d_F}(\e^\chi), \bar \Q_l)$ for any $\tw \in \tW$ and $j \in \ZZ$. 
\end{cor}

\begin{rmk}
The special case for $G=GL_2$ or $GL_3$ is obtained by Zbarsky \cite{Zb} in a different way. 
\end{rmk}

By Proposition \ref{00}, $J_{\e^\chi, \d_F} \subset T(L) I$. Since $\chi$ is dominant regular, it is easy to see that $J_{\e^\chi, \d_F} \subset T(L)$. Hence $J_{\e^\chi, \d_F}=T(L)^{\d_F}$ and $J_{\e^\chi, \d_F} \cap I=T(o)^{\d_F}$ is an open compact subgroup of $J_{\e^\chi, \d_F}$. Since $T(L) \cap I$ acts trivially on $T(L)/(T(L) \cap I)$, $J_{\e^\chi, \d_F} \cap I$ acts trivially on $X_{t^\chi, \s}(\e^\chi)$ and also trivially on $H^{BM}_j(X_{t^\chi, \s}(\e^\chi), \bar \Q_l)$. By Corollary \ref{unique}, any simple module of $J_{\e^\chi, \d_F} \cap I$ that appears as a direct summand of $H^{BM}_j(X_{\tw, \d_F}(\e^\chi), \bar \Q_l)$ is also trivial. Hence $J_{\e^\chi, \d_F} \cap I$ acts trivially on $H^{BM}_j(X_{\tw, \d_F}(\e^\chi), \bar \Q_l)$.

\begin{cor}\label{superbasic}
Let $G=GL_n$ and $\t=t^{\o_r} w_{{S-\{r\}}} w_S$ for some $0<r<n$ such that $n$ and $r$ are coprime. Let $\Om=\{\tx \in \tW; l(\tx)=0\}$. Then $J_{\dot \t, \d_F}/(J_{\dot \t, \d_F} \cap I) \cong \Om$ and $J_{\dot \t, \d_F} \cap I$ acts trivially on $H^{BM}_j(X_{\tw, \d_F}(\dot \t), \bar \Q_l)$ for any $\tw \in \tW$ and $j \in \ZZ$. 
\end{cor}

Let $g \in G(L)$ with $g \i \dot \t \d_F(g) \in I \dot \t$. We may assume that $g \in I \dot \tx I$ for some $\tx \in \tW$. Now $g \in \dot \t \d_F(g) \dot \t \i I=I \dot \t \dot \tx \dot \t \i I$ and $I \dot \tx I \cap I \dot \t \i \dot \tx \dot \t I \neq \emptyset$. Therefore $\tx=\t \tx \t \i$. If $l(\tx) \neq 0$, then there exist $i \in \tS$ such that $s_i \tx<\tx$. Since conjugation by $\t$ preserve the Bruhat order, we have that $s_{\t^m(i)} \tx=s_{\t^m(i)} \t^m \tx \t^{-m} \le \t^m \tx \t^{-m}=\tx$ for any $m \in \ZZ$. However, $\tS$ is a single $\t$-orbit. Therefore $s_j \tx<\tx$ for any $j \in \tS$. That is impossible. Hence $l(\tx)=0$ and $g \in I \dot \tx$. On the other hand, for any $\tx \in \Om$, $\d_F(\dot \tx)=\dot \tx t$ for some $t \in T$. Then $(\dot \tx t') \i \dot \t \d_F(\dot \tx t')=(t') \i \dot \tx \i \dot \t  \dot \tx t \d_F(t')=(t') \i \dot \t t \d_F(t')$ for any $t' \in T$. By Lang's theorem, there exists $t' \in T$ such that $(t') \i \dot \t t \d_F(t')=\dot \t$ and $\dot \tx t' \in J_{\dot \t, \d_F}$. 

Therefore, $J_{\dot \t, \d_F}/(J_{\dot \t, \d_F} \cap I) \cong \Om$ and $X_{\t, \d_F} \cong \Om$. In particular, $J_{\dot \t, \d_F} \cap I$ acts trivially on $X_{\t, \d_F}(\dot \t)$. Hence $J_{\dot \t, \d_F} \cap I$ acts trivially on $H^{BM}_j(X_{\t, \d_F}(\dot \t), \bar \Q_l)$. 

By Corollary \ref{unique}, any simple module of $J_{\dot \t, \d_F} \cap I$ that appears as a direct summand of $H^{BM}_j(X_{\tw, \d_F}(\dot \t), \bar \Q_l)$ is also trivial. Hence $J_{\dot \t, \d_F} \cap I$ acts trivially on $H^{BM}_j(X_{\tw, \d_F}(\dot \t), \bar \Q_l)$.

\section*{Acknowledgement} 
We thank Victor Ginzburg, Ulrich G\"ortz, George Lusztig, Michael Rapoport, Eva Viehmann and Yong-Chang Zhu for some helpful discussions. Some application on affine Deligne-Lusztig varieties was done during the visit at University of Bonn. We thank them for the warm hospitality.

\end{document}